\documentclass[11pt]{amsart}

\usepackage[T1]{fontenc}
\usepackage{amsmath,amssymb,amsthm,mathrsfs}
\usepackage{mathtools}

\usepackage{microtype}
\usepackage{needspace}
\usepackage{xcolor}
\usepackage{enumitem}
\usepackage[colorlinks=true,linkcolor=blue,
citecolor=blue,urlcolor=blue]{hyperref}
  \usepackage{parskip}

\numberwithin{equation}{section}
\allowdisplaybreaks

\hypersetup{
pdftitle={Dimension expansion for analytic maps and spectral averaging for Jacobi operators},
pdfauthor={Ji Li, Minh-Quy Pham, Thang Pham, Chun-Yen Shen}
}

\theoremstyle{plain}
\newtheorem{theorem}{Theorem}[section]
\newtheorem{lemma}[theorem]{Lemma}
\newtheorem{proposition}[theorem]{Proposition}
\newtheorem{corollary}[theorem]{Corollary}
\newtheorem{claim}{Claim}[section]
\theoremstyle{definition}
\newtheorem{definition}[theorem]{Definition}

\theoremstyle{remark}
\newtheorem{remark}[theorem]{Remark}

\newcommand{\R}{\mathbb{R}}
\newcommand{\Leb}{\mathcal{L}^{1}}
\newcommand{\supp}{\operatorname{supp}}
\newcommand{\reg}{\mathrm{reg}}
\newcommand{\e}{\varepsilon}
\newcommand{\dimH}{\dim_H}
\newcommand{\Gaux}{\mathcal G}
\newcommand{\dd}{\,d}

\title[Dimension expansion $\&$ spectral averaging for Jacobi operators]{Explicit dimension expansion for analytic maps and spectral averaging for Jacobi operators\\ with singular site laws}

\author{Ji Li}
\address{School of Mathematical and Physical Sciences,
Macquarie University, Sydney, Australia}
\email{ji.li@mq.edu.au}

\author{Minh-Quy Pham}
\address{Department of Mathematics, Virginia Tech, USA}
\email{mqpham@vt.edu}

\author{Thang Pham}
\address{Institute of Mathematics and Interdisciplinary Sciences,
Xidian University, China}
\email{thangphammath@xidian.edu.cn}

\author{Chun-Yen Shen}
\address{Department of Mathematics, National Taiwan University,
Taiwan}
\email{cyshen@math.ntu.edu.tw}

\date{November 23, 2026}

\begin{document}
\maketitle

\begin{abstract}
We prove explicit dimension expansion estimates for images of Cartesian
powers of a Borel set under real analytic maps. Let $k\ge2$, and let $f$ be real
analytic on a connected open neighborhood of $[0,1]^k$. Assume that
$f$ depends on each variable and admits no local representation
$f(x_1,\ldots,x_k)=g(h_1(x_1)+\cdots+h_k(x_k))$ by real analytic
functions of one variable. Given a Borel set $A\subset[0,1]$ of
Hausdorff dimension $\alpha\in (0,1]$, we show that for $k\ge3$,
\[
\dim_H f(A^k)\ge\min\{3\alpha/2,1\},
\]
and that $f(A^k)$ has positive Lebesgue measure whenever $\alpha>2/3$.
For $k=2$, we obtain a piecewise lower bound for $\dim_H f(A^2)$.
In particular, $\dim_H f(A^2)\ge 4\alpha/3$ for $0<\alpha\le 2/3$,
and $\dim_H f(A^2)=1$ for $\alpha\ge 3/4$. If $\alpha>3/4$,
then $f(A^2)$ has positive Lebesgue measure.

As an application to spectral theory, let $m\geq 4$ be even, and let $J$
be a bounded half-line Jacobi operator with fixed nonzero off-diagonal
entries. Suppose that the first $m$ diagonal entries are independent,
that the law of the $j$th entry is compactly supported and
$s_j$-Frostman for every $1\leq j\leq m$, and that the remaining
diagonal entries are fixed. We establish that if
$\sum_{j=1}^m s_j>m/2+1$, then for every $u\in\ell^2(\mathbb N)$
the averaged scalar spectral measure
\[
B\longmapsto \mathbb E\bigl[\langle \mathbf{1}_B(J)u,u\rangle\bigr]
\]
is absolutely continuous with respect to Lebesgue measure on
$\mathbb R$, even when all these laws are singular.

Furthermore, let $H$ be an $n\times n$ Jacobi matrix, $n\geq 3$, with
fixed nonzero off-diagonal entries and independent diagonal entries
having compactly supported nonatomic laws. Suppose that the laws at
three consecutive sites are $s_1$-, $s_2$-, and $s_3$-Frostman,
respectively, with $s_1+s_2+s_3>2$. Writing
$\lambda_1<\cdots<\lambda_n$ for the ordered eigenvalues of $H$, we
show that the law of each $\lambda_i$ and the law of each gap
$\lambda_j-\lambda_i$ with $1\leq i<j\leq n$ are absolutely continuous
with respect to Lebesgue measure.
\end{abstract}

{\small

}
\tableofcontents
    \section{Introduction}\label{sec:introduction}
{
We study the size and regularity of nonlinear images of Cartesian
products. Our geometric results give explicit Hausdorff dimension
bounds for real analytic maps and positive Lebesgue measure above
suitable dimension thresholds. The Fourier estimates also yield
absolute continuity results for images of product measures, which
we apply to eigenvalue distributions and spectral averaging for
Jacobi operators with independent singular site laws.

Theorem~\ref{thm:intro-positive-bivariate} below illustrates the
positive-measure conclusion. The general dimension bounds follow in
the next subsection; we then explain the stronger conclusions for product
measures before turning to the spectral applications.
\par}

We begin with the simplest example, the sum map
$f(x,y)=x+y$, $x,y\in \R$. For this function, a set of full Hausdorff dimension need not have
an image of positive Lebesgue measure. 
Indeed, an affine copy $C\subset[0,1]$ of Beck's compact symmetric set \cite{Bec84} satisfies
\[
\dimH C=\dimH(C+C)=1,
\qquad
\Leb(C+C)=\Leb(C-C)=0.
\]
For general smooth maps, the image can shrink to a single point. Indeed, there is a smooth function
$f:\R^2\to[0,\infty)$ that vanishes precisely on $C\times C$,
so that $f(C\times C)=\{0\}$ by Whitney's construction \cite{Whitney34}; see also \cite[Theorem~2.29]{Lee13}.
Thus, neither full Hausdorff dimension nor smoothness alone
guarantees that the image has positive Lebesgue measure.

\textit{\textbf{Question 1.} For which real analytic functions $f$ does the image
$f(A\times B)=\{f(x,y):x\in A,\ y\in B\}$
have positive Lebesgue measure whenever $A,B\subset[0,1]$ are
Borel sets of sufficiently large Hausdorff dimension?
What dimension thresholds guarantee this conclusion?}

{
The sumset obstruction persists under local analytic changes of coordinates. The following illustrative
consequence of our bivariate theory excludes this additive structure.
When both first partial derivatives are nowhere zero, it gives
positive Lebesgue measure if both input dimensions exceed $3/4$,
improving the sufficient threshold $5/6$ in \cite{Pha26}.
Theorem~\ref{thm:bivariate-dimension-expansion} below gives the full
dimension-expansion statement for Cartesian products of sets of small dimensions.
\par}

\begin{theorem}
\label{thm:intro-positive-bivariate}
Let $f$ be real analytic on a connected open neighborhood of
$[0,1]^2$, with $f_1=\partial_x f$ and $f_2=\partial_y f$
nowhere zero on $[0,1]^2$.
Suppose that $f$ does not admit a representation
$f(x,y)=G(H_1(x)+H_2(y))$ on $[0,1]^2$, with univariate
real analytic functions $G,H_1,H_2$ defined on suitable open intervals.
Then, for any Borel sets $A,B\subset[0,1]$ satisfying 
\begin{equation*}
    \dimH A,\dimH B>\frac34,
\end{equation*}
we have $\Leb(f(A\times B))>0$.
\end{theorem}

{
Beyond this positive-measure conclusion, we obtain explicit dimension
gains for every $0<\alpha<1$, strengthening the qualitative expansion
result of Raz and Zahl \cite{RZ24}. We also treat functions of three
or more variables. The role of the number of variables, and the
additional hypotheses needed to exploit all of them, will be important
for the spectral applications.
\par}

\subsection{Dimension expansion}

{
The positive-measure question is part of a quantitative problem over
the full range of input dimensions. Given Borel sets
$A_1,\ldots,A_k\subset\R$ of common Hausdorff dimension
$0<\alpha<1$, we seek bounds of the form
\[
\dimH f(A_1\times\cdots\times A_k)
\geq\alpha+\e(\alpha),
\qquad \e(\alpha)>0.
\]
Such a gain is called \emph{dimension expansion}. The aim is to
identify the functions that force expansion and to quantify the gain.
\par}

The discrete sum-product phenomenon offers a useful starting point. In 1983,
Erd\H{o}s and Szemer\'edi \cite{ES83} proved that a finite set $A$ of
positive integers cannot have both a small sumset and a small product set:
there exists $\e>0$ such that
\[
 |A+A|+|A\cdot A|\gtrsim |A|^{1+\e}.
\]
This raises a natural question: which polynomial maps can fail to expand
Cartesian products? Elekes and R\'onyai \cite{ER00} showed that for a real
bivariate polynomial $f$, such failure is governed by special algebraic
structure. Either $f$ has one of the \emph{special forms}
\[
 f(x,y)=g(h(x)+k(y))
 \quad\text{or}\quad
 f(x,y)=g(h(x)k(y))
\]
for univariate polynomials $g,h,k$; or, for all finite sets
$A,B\subset\R$ with $|A|=|B|$ sufficiently large,
\[
 |f(A\times B)|=\omega(|A|).
\]
Subsequent work quantified this expansion. Raz, Sharir, and Solymosi
\cite{RSS16} established the lower bound $|f(A\times B)|\gtrsim |A|^{4/3}$
in the nonspecial case. Solymosi and Zahl \cite{SZ24} subsequently
improved the exponent to $3/2$. Higher-dimensional extensions and related
results were developed by Schwartz, Solymosi and de Zeeuw \cite{SSdZ13},
Raz, Sharir and de Zeeuw \cite{RSdZ18}, and Raz and Shem-Tov \cite{RST20}.
Jahn and Raz \cite{JR26} obtained a rank-sensitive improvement for
$k$-variate polynomials, with an application to distinct simplex volumes
on the moment curve. The related Elekes--Szab\'o problem is studied in
\cite{ESz12}; see also the survey \cite{Z18}.

For infinite sets, Hausdorff dimension takes the place of cardinality. A fundamental result in this direction is due
to Bourgain \cite{Bou10}, who proved that for every
$\alpha\in(0,1)$ there exists $\e(\alpha)>0$ such that
every Borel set $A\subset\R$ with $\dimH A=\alpha$
admits some $t\in A$ for which
\[
\dimH(A+t A)\geq\alpha+\e(\alpha).
\]
Recently, Raz and Zahl \cite[Theorem~1.20]{RZ24} established
an Elekes--R\'onyai type dichotomy for Hausdorff
dimension, characterizing the real analytic functions
that necessarily produce dimension expansion.
Let $\Omega\subset\R^2$ be a connected open
set containing $[0,1]^2$, and let $f\in C^\omega(\Omega)$
satisfy $\partial_xf,\partial_yf\not\equiv0$.
Raz and Zahl showed that either $f$ is locally of
\emph{analytic special form} on its regular set, $f(x,y)=g(h(x)+k(y))$,
with $g,h,k$ univariate real analytic functions,
or $f(A\times B)$ has Hausdorff dimension strictly
greater than $\alpha$ whenever $A,B\subset[0,1]$
are Borel sets with dimensions at least
$\alpha\in(0,1)$.
 Their approach combines ideas from
additive combinatorics with Shmerkin's nonlinear
projection theorem \cite{Shm23}. They also established
a discretized version of the Elekes--R\'onyai theorem.

A complementary approach to dimension expansion exploits
the $L^2$-Sobolev mapping properties of generalized Radon
transforms associated with the defining function, viewed
as Fourier integral operators (FIO). In this direction,
Koh, T. Pham, and Shen \cite{KPS24} established
a positive-measure theorem for nonspecial quadratic
polynomials depending on all three variables, building on the FIO approach of Eswarathasan, Iosevich, and Taylor \cite{EIT11}.
M.-Q. Pham \cite{Pha26} extended this result
to real analytic functions by adapting the microlocal
partition method of Greenleaf, Iosevich, and Taylor
\cite{GIT22,GIT24}. In the analytic setting, this approach
gives the lower bound $\min\{3\alpha-1,1\}$ and positive
Lebesgue measure when $\alpha>2/3$. The same work also
treats bivariate analytic functions, establishing the
lower bound $\max\{2\alpha-2/3,0\}$ for
$0<\alpha\leq5/6$ and positive measure when
$\alpha>5/6$. The dimension bound gives strict expansion
when $2/3<\alpha\leq5/6$.

These results motivate a natural quantitative question.
Raz and Zahl established dimension expansion for every
$0<\alpha<1$, but without an explicit gain, while the
Fourier integral operator estimates in \cite{Pha26}
yield explicit gains only for sufficiently large input
dimensions. Our aim is to obtain explicit dimension
expansion estimates for nonspecial real analytic functions
throughout the full range $0<\alpha<1$.

Before giving the quantitative dimension bounds, we fix some notation.
Let $k\geq2$, and let $\Omega\subset\R^k$ be a
connected open neighborhood of $[0,1]^k$.
We write $f_i=\partial_{x_i}f$ and $f_{ij}=\partial_{x_i}\partial_{x_j}f$.
See Section~\ref{sec:prelim} for the full notation conventions.
Given Borel sets $A_1,\ldots,A_k\subset[0,1]$, we write
\[
f(A_1\times\cdots\times A_k)
=\{f(x_1,\ldots,x_k):x_i\in A_i,\ 1\leq i\leq k\}.
\]
We denote $A^k=A\times\cdots\times A$ ($k$ factors). The regular set $\Omega^{\reg}$ is where all first
partial derivatives of $f$ are nonzero.
\begin{definition}
Let $f\in C^\omega(\Omega)$. We say that $f$ is of
\emph{analytic special form} if there is a nonempty open set
$W\subset\Omega^{\reg}$ and real analytic maps $H_i:I_i\to\R$ for
$1\le i\le k$ and $G:J\to\R$, where $I_i$ is an open interval containing
the $i$th coordinate projection of $W$ and $J$ is an open interval containing
all sums $\sum_{i=1}^kH_i(x_i)$ for $\mathbf x\in W$, such that
\begin{equation}\label{eq:special form}
 f(x_1,\ldots,x_k)=G\bigl(H_1(x_1)+\cdots+H_k(x_k)\bigr),
 \qquad (x_1,\ldots,x_k)\in W.
\end{equation}
We call $f$ \emph{nonspecial} if it is not of analytic special form.
\end{definition}

\Needspace{12\baselineskip}
We begin with real analytic functions of two variables.
The quantitative bivariate theorem gives an explicit dimension gain for
every $0<\alpha<1$ and contains the positive-measure conclusion of
Theorem~\ref{thm:intro-positive-bivariate}.

\begin{samepage}
\begin{theorem}\label{thm:bivariate-dimension-expansion}
Let $\Omega\subset\R^2$ be a connected open neighborhood
of $[0,1]^2$, and let $f\in C^\omega(\Omega)$ satisfy
$f_1,f_2\not\equiv0$.
Then one of the following holds:
\begin{itemize}
\item[(i)] $f$ is of analytic special form
\eqref{eq:special form}.
\item[(ii)] For every pair of Borel sets
$A,B\subset[0,1]$ with
$\dimH A,\dimH B\geq\alpha$, where $0<\alpha\leq1$, one has
\begin{equation}\label{eq:bivariate-bound}
\dimH f(A\times B)\geq
\begin{cases}
4\alpha/3,&0<\alpha\leq2/3,\\[2pt]
2\alpha-\alpha^2,&2/3\leq\alpha\leq1/\sqrt2,\\[2pt]
2\alpha-1/2,&1/\sqrt2\leq\alpha\leq3/4,\\[2pt]
1,&3/4\leq\alpha\leq1.
\end{cases}
\end{equation}
Moreover, for any Borel sets $A,B\subset[0,1]$ with
$\dimH A,\dimH B>3/4$, one has $\Leb(f(A\times B))>0$.
\end{itemize}
\end{theorem}
\end{samepage}
In \eqref{eq:bivariate-bound}, the first two branches follow from incidence estimates
for transversal families of collision curves, while the
last two arise from FIO estimates
associated with the same collision geometry.
\begin{remark}
Theorem~\ref{thm:bivariate-dimension-expansion} improves the bivariate dimension
estimate in \cite[Theorem~1.4]{Pha26} when the input sets
have a common lower dimension $\alpha$. The earlier result
gives the bound $2\alpha-2/3$ for $2/3<\alpha\leq5/6$,
whereas \eqref{eq:bivariate-bound} is strictly stronger
for $2/3<\alpha<5/6$, with equality at $\alpha=5/6$.
Our theorem also provides an explicit dimension gain
for every $0<\alpha\leq2/3$ and lowers the sufficient
threshold for positive Lebesgue measure from
$\alpha>5/6$ to $\alpha>3/4$.

Demeter and O'Regan \cite[Theorem~1.3]{DO26}
established the estimate
\[
\dimH\{x(x+y):x\in A,\ y\in B\}
\geq\frac{2}{3}\bigl(\dimH A+\dimH B\bigr)
\]
whenever $0<\dimH A,\dimH B\leq2/3$.
When both dimensions equal $\alpha$, this gives
$4\alpha/3$, agreeing with the first branch of
\eqref{eq:bivariate-bound}.
Our theorem extends this exponent to arbitrary
nonspecial real analytic functions under the stated
hypotheses, while their result also captures the
dependence on unequal input dimensions for the
particular polynomial.
\end{remark}

\medskip

We next consider analytic functions of three or more variables. If no first partial derivative vanishes identically, nonspeciality
gives a dimension bound independent of the number of variables, with a
gain over the common input dimension when it lies in $(0,1)$.

\begin{samepage}
\begin{theorem}\label{thm:multivariate-dimension-expansion}
Let $k\geq3$, let $\Omega\subset\R^k$ be a connected
open neighborhood of $[0,1]^k$, and let
$f\in C^\omega(\Omega)$ satisfy
$f_i\not\equiv0$ for every $i\in\{1,\ldots,k\}$.
Then one of the following holds:
\begin{itemize}
\item[(i)] $f$ is of analytic special form
\eqref{eq:special form}.
\item[(ii)] For every Borel set $A\subset[0,1]$ with
$\dimH A=\alpha\in(0,1]$, one has
\begin{equation}\label{eq:main-bound}
\dimH f(A^k)\geq
\min\left\{\frac{3\alpha}{2},1\right\}.
\end{equation}
Moreover, if $\alpha>\frac{2}{3}$, then
$\Leb(f(A^k))>0$.
\end{itemize}
\end{theorem}
\end{samepage}
When $f$ is not of special form, Theorem~\ref{thm:multivariate-dimension-expansion} gives a uniform
dimension expansion bound throughout $0<\alpha<1$.
In particular, it provides an explicit gain for
$0<\alpha\le1/2$ and improves the bound $3\alpha-1$
from \cite{Pha26} when $1/2<\alpha<2/3$.
At $\alpha=2/3$, the image has full Hausdorff dimension,
while positive Lebesgue measure follows for $\alpha>2/3$ from
the balanced Fourier shell estimate after localization to three
variables, recovering the threshold of \cite[Theorem~5.1(ii)]{Pha26}.

Both the dimension bound and the positive-measure threshold
are independent of the number of variables $k\ge3$.
The reason is that nonspeciality can be detected by a
nonvanishing differential condition involving only three
coordinates. Analyticity then allows us to fix the remaining
variables so that the resulting three-variable function
satisfies the required nondegeneracy condition, and the
trivariate estimates apply.

{
\begin{remark}
The nonspeciality assumption is essential to a
uniform dimension expansion statement. Indeed,
the additive map
$f(x_1,\ldots,x_k)=x_1+\cdots+x_k$
need not increase the Hausdorff dimension of every
compact input set; see \cite{SS10}.
\end{remark}
\par}

{
\subsection{Absolute continuity of image measures}

The positive-measure conclusions above arise from estimates for
images of product measures. This is also the point of contact with
spectral theory: the site laws are prescribed, so we need regularity
of their pushforwards, rather than only a lower bound for the size
of the image set. A nonzero absolutely continuous measure carried
by an image gives positive Lebesgue measure of that image; the
converse does not hold for a prescribed pushforward.
\par}

{
A probability measure $\nu$ on $\R$ is $s$-Frostman if
$\nu(J)\le C_\nu|J|^s$ for every interval $J\subset\R$, where
$0<s\le1$ and $C_\nu$ is independent of $J$.
\par}

{
For functions of three or more variables, nonspeciality allows a
reduction to three coordinates. An additional collision
nondegeneracy condition lets the Fourier argument use all $k$
variables. It gives both a dimension bound depending on $k$ and an
absolute continuity criterion in terms of the sum of the Frostman
exponents. We first state its set-theoretic consequences and then
record the measure conclusion needed for the applications.

To formulate the condition, consider the collision function
\[
\Phi(\mathbf x,\mathbf y)=f(\mathbf x)-f(\mathbf y),
\]
whose zero set consists of pairs with the same image. We use the
partition approach developed by Greenleaf, Iosevich, and Taylor
\cite{GIT22,GIT24} for nonempty interior problems and extended by
Pham \cite[Theorem~3.4]{Pha26} to positive Lebesgue measure.
\par}

Partition the index set $\{1,\ldots,2k\}=\sigma_L\sqcup\sigma_R$
into two sets of size $k$, and write $\sigma=(\sigma_L,\sigma_R)$. If
$\mathbf w=(x_1,\ldots,x_k,y_1,\ldots,y_k)$, write
\begin{align*}
    x_{\sigma_L}=(w_j)_{j\in\sigma_L}, \quad x_{\sigma_R}=(w_j)_{j\in\sigma_R},
\end{align*} for the corresponding coordinate
vectors, with indices in increasing order unless stated otherwise.
Let $\pi_\sigma$ denote the coordinate permutation
$(\mathbf x,\mathbf y)\mapsto(x_{\sigma_L},x_{\sigma_R})$, and write
$\Phi_\sigma=\Phi\circ\pi_\sigma^{-1}$ for the collision function in the reordered coordinates.
Gradients are column vectors. For ordered coordinate blocks $u$ and $v$,
we use $(D_uD_v\Phi)_{ij}=\partial_{u_i}\partial_{v_j}\Phi$:
rows correspond to $u$ and columns to $v$.
The required nondegeneracy is the \textit{Phong--Stein rotational
curvature} condition \cite[Section~1, (1.5)]{PS86}: for some such partition and some
$\mathbf z\in(0,1)^k$, the Monge--Amp\`ere determinant satisfies
\begin{equation}\label{shell:bordered-condition}
\left.\det\begin{pmatrix}
D_{x_{\sigma_L}}D_{x_{\sigma_R}}\Phi_\sigma & \nabla_{x_{\sigma_L}}\Phi_\sigma\\
(\nabla_{x_{\sigma_R}}\Phi_\sigma)^T & 0
\end{pmatrix}\right|_{(x_{\sigma_L},x_{\sigma_R})=\pi_\sigma(\mathbf z,\mathbf z)}
\ne0.
\end{equation}
\begin{samepage}
\begin{theorem}
\label{thm:balanced-collision-image}
Let $k\ge3$ and let $f$ be real analytic on a connected open
neighborhood $\Omega$ of $[0,1]^k$.
Suppose that \eqref{shell:bordered-condition} holds for some
partition $\sigma=(\sigma_L,\sigma_R)$ and some $\mathbf z\in(0,1)^k$.
Then, for every Borel set $A\subset[0,1]$ with
$\dimH A=\alpha\in (0,1]$, one has
\[
\dimH f(A^k)
\ge \min\left\{k\left(\alpha-\frac{1}{2}\right)+\frac{1}{2},1\right\}.
\]
Moreover, if $$\alpha>\frac12+\frac{1}{2k},$$ then $\Leb(f(A^k))>0$.
\end{theorem}
\end{samepage}
The sufficient threshold $\alpha>1/2+1/(2k)$ for positive Lebesgue
measure of $f(A^k)$ is $2/3$ when $k=3$, $5/8$ when $k=4$,
and $3/5$ when $k=5$, approaching $1/2$ as $k\to\infty$.
Both conclusions follow from the balanced, full-rank case of
\cite[Theorem~3.4(i)--(ii)]{Pha26} after localization.
\begin{remark}
    Nonspeciality alone does not imply the Monge--Amp\`ere determinant condition
\eqref{shell:bordered-condition}. To see this, consider
\[
f(\mathbf x)=x_1+x_2x_3+x_4+\cdots+x_k,
\qquad k\ge4.
\]
This function is nonspecial because $f_2/f_1=x_3$ depends
on a third variable. However, the only nonzero mixed second
derivatives involve $x_2$ and $x_3$. Hence $D_{x_{\sigma_L}}D_{x_{\sigma_R}}\Phi_\sigma$
has rank at most two for every partition $(\sigma_L,\sigma_R)$, and the
matrix in the Monge--Amp\`ere determinant condition \eqref{shell:bordered-condition} has rank
at most four. Since this matrix has size $(k+1)\times(k+1)$,
its determinant vanishes for every partition.
\end{remark}

{
The corresponding conclusion for product measures is stated in
Theorem~\ref{thm:balanced-collision-products}(i). In particular,
under the hypotheses of
Theorem~\ref{thm:balanced-collision-image}, if $\nu_i$ is an
$s_i$-Frostman probability measure supported on $[0,1]$ for every
$1\le i\le k$, then
\[
\sum_{i=1}^k s_i>\frac{k+1}{2}
\quad\Longrightarrow\quad
f_*\left(\bigotimes_{i=1}^k\nu_i\right)\ll\Leb.
\]
The Fourier estimate in Proposition~\ref{shell:balanced-local}
gives $L^2$ densities after restriction to suitable product boxes.
Analytic localization covers a set of full product measure by such
boxes and yields absolute continuity of the full pushforward.
For two variables, Lemma~\ref{spectral:bivariate-frostman-ac}
gives the corresponding conclusion when $s_1+s_2>3/2$, under its
bivariate nondegeneracy hypothesis. These measure estimates provide
the input for the finite spectral applications below. The half-line
result uses related Fourier estimates together with a separate
argument for boundary ratios.
\par}

{
\subsection{Spectral applications}
\label{sec:intro-spectral}

We now apply the product-measure estimates to Jacobi operators.
We first study the distributions of eigenvalues and their linear
combinations, then the averaged scalar spectral measures of finite
and half-line operators. The finite results use the regularity of
eigenvalue maps; the half-line argument also exploits the recurrence
for solutions at a fixed energy.
\par}

\subsubsection{Level spacings and linear combinations of eigenvalues}
\label{spectral:stats-section}

{
The compact set $C$ used at the beginning of the introduction also
provides an obstruction for eigenvalue spacings.
Consider the $2\times2$ Jacobi matrix
\[
H=
\begin{pmatrix}
X_1 & t_1\\
t_1 & X_2
\end{pmatrix},
\qquad t_1\in\R\setminus\{0\}.
\]
Its eigenvalues $\lambda_1<\lambda_2$ have spacing
(also called level spacing) $\lambda_2-\lambda_1
=\sqrt{(X_1-X_2)^2+4t_1^2}$.
If the diagonal entries $X_1$ and $X_2$ (the site potentials)
are independent random variables with laws supported on $C$,
then $X_1-X_2$ belongs to the Lebesgue null set $C-C$
almost surely. Since the eigenvalue spacing is a Lipschitz
function of this difference, its law is singular with respect
to Lebesgue measure.

\textit{\textbf{Question 2.} Can a random Jacobi matrix with
independent site potentials having singular laws nevertheless
have absolutely continuous distributions for each of its
eigenvalues and every pairwise level spacing?}
\par}

{
An additional site changes the conclusion under suitable Frostman
assumptions. The next theorem illustrates the three-variable
measure estimate. It is a special case of
Theorem~\ref{spectral:stats-all}, which treats linear combinations
of eigenvalues whose coefficients are not all equal, in every
matrix size $n\ge3$.
\par}

{
\begin{samepage}
\begin{theorem}[]
\label{spectral:stats-three}
Consider the $3\times3$ random Jacobi matrix
\[
H=
\begin{pmatrix}
X_1 & t_1 & 0\\
t_1 & X_2 & t_2\\
0 & t_2 & X_3
\end{pmatrix},
\qquad t_1,t_2\in\R\setminus\{0\},
\]
where $X_1,X_2,X_3$ are independent real-valued random variables
with compactly supported $s_i$-Frostman laws, respectively,
for $0<s_i\le1$, $1\le i\le3$.
Suppose that $$s_1+s_2+s_3>2.$$
Write $\lambda_1<\lambda_2<\lambda_3$ for the ordered eigenvalues of $H$.
Then every pairwise level spacing
$\lambda_j-\lambda_i$
has an absolutely continuous distribution, for $1\le i<j\le3$.
\end{theorem}
\end{samepage}
\par}

{
To state the general result, we introduce a Jacobi matrix of
arbitrary finite size and its ordered eigenvalues.
\par}

Let $n\ge2$ and fix nonzero real numbers $t_1,\ldots,t_{n-1}$.
Given $\mathbf x=(x_1,\ldots,x_n)\in\R^n$, let
\begin{equation}\label{eq:finite-jacobi-matrix}
H(\mathbf x)=
\begin{pmatrix}
x_1 & t_1 & 0 & \cdots & 0\\
t_1 & x_2 & t_2 & \ddots & \vdots\\
0 & t_2 & \ddots & \ddots & 0\\
\vdots & \ddots & \ddots & x_{n-1} & t_{n-1}\\
0 & \cdots & 0 & t_{n-1} & x_n
\end{pmatrix}.
\end{equation}
Since $H(\mathbf x)$ is real symmetric and tridiagonal
with nonzero off-diagonal entries, its eigenvalues are real
and simple. We denote them by $\lambda_1(\mathbf x)<\cdots<\lambda_n(\mathbf x)$.
Let $X_1,\ldots,X_n$ be independent real-valued random
variables (site potentials) on a probability space
$(\Omega_{\mathrm{prob}},\mathcal F,\mathbb P)$,
with respective laws $\nu_1,\ldots,\nu_n$.
Write $\mathbf X=(X_1,\ldots,X_n)$.
For $c=(c_1,\ldots,c_n)\in\R^n$, define $F_c:\R^n\to\R$ by $F_c(\mathbf x)=\sum_{j=1}^n c_j\lambda_j(\mathbf x)$.
If all coefficients equal $a\in\R$, then $F_c(\mathbf X)=a\,\operatorname{tr}H(\mathbf X)
=a\sum_{i=1}^n X_i$,
whose law may remain singular.

{
The two-site spacing discussed above may have a singular law even when
the individual eigenvalue laws are absolutely continuous: take
$3/4<s<1$ in Remark~\ref{rem:stats-exclusions}(i) and apply
Theorem~\ref{spectral:finite} with $k=n=2$. The general result below
identifies a class of linear statistics for which Frostman control
at three consecutive sites suffices under the nonatomicity
assumptions stated below.
\par}

Related work includes spectral measures of random Jacobi operators
\cite{dRS11}, eigenvalue counting estimates for nonatomic site
laws \cite{CGK09}, and Poisson local eigenvalue statistics
\cite{Bou14}. Tao and Vu \cite[Theorem~15]{TV11} prove a four-moment
comparison theorem for local eigenvalue statistics in the bulk of
Wigner matrices,
while Valk\'o and Vir\'ag \cite{VV09} identify the bulk
$\mathrm{Sine}_\beta$ limit for Gaussian $\beta$-ensembles.

{
Our result concerns absolute continuity at fixed matrix size. Under
the hypotheses below, every linear eigenvalue statistic with
coefficients that are not all equal has an absolutely continuous
law. Its proof applies the three-variable product-measure theorem
after conditioning on the remaining sites.
\par}

\begin{theorem}[]
\label{spectral:stats-all}
Let $n\ge3$ and let $H(\mathbf X)$ be a random Jacobi matrix
with fixed nonzero off-diagonal entries and independent
diagonal entries $X_1,\ldots,X_n$ having compactly supported
nonatomic laws $\nu_1,\ldots,\nu_n$.
Suppose that for some $1\le r\le n-2$, the law $\nu_i$ is
$s_i$-Frostman for every $i\in\{r,r+1,r+2\}$ with $0<s_i\le1$, and $s_r+s_{r+1}+s_{r+2}>2$.
Then, for every $c=(c_1,\ldots,c_n)\in\R^n$
whose entries are not all equal, the statistic $F_c(\mathbf X)=\sum_{j=1}^n c_j\lambda_j(\mathbf X)$
has an absolutely continuous distribution.
In particular, every individual eigenvalue and every
pairwise eigenvalue difference has an absolutely continuous law.
\end{theorem}

\begin{remark}\label{rem:stats-exclusions}
    \begin{itemize}\item[]
        \item[(i)] The restriction $n\ge3$ is essential for the general result.
Indeed, the two-site example in the introduction shows that
\[
\lambda_2-\lambda_1
=\sqrt{(X_1-X_2)^2+4t_1^2}
\]
can have a singular law even when the independent site laws
are $s$-Frostman for any prescribed $0<s<1$.
Indeed, choose integers $b$ and $m$ with $2m-1<b$
and $\log m/\log b>s$. Let $X_1,X_2$ be independent with the
law of $\sum_{j\ge1}D_j b^{-j}$, where the $D_j$ are independent
and uniform on $\{0,\ldots,m-1\}$. An interval of length $b^{-j}$
meets at most three base-$b$ intervals at depth $j$, each of mass
at most $m^{-j}$, which gives the $s$-Frostman bound. The difference set is covered
at scale $b^{-j}$ by $(2m-1)^j$ intervals of length at most
$2b^{-j}$, so it is Lebesgue null.

\item[(ii)] The scalar conclusion does not extend to the trace or the joint law. Let $\eta$ be the law of $\sum_{h\ge1}D_h100^{-h}$, with independent digits uniform on $\{0,\ldots,39\}$. Each depth-$h$ cylinder has mass $40^{-h}$, and an interval of length $100^{-h}$ meets at most three corresponding base-$100$ intervals. Hence $\eta$ is $s_0$-Frostman for $s_0=\log40/\log100$ and $3s_0>2$. If $n\ge3$ and all diagonal entries have law $\eta$ independently, Theorem~\ref{spectral:stats-all} gives a density for every linear eigenvalue statistic with nonconstant coefficients.

The trace has the nonatomic law $\eta^{*n}$. Its projection modulo one is ergodic for multiplication by $100$, as a factor of the joint Bernoulli digit shift. Its first Fourier coefficient is $\widehat\eta(1)^n\ne0$, since $\operatorname{supp}\eta\subset[0,39/99]$ and the corresponding Fourier phases lie in an open semicircle. Thus the projected law differs from Haar measure and is singular with respect to it, by mutual singularity of distinct ergodic invariant probabilities \cite[Proposition~2.15 and Lemma~4.6]{EW11}. Pulling back a null carrier shows that the trace law is singular continuous. The joint eigenvalue law is also singular: the inverse image of that carrier under $(z_1,\ldots,z_n)\mapsto\sum_{j=1}^n z_j$ is Lebesgue-null by Fubini.
    \end{itemize}
\end{remark}

\subsubsection{Spectral averaging for finite Jacobi matrices}
\label{spectral:random-jacobi}
\label{spectral:many-sites}

Let $2\le k\le n$ and fix
$\mathbf b=(b_{k+1},\ldots,b_n)\in\R^{n-k}$.
For $\mathbf x=(x_1,\ldots,x_k)\in\R^k$, define
\[
H_{\mathbf b}(\mathbf x)
=H(x_1,\ldots,x_k,b_{k+1},\ldots,b_n),
\]
where $H$ is the Jacobi matrix in
\eqref{eq:finite-jacobi-matrix}.
Denote its ordered eigenvalues by
$\lambda_1^{\mathbf b}(\mathbf x)<\cdots<
\lambda_n^{\mathbf b}(\mathbf x)$
and the corresponding normalized eigenvectors by
$v_j^{\mathbf b}(\mathbf x)\in\R^n$.

Let $X_1,\ldots,X_k$ be independent real-valued random
variables on
$(\Omega_{\mathrm{prob}},\mathcal F,\mathbb P)$,
with respective laws $\nu_1,\ldots,\nu_k$.
Write $\mathbf X=(X_1,\ldots,X_k)$, whose joint law is $\boldsymbol\nu=\bigotimes_{i=1}^k\nu_i$.
For each $1\le j\le n$, the distribution of the random eigenvalue
$\lambda_j^{\mathbf b}(\mathbf X)$ is the pushforward measure
$\rho_j=(\lambda_j^{\mathbf b})_*\boldsymbol\nu$.

Given $u\in\mathbb C^n$, the scalar spectral measure of
$H_{\mathbf b}(\mathbf x)$ is defined by
\begin{equation}\label{spectral:scalar-spectral-weights}
\sigma_{H_{\mathbf b}(\mathbf x),u}(B)
=\langle \mathbf 1_B(H_{\mathbf b}(\mathbf x))u,u\rangle
=\sum_{j=1}^n
|\langle u,v_j^{\mathbf b}(\mathbf x)\rangle|^2
\delta_{\lambda_j^{\mathbf b}(\mathbf x)}(B)
\end{equation}
for every Borel set $B\subset\R$.
Here $\mathbf 1_B(H)$ is the spectral projection onto $B$ and
$\delta_a$ is the unit mass at $a$. The inner product is linear in
its first argument.
Following the spectral averaging framework
\cite{dRMSB08,KS13}, we define
\[
\overline\sigma_u(B)
=\int_{\R^k}
\sigma_{H_{\mathbf b}(\mathbf x),u}(B)
\dd\boldsymbol\nu(\mathbf x), \quad B\subset \R \text{Borel}.
\]
Spectral averaging is an important tool in spectral theory,
particularly in the study of localization by Simon and Wolff
\cite[Section~7]{SW86}. Its key feature is that averaging over
a coupling parameter can eliminate spectral mass from a fixed
Lebesgue null set. To see this in our notation, fix all diagonal
entries except the first and write
$H(x)=H(0)+xe_1e_1^*$, where $e_1$ is the first coordinate
vector of $\mathbb C^n$ and the first diagonal entry of $H(0)$ is zero. Simon and Wolff \cite[Theorem~5]{SW86} showed
that the weighted average of $\sigma_{H(x),e_1}$ against
$(1+x^2)^{-1}\,\dd x$ has the same null sets as Lebesgue measure.
In particular, for every Lebesgue null Borel set $B\subset\R$,
\[
\sigma_{H(x),e_1}(B)=0
\quad\text{for Lebesgue almost every }x.
\]
In fact, the unweighted average satisfies
\[
\int_{\R}\sigma_{H(x),e_1}(B)\,\dd x=\Leb(B);
\]
see \cite[Theorem~1.1]{Marx11}.
The exceptional set of coupling parameters in the preceding
statement is Lebesgue null, but it may carry the entire mass of
a singular probability law. Marx \cite[Theorem~1.2]{Marx11} studied the
continuity of rank-one spectral averages with respect to
Hausdorff measures under singular averaging laws.
Multi-parameter spectral averaging was developed in
\cite{BS01,dRMSB08}, while del Rio and Silva \cite{dRS11}
investigated spectral measures of random Jacobi operators.
Related questions concerning ensemble-averaged quantum dynamics
under static disorder were studied in \cite{KGB16}.

For finite Jacobi matrices, averaging over a single site
against a singular law cannot produce an absolutely continuous
spectral measure for a nonzero vector. Indeed, suppose the variable diagonal entry
$x$ has a probability law $\nu$ supported on a Lebesgue null
set $E\subset\R$. By the min--max principle, each ordered
eigenvalue $\lambda_j(x)$ is Lipschitz in $x$. Let $N\subset\R$ be a Borel null set containing
$\bigcup_{j=1}^n\lambda_j(E)$. Then
$\sigma_{H(x),u}(N)=\|u\|^2$ for $\nu$-almost every $x$.
Thus, the averaged scalar spectral measure remains singular.

This obstruction motivates averaging over several independent
site potentials. We ask whether Frostman regularity of the first
$k$ singular site laws can guarantee
$\overline\sigma_u\ll\Leb$ for every $u\in\mathbb C^n$,
and what dimensional threshold suffices.
Theorem~\ref{spectral:finite} answers this question.
Its half-line counterpart, Theorem~\ref{spectral:halfline},
gives zero averaged spectral mass to any fixed Lebesgue-null set,
including a fixed Lebesgue-null essential spectrum.
\begin{samepage}
    \begin{theorem}[]
\label{spectral:finite}
Let $2\le k\le n$ and fix
$\mathbf b=(b_{k+1},\ldots,b_n)\in\R^{n-k}$.
Suppose that $X_1,\ldots,X_k$ are independent real-valued
random variables with compactly supported $s_i$-Frostman
probability laws $\nu_i$, where $0<s_i\le1$, for $1\le i\le k$.
If $$\sum_{i=1}^k s_i>\frac{k+1}{2},$$
then the law $\rho_j\ll\Leb$ for every $1\le j\le n$ and
$\overline\sigma_u\ll\Leb$ for every $u\in\mathbb C^n$.
\end{theorem}
\end{samepage}
When all site laws have a common Frostman exponent $s$, the
sufficient condition becomes $s>1/2+1/(2k)$, giving thresholds
$3/4$, $2/3$, $5/8$, and $3/5$ for $k=2,3,4,5$, respectively.
If $s_1+s_2>3/2$, two sites already suffice. If the remaining
potentials are random but independent of $(X_1,X_2)$, the same conclusion
follows by conditioning on them. The lower thresholds for larger
$k$ arise from exploiting all $k$ variables in the Fourier
estimate for each eigenvalue law. After localization to small
rectangles satisfying the relevant collision nondegeneracy condition, the pushforward measures have $L^2$
densities. These rectangles cover a set of full product
measure, yielding absolute continuity of the eigenvalue laws.

\begin{remark}
The conclusion of Theorem~\ref{spectral:finite} can fail
if the site laws are dependent or if the randomized sites are
chosen without a placement condition.
\begin{itemize}
    \item[(i)] To see why independence is needed, consider the two-site Jacobi matrix
\[
H=
\begin{pmatrix}
X&t_1\\
t_1&X
\end{pmatrix},
\qquad t_1\ne0,
\]
where $X$ has a singular probability law.
Its eigenvalues are $X\pm|t_1|$, so both eigenvalue laws
and every averaged scalar spectral measure are singular.
This remains possible even when the law of $X$ is
$s$-Frostman for any prescribed $0<s<1$.

\item[(ii)] The randomized sites must also be chosen appropriately.
For $k\ge2$, consider a $(2k+1)\times(2k+1)$ Jacobi matrix
with unit off-diagonal entries, whose diagonal is $(b,X_1,b,X_2,\ldots,b,X_k,b)$,
where $b\in\R$ is fixed and $X_1,\ldots,X_k$ are independent
random variables. The unit vector
\[
u=\frac{1}{\sqrt{k+1}}
(1,0,-1,0,\ldots,0,(-1)^k)
\]
satisfies $Hu=bu$ for every realization.
Indeed, at each even site the two neighboring coordinates
of $u$ cancel, while at each odd site the diagonal entry
is $b$. Consequently, $\sigma_{H,u}=\overline\sigma_u=\delta_b$,
even when all randomized site laws have smooth densities.
Thus, the first $k$ sites in Theorem~\ref{spectral:finite}
cannot be replaced by an arbitrary set of $k$ sites.
By reversing the order of the matrix, the theorem also
holds when the last $k$ sites are randomized.
\end{itemize}
\end{remark}

\subsubsection{Spectral averaging for half-line Jacobi operators}

Given a bounded real sequence $\mathbf{x}=(x_j)_{j\ge1}$ and a fixed
bounded sequence of nonzero real off-diagonal entries $(t_j)_{j\ge1}$,
let $J(\mathbf{x})$ be the bounded self-adjoint operator on
$\ell^2(\mathbb N)$ with matrix
\[
J(\mathbf{x})=
\begin{pmatrix}
x_1 & t_1 & 0 & 0 & \cdots\\
t_1 & x_2 & t_2 & 0 & \cdots\\
0 & t_2 & x_3 & t_3 & \ddots\\
0 & 0 & t_3 & x_4 & \ddots\\
\vdots & \vdots & \ddots & \ddots & \ddots
\end{pmatrix}.
\]
Unlike a finite Jacobi matrix, $J(\mathbf{x})$ need not have discrete spectrum.
For $u\in\ell^2(\mathbb N)$ and a Borel set $B\subset\R$, its scalar
spectral measure and the average for a random diagonal $\mathbf X$
are
\[
\sigma_{J(\mathbf{x}),u}(B)
=\langle \mathbf 1_B(J(\mathbf{x}))u,u\rangle,
\qquad
\overline\sigma_u(B)=\mathbb E\,\sigma_{J(\mathbf X),u}(B).
\]
Our goal is to establish absolute continuity of
$\overline\sigma_u$ when only finitely many diagonal entries
near the boundary are independently randomized with singular laws.
Del Rio, Martinez, and Schulz-Baldes
\cite[Section~4, Theorem~2]{dRMSB08} proved that the averaged spectral
measure of the boundary vector $e_1$ is equivalent to Lebesgue measure on
compact intervals inside $(-2+\lambda/2,2-\lambda/2)$ for $0<\lambda<4$.
They average sufficiently many independent uniform potentials on
$[-\lambda/2,\lambda/2]$ on an initial block with unit internal hopping,
coupled to a fixed limit-point Jacobi tail.

 Our theorem allows singular Frostman laws and an almost surely bounded
tail independent of the first $m$ sites. The proof uses
the ratio of consecutive coordinates of a solution at
fixed energy, yielding estimates independent of the
remaining tail.

\begin{samepage}
\begin{theorem}[]
\label{spectral:halfline}
Let $J(\mathbf X)$ be the bounded self-adjoint Jacobi operator on
$\ell^2(\mathbb N)$ with almost surely bounded random diagonal
$\mathbf X=(X_j)_{j\ge1}$ and a fixed bounded sequence of nonzero real
hoppings $(t_j)_{j\ge1}$. Let $m=2\ell$ with $\ell\in\mathbb N$ and $\ell\ge2$ and suppose that the first $m$ diagonal entries
$X_1,\ldots,X_m$ are independent real-valued random variables
with compactly supported $s_i$-Frostman laws $\mu_i$,
where $0<s_i\le1$, for $1\le i\le m$. Assume that the remaining diagonal entries
are either fixed or form a random sequence independent of
$(X_1,\ldots,X_m)$.
If
\begin{equation}\label{spectral:intro-infinite-sum}
\sum_{i=1}^m s_i>\frac m2+1,
\end{equation}
then the averaged scalar spectral measure satisfies
\[
\overline\sigma_u\ll\Leb
\qquad\text{for every }u\in\ell^2(\mathbb N).
\]
In particular, when $s_i=s>1/2$ for $1\le i\le m$, the conclusion holds
for every even $m>2/(2s-1)$.
\end{theorem}
\end{samepage}
Our result concerns the averaged scalar measure
$\overline\sigma_u=\mathbb E\,\sigma_{J,u}$, not the bulk density of
states. For discrete Anderson models, conditioning on the other sites
and summing the rank-one estimate gives a finite-volume Wegner bound
\cite[(2.2), (2.4)--(2.5)]{CGK09}. For i.i.d.\ $s$-Frostman site laws,
volume normalization and passage to the density-of-states measure give
$O(\varepsilon^s)$ on intervals of length $\varepsilon$.
For the continuum analogue, see \cite[Theorems~1.1 and~1.3]{CHK07}. For the discrete Schr\"odinger operator on $\mathbb Z$ with i.i.d.\ compactly supported site laws satisfying a positive Frostman bound, Bourgain proved that the bulk density-of-states measure has a $C^\infty$ density \cite[Theorem~3]{Bou12}. The next example shows that an absolutely
continuous boundary average can coexist with a singular bulk density
of states.

Fix an irrational frequency $\alpha$ and a phase $\theta$.
Let $J_0$ be the half-line Jacobi operator on $\ell^2(\mathbb N)$
with unit off-diagonal entries and diagonal potential $d_j=2\cos(2\pi(j\alpha+\theta)), \quad j\ge1$.
The spectrum $\Sigma_\alpha$ of the corresponding full-line
operator has zero Lebesgue measure
\cite[Theorem~1.5]{AK06}.
We perturb the first eight sites by independent random variables
$X_1,\ldots,X_8$ having the middle-thirds Cantor law, and write
\[
J_{\mathbf X}
=J_0+\sum_{j=1}^8 X_j\langle\,\cdot\,,e_j\rangle e_j.
\]
The Cantor law is $s$-Frostman with
$s=\log2/\log3$, and $8s=\frac{8\log2}{\log3}>5=\frac82+1$.
Theorem~\ref{spectral:halfline} therefore gives $\mathbb E\,\sigma_{J_{\mathbf X},u}\ll\Leb$ for every $u\in\ell^2(\mathbb N)$.
On the other hand, since $J_{\mathbf X}-J_0$ has rank at most eight,
Weyl's theorem on the invariance of the essential spectrum
under compact perturbations gives
\[
\sigma_{\mathrm{ess}}(J_{\mathbf X})
=\sigma_{\mathrm{ess}}(J_0)=\Sigma_\alpha.
\]
For the second equality, let $J_{\mathbb Z}$ be the full-line
operator with the same phase. Extending a singular Weyl sequence for
$J_0$ by zero adds only a residual at site $0$, equal to the
first coordinate of the sequence and tending to zero by weak
convergence. Thus $\sigma_{\mathrm{ess}}(J_0)\subset
\sigma(J_{\mathbb Z})$.
Conversely, if $\tau\in\sigma(J_{\mathbb Z})$, take finitely supported
unit vectors $u_n$ with $\|(J_{\mathbb Z}-\tau\mathrm{Id})u_n\|\le n^{-1}$.
Choose integers $q_n>0$ that move $\supp u_n$ into
$\{n,n+1,\ldots\}$ and satisfy
$\operatorname{dist}(q_n\alpha,\mathbb Z)<n^{-1}$.
Translating by $q_n$ changes the potential by at most $4\pi/n$
uniformly in the site index. The translated vectors thus form a
singular Weyl sequence for $J_0$, which gives the reverse inclusion.

{
Since $\Leb(\Sigma_\alpha)=0$, absolute continuity of the
averaged spectral measure implies
\[
\mathbb E\,\sigma_{J_{\mathbf X},e_1}(\Sigma_\alpha)=0.
\]
\par}
By positivity and the cyclicity of $e_1$ for half-line
Jacobi operators, $\mathbf 1_{\Sigma_\alpha}(J_{\mathbf X})=0$ almost surely.
Thus, almost surely $J_{\mathbf X}$ has a complete
orthonormal eigenbasis with eigenvalues outside
$\Sigma_\alpha$, even though its averaged scalar
spectral measures are absolutely continuous.

Finally, let $N_{N,\mathbf X}(E)$ and $N_{N,0}(E)$ denote
the normalized eigenvalue counting functions of the
$N\times N$ Dirichlet restrictions of $J_{\mathbf X}$
and $J_0$, respectively. By unique ergodicity, the unperturbed
Dirichlet counting measures converge weakly to the density-of-states
measure $\nu_\alpha$ supported on $\Sigma_\alpha$; see
\cite[Section~1.2, immediately after (1.11)]{ALSZ24}.
Finite-rank interlacing gives
\[
\sup_{E\in\R}
|N_{N,\mathbf X}(E)-N_{N,0}(E)|
\le\frac8N.
\]
The perturbed counting measures therefore have the same weak limit
$\nu_\alpha$, which is singular since $\Leb(\Sigma_\alpha)=0$.
We have thus obtained
\[
\mathbb E\,\sigma_{J_{\mathbf X},u}\ll\Leb,
\qquad
\sigma_{J_{\mathbf X},u}\text{is pure point almost surely},
\qquad
\nu_\alpha\perp\Leb.
\]

\begin{remark}
The sufficient thresholds in Theorem~\ref{spectral:finite}
and Theorem~\ref{spectral:halfline} both approach $s=1/2$
as the number of random sites increases. The finite-matrix
condition is $k(s-1/2)>1/2$, whereas the half-line result
requires $m(s-1/2)>1$. In the finite case, Proposition~\ref{bi:radon-shell}
for $k=2$ and Proposition~\ref{shell:balanced-local} for $k\ge3$
give the Fourier shell exponent $(k+1)/2-ks$. The half-line boundary-ratio
estimate \eqref{spectral:boundary-ratio-shell-bound} has exponent $(m/2+1-ms)/2$.
Both exponents equal $1/2$ at $s=1/2$, so neither estimate gives dyadic decay there.
These are sufficient bounds; neither sharpness nor failure at the
half-line endpoint is shown here. The nonzero link to the tail
cannot simply be dropped: if $t_m=0$, a deterministic tail with
unit off-diagonal entries and diagonal $(2,0,0,\ldots)$ has the
square-summable eigenvector $(1,1/2,1/4,\ldots)$ at energy $5/2$,
so randomizing the first $m$ sites leaves an atom for this vector.
\end{remark}

{
\subsection{Main ideas of the proofs}

\noindent\textbf{Functions of two variables.}
The dimension estimates and the positive-measure conclusion use the same
collision relation, $f(x,y)=f(x',y')$, in different ways. Analytic
nonspeciality allows us to localize where
$\partial_{12}\log|f_1/f_2|\ne0$ and the other required derivatives do
not vanish. After normalizing the separation between the two input
points, this condition gives a transversal family of collision curves.
The incidence estimates of O'Regan, Wu, and Yi \cite{OWY26}, applied
with the relevant Frostman weights, then bound near-collisions and give
the first two branches of \eqref{eq:bivariate-bound}. The Fourier
argument instead uses oscillatory cancellation in the collision
relation. The classical $L^2$ estimate for nondegenerate Fourier
integral operators \cite{GS02}, combined with bounds for the normalized
pair measures, gives the remaining branches. Above the threshold
$3/4$, the resulting Fourier shell estimates are summable and give
$L^2$ densities on the localized boxes.

\noindent\textbf{Functions of three or more variables.}
For Theorem~\ref{thm:multivariate-dimension-expansion}, we characterize
analytic special form by the simultaneous vanishing of
$\Gaux_{i;j\ell}[f]=f_jf_{i\ell}-f_\ell f_{ij}$ for pairwise distinct
indices. If $f$ is nonspecial, one of these functions is nonzero on a
suitable part of the input product. Fixing the remaining coordinates
gives a three-variable function whose graphs in one variable form a
transversal family parametrized by the other two. We prove a common
image estimate for this family using the discretized Furstenberg
theorem of Orponen, Py\"or\"al\"a, and Yi \cite{OPY25}, building on
\cite{RW25}. The key step is to compare its incidence lower bound with
a single Hausdorff cover of the image, shared by all the curves. This
gives the bound $\min\{3\alpha/2,1\}$ without introducing packing
dimension. The same nonvanishing auxiliary function gives the
three-variable collision condition needed for positive Lebesgue
measure when $\alpha>2/3$.

\noindent\textbf{Using all the variables.}
The stronger hypothesis in
Theorem~\ref{thm:balanced-collision-image} permits a Fourier estimate
for the full collision relation, rather than a three-variable slice.
On each suitable product box, the balanced partition and the
nondegenerate Fourier integral operator estimate give a Fourier shell
bound for the image of a product of Frostman measures. The bound
depends on the sum of their exponents. In particular, it is summable
when this sum exceeds $(k+1)/2$, giving an $L^2$ density for each such restriction. Analytic
localization then yields absolute continuity of the full pushforward,
as stated in Theorem~\ref{thm:balanced-collision-products}. This
measure conclusion is the bridge to the finite spectral applications.

\noindent\textbf{Finite spectral applications.}
For linear eigenvalue statistics, we use covariance under common
shifts of the diagonal, algebraicity, and large-diagonal limits to
exclude additive special form when the coefficients are not all
equal. The required nondegeneracy is reduced to a three-site Jacobi
block, and the trivariate Fourier estimate applies after conditioning
on the remaining sites. For the individual eigenvalue laws in
Theorem~\ref{spectral:finite}, the Schur complement handles two
variable sites. With three or more variable sites, eigenvalue
derivative identities and the eigenvector recurrence establish the
collision determinant condition for an alternating partition.
The resulting local $L^2$ estimates yield absolute continuity of the
eigenvalue laws and hence of the averaged scalar spectral measures.

\noindent\textbf{Half-line spectral averaging.}
Theorem~\ref{spectral:halfline} requires an additional argument using
ratios of consecutive coordinates of a solution at fixed energy.
Their recurrence, together with the partition into odd and even
sites, gives the required Monge--Amp\`ere nondegeneracy. The Fourier
estimates yield bounded densities for restricted boundary-ratio laws,
uniformly in energy on suitable boxes. A Schur complement transfers
these bounds to localized finite-volume spectral averages, uniformly
in the matrix size, so that they pass to the half-line limit.
Averaging in the first site shows that the joint measure of energy
and site parameters gives no mass to the set where a denominator
vanishes. Finally, cyclicity extends absolute continuity to every
averaged scalar spectral measure.
\par}

\paragraph{\textbf{Organization.}}
Section~\ref{sec:prelim} collects the analytic and incidence tools,
and Section~\ref{sec:fourier-shell} proves the local Fourier shell estimates.
Section~\ref{sec:analytic} develops the reduction to three variables and
the balanced product theorem, both used later in the spectral arguments.
Section~\ref{sec:bivariate} combines a near-collision incidence bound
with the bivariate Fourier shell estimate.
Section~\ref{sec:spectral-averaging} applies the pushforward results to
linear eigenvalue statistics and spectral averaging for finite Jacobi
matrices. It then proves half-line averaging by a separate boundary-ratio
argument.

    \section{Preliminaries}\label{sec:prelim}
We recall the tools needed in our proofs, including
Frostman measures, Furstenberg and incidence estimates
for transversal families of curves, and $L^2$ estimates
for nondegenerate FIOs. We also recall a trivariate
positive-measure theorem and establish a discretization
lemma for Frostman measures.

We write $\dimH$ for Hausdorff dimension, $\mathcal H^s$ for
\emph{$s$-dimensional Hausdorff measure}, and $\Leb$ for Lebesgue measure on $\R$ (completed when applied to image
sets). Absolute continuity of Borel measures is understood on Borel sets. We write $X\lesssim Y$ when $X\leq CY$ for a
constant independent of the scale parameters under consideration,
and $X\approx Y$ when both $X\lesssim Y$ and $Y\lesssim X$ hold.
Dependence on fixed geometric data, exponents, and Frostman
constants may be suppressed; subscripts highlight particular dependencies.

For a function $f$ on an open set $\Omega\subset\R^k$, write
$f_i=\partial_{x_i}f$ and $f_{ij}=\partial_{x_i}\partial_{x_j}f$ for
$1\le i,j\le k$. Higher subscripts denote iterated coordinate derivatives
whenever they exist. In two variables we use $x_1=x$ and $x_2=y$.
For $f\in C^1(\Omega)$, the regular set is
$\Omega^{\reg}=\{\mathbf x\in\Omega:f_i(\mathbf x)\ne0\text{ for every }1\le i\le k\}$.
A regular box is an open product of intervals contained in this set.

All measures are nonnegative Borel measures unless otherwise stated.
We write $\boldsymbol\mu=\bigotimes_{i=1}^k\mu_i$ for a product
measure, $\mu|_A$ for the restriction to a Borel set $A$, and
$(f_*\mu)(B)=\mu(f^{-1}(B))$ for the pushforward under a measurable map $f$.
If $\mu(A)>0$ and $\mu$ is finite, the normalized restriction is
$\mu|_A/\mu(A)$. For a box $Q=\prod_{i=1}^k I_i$ with
$0<\boldsymbol\mu(Q)<\infty$ and finite marginals, this convention gives
$\boldsymbol\mu|_Q/\boldsymbol\mu(Q)=\bigotimes_{i=1}^k(\mu_i|_{I_i}/\mu_i(I_i))$.

For a finite measure $\mu$ on a measurable space $X$, write
$\|\mu\|=\mu(X)$ for its total mass. If $X=\R^d$, its Fourier transform is
\[
\widehat\mu(\xi)=\int e^{-2\pi i x\cdot\xi}\,\dd\mu(x),
\]
and similarly for integrable functions.
For a finite Borel measure $\mu$ on $\R^d$ and $0<t<d$, its
$t$-Riesz energy is
\[
I_t(\mu)=\iint |x-y|^{-t}\,\dd\mu(x)\dd\mu(y).
\]
We use $B(x,r)$ for
the open ball in the ambient Euclidean or metric space.
At a dyadic scale $\delta=2^{-n}$, $n\in\mathbb N_0$,
a \emph{dyadic $\delta$-cube} in $\R^d$ is a set of the form
\[
Q=\prod_{j=1}^d
[m_j\delta,(m_j+1)\delta),
\qquad (m_1,\ldots,m_d)\in\mathbb Z^d.
\]
When $d=1$ or $d=2$, we call these \emph{dyadic
intervals} or \emph{dyadic squares}, respectively.
The half-open convention ensures that dyadic
$\delta$-cubes partition $\R^d$ at each fixed scale.

A finite Borel measure $\mu$ on $\R^d$ is called \emph{$s$-Frostman} if there is a constant
$C_\mu>0$ such that
\[
 \mu(B(x,r))\leq C_\mu r^s
 \qquad\text{for all }x\in\R^d\text{and }r>0.
\]
By Frostman's lemma, if $E\subset\R^d$ is Borel and $0<s<\dimH E$, then there is a compact set $K\subset E$ supporting an $s$-Frostman probability measure; see
\cite[Theorem~2.7]{Mat15}. A restriction of such a measure to a Borel set of
positive mass remains $s$-Frostman after normalization. The product
of Frostman measures of exponents $s_1,\ldots,s_d$ on $\R$ is
$(s_1+\cdots+s_d)$-Frostman on $\R^d$.

\begin{lemma}[]
\label{lem:local-to-global-ac}
Let $\mu$ be a finite positive measure on a measurable space $X$ and let $F:X\to\R$ be measurable. Suppose that $(w_j)_{j\geq1}$ are nonnegative $\mu$-integrable functions with $\sum_j w_j\geq1$ $\mu$-almost everywhere. If $F_*(w_j\mu)\ll\Leb$ for every $j\geq1$, then $F_*\mu\ll\Leb$.
\end{lemma}
\begin{proof}
For any Lebesgue-null Borel set $B\subset\R$, monotone convergence gives
\[
F_*\mu(B)\leq\sum_j F_*(w_j\mu)(B)=0.
\]
\end{proof}
For a finite measure $\mu$ on $\R^d$, define its Sobolev dimension by
\[
\dim_S\mu=\sup\left\{s\in\R:\int_{\R^d}
 |\widehat\mu(\xi)|^2(1+|\xi|)^{s-d}\,\dd\xi<\infty\right\}.
\]
By \cite[Theorem~5.4(a)--(b)]{Mat15}, for a nonzero compactly
supported measure with $0<\dim_S\mu<d$, its Sobolev dimension is
$\sup\{t\in(0,d):I_t(\mu)<\infty\}$. If $\dim_S\mu>d$, the measure
has an $L^2(\R^d)$ density. We will use the following consequence of
a Fourier shell bound on the line.

\begin{lemma}[]
\label{lem:fourier-shell-consequences}
Let $\omega$ be a nonzero finite compactly supported measure on $\R$.
Suppose that for some $\beta>0$,
\[
\int_{R\le |\tau|\le 2R}|\widehat\omega(\tau)|^2\,\dd\tau
\lesssim R^{1-\beta}\qquad\text{for every dyadic }R\ge2.
\]
Then $\dim_S\omega\ge\beta$ and $I_t(\omega)<\infty$ for every
$0<t<\min\{1,\beta\}$. In particular,
$\dimH\supp\omega\ge\min\{1,\beta\}$. If $\beta>1$, then
$\omega$ has an $L^2(\R)$ density.
\end{lemma}

\begin{proof}
The Sobolev shell bounds $O(2^{j(s-\beta)})$ sum for $s<\beta$,
and boundedness of $\widehat\omega$ handles low frequencies.
For $0<t<\min\{1,\beta\}$, the same sum gives finite $t$-energy
and the support dimension bound by \cite[Theorems~3.10 and~2.8]{Mat15};
the energy weight $|\tau|^{t-1}$ is integrable at zero.
If $\beta>1$, take $s=1$ and apply Plancherel.
\end{proof}

\subsection{Discretization of Frostman measures}
We will frequently pass from Frostman measures to finite
sets at a fixed scale while preserving their nonconcentration
properties. Since dyadic cubes of the same size may carry very different masses, a Frostman bound does not directly yield a counting estimate. The following lemma therefore selects cubes of comparable mass. Their representatives satisfy the desired
nonconcentration condition, with only a logarithmic loss in $1/\delta$.

We first recall the relevant notion of
nonconcentration. Let $0<\delta\leq1$, $u\geq0$,
and $C>0$. Given a bounded set $E$ in a metric
space $X$, let $|E|_\delta$ denote the
smallest number of balls of radius $\delta$
needed to cover $E$. A nonempty bounded set
$E$ with $|E|_\delta<\infty$ is called a
\emph{$(\delta,u,C)$-set} if
\begin{equation}\label{eq:delta-set}
|E\cap B(x,r)|_\delta
\leq Cr^u|E|_\delta
\qquad\text{for all }x\in X
\text{and }\delta\leq r\leq1.
\end{equation}
Thus, the proportion of $E$ contained in a ball
of radius $r$, measured by covering numbers,
is controlled by $Cr^u$.

For a collection of dyadic $\delta$-squares, we interpret
the $(\delta,u,C)$-set condition through the set of their
centers. At scales $r\geq\delta$, this agrees, up to fixed
dimensional constants, with the covering-number condition
on the union of the squares.

\begin{lemma}[]
\label{lem:dyadic-cell-masses}
Let $0<\delta\leq 1/2$ be dyadic, and suppose that at most $C\delta^{-d}$ cells carry nonnegative masses $m_Q\leq A$, where $d\geq1$ and $A,C>0$ are fixed. For every $\gamma>0$, the cells with $m_Q<\delta^{d+\gamma}$ carry total mass at most $C\delta^\gamma$. The remaining masses lie in $O(\log(1/\delta))$ classes $[2^{-j},2^{1-j})$ indexed by $j\in\mathbb Z$, with the implicit constant depending on $A,d,\gamma$.
\end{lemma}
\begin{proof}
The discarded mass is at most $C\delta^{-d}\delta^{d+\gamma}$. The number of dyadic intervals meeting $[\delta^{d+\gamma},A]$ is $O(1+\log(1/\delta))$.
\end{proof}

The following lemma converts Frostman estimates
into this discrete nonconcentration condition.

\begin{lemma}\label{lem:discretization}
Let $d\geq1$ be an integer, $0<u\leq d$,
$\kappa\geq0$, and $c,C_0>0$.
Let $Q_0\subset\R^d$ be a bounded box.
For a sufficiently large integer $n$, set
$\delta=2^{-n}$. Suppose that $\sigma$ is a
finite Borel measure supported in $Q_0$ satisfying
\[
c n^{-\kappa}\leq\sigma(\R^d)\leq1,
\qquad
\sigma(B(x,r))\leq C_0r^u
\quad\text{for all }x\in\R^d,\ r>0.
\]
Let $E\subset Q_0$ be a Borel set with
$\sigma(Q_0\setminus E)=0$.

Then there exist a nonempty family $\mathcal Q$
of dyadic $\delta$-cubes and a number $w>0$
such that
\[
w\leq\sigma(Q)<2w
\quad\text{for all }Q\in\mathcal Q,
\qquad
\tau:=\sum_{Q\in\mathcal Q}\sigma(Q)
\geq c_1n^{-\kappa-1}.
\]
Choosing a point $p_Q\in E\cap Q$ for each
$Q\in\mathcal Q$ and setting
$P_\delta=\{p_Q:Q\in\mathcal Q\}$, we have
\begin{equation}\label{eq:discretization-size}
|P_\delta|
\geq\frac{\tau}{2w}
\geq c_2n^{-\kappa-1}\delta^{-u}.
\end{equation}
Moreover, $P_\delta$ is a
$(\delta,u,C_1n^{\kappa+1})$-set, with
$|P_\delta|_\delta\approx|P_\delta|$.

If $\Psi:E\to X$ is an $L$-bi-Lipschitz map
into a metric space $X$, with $L\geq1$, then
$\Psi(P_\delta)$ satisfies the same cardinality
lower bound and is a
$(\delta,u,C_Ln^{\kappa+1})$-set, with
\[
|\Psi(P_\delta)|_\delta
\approx_L|P_\delta|.
\]
The threshold for $n$ and the constants depend
only on $d,u,\kappa,c,C_0$, and $Q_0$, with
additional dependence on $L$ for the image
estimates.
\end{lemma}

\begin{proof}
The box $Q_0$ meets $O(\delta^{-d})$ cells. Lemma~\ref{lem:dyadic-cell-masses} with $\gamma=1$ gives discarded mass $O(\delta)$ below $\delta^{d+1}$. For large $n$, the remaining cells carry at least $\frac{c}{2}n^{-\kappa}$ and fall into at most $(d+1)n+1$ classes, since their masses lie in $[\delta^{d+1},1]$. By the pigeonhole
principle, there exists one class $\mathcal Q$
satisfying
\[
w\leq\sigma(Q)<2w
\qquad\text{for all }Q\in\mathcal Q,
\]
and
\[
\tau:=\sum_{Q\in\mathcal Q}\sigma(Q)
\geq c_1n^{-\kappa-1}.
\]
Since $\sigma(Q_0\setminus E)=0$ and
$\sigma(Q)\geq w>0$, each $Q\in\mathcal Q$
intersects $E$. Choose $p_Q\in E\cap Q$ and set $P_\delta=\{p_Q:Q\in\mathcal Q\}$.
The selected cubes are disjoint, so
$|P_\delta|=|\mathcal Q|$. Moreover, each cube
is contained in a ball of radius $C_d\delta$,
and the Frostman bound gives $w\leq\sigma(Q)\leq C\delta^u$.
Together with the lower bound for $\tau$, this yields
\[
|P_\delta|
=|\mathcal Q|
\geq\frac{\tau}{2w}
\geq c_2n^{-\kappa-1}\delta^{-u},
\]
proving \eqref{eq:discretization-size}.

We now show that $P_\delta$ inherits the
nonconcentration property of $\sigma$.
Fix $x\in\R^d$ and $r\geq\delta$.
Every cube $Q\in\mathcal Q$ with
$p_Q\in B(x,r)$ is contained in
$B(x,C_dr)$. Since the selected cubes are
disjoint and each has $\sigma$-mass at least
$w$, the Frostman bound gives
\[
|P_\delta\cap B(x,r)|\,w
\leq\sigma(B(x,C_dr))
\leq Cr^u.
\]
Moreover, $\tau<2w|P_\delta|$ and
$\tau\geq c_1n^{-\kappa-1}$. Hence
\begin{equation}\label{eq:discretization-ball-count}
\begin{aligned}
|P_\delta\cap B(x,r)|
\leq\frac{Cr^u}{w}\leq\frac{2C}{\tau}r^u|P_\delta|\leq Cn^{\kappa+1}r^u|P_\delta|.
\end{aligned}
\end{equation}
Since the points of $P_\delta$ lie in distinct
dyadic $\delta$-cubes, every ball of radius
$\delta$ contains at most $C_d$ such points.
Thus, $|P_\delta|_\delta\approx|P_\delta|$.
Together with the ball count \eqref{eq:discretization-ball-count},
this shows that $P_\delta$ is a
$(\delta,u,C_1n^{\kappa+1})$-set.

Finally, let $\Psi:E\to X$ be $L$-bi-Lipschitz. If $B_X(y,r)$ meets $\Psi(P_\delta)$, choose $p_0\in P_\delta$ with $\Psi(p_0)\in B_X(y,r)$. Every selected preimage in this ball lies in $B(p_0,2Lr)$. For $\delta\leq r\leq1$, \eqref{eq:discretization-ball-count} bounds their number by $C_Ln^{\kappa+1}r^u|P_\delta|$. At $r=\delta$, the preimage ball meets only $O_{d,L}(1)$ dyadic $\delta$-cubes. Thus each target $\delta$-ball covers at most $O_{d,L}(1)$ images, and $|\Psi(P_\delta)|_\delta\approx_L|P_\delta|$. These estimates give the required $(\delta,u,C_Ln^{\kappa+1})$-set property. Since $\Psi$ is injective, the cardinality bound is unchanged.
\end{proof}

\subsection{A discretized Furstenberg set theorem for transversal families of curves}
We now recall a
discretized Furstenberg set theorem that gives a
lower bound for the size of a union of squares
selected along a family of curves. Besides
nonconcentration of the family and the selected
squares, the theorem requires transversality
between distinct curves.

Transversality controls how two curves can meet.
Although their graphs may intersect, they cannot
agree in both value and slope at any point.
Quantitatively, the separation of their values
and slopes must control their distance in the
$C^2$ norm.

Let $I\subset\R$ be a nondegenerate compact
interval, write $C^m(I)$ for the space of real-valued $C^m$ functions on $I$, and set
\[
\|h\|_{C^m(I)}
=\sup_{z\in I}\sum_{j=0}^{m}|h^{(j)}(z)|,
\qquad m=1,2.
\]
Following \cite[Definition~2.10]{OPY25}, a family
$\mathcal H\subset C^2(I)$ is called
\emph{transversal over $I$} with transversality
constant $T\geq1$ if
\begin{equation}\label{eq:transversality}
\inf_{z\in I}
\bigl(
|h(z)-\widetilde h(z)|
+|h'(z)-\widetilde h'(z)|
\bigr)
\geq T^{-1}\|h-\widetilde h\|_{C^2(I)}
\end{equation}
for all $h,\widetilde h\in\mathcal H$.

We write $B_{C^2(I)}(a)$ for the open ball of
radius $a>0$ centered at zero in $C^2(I)$.
For $h\in C^2(I)$, a \emph{shading} $\mathcal P(h)$ is a finite collection
of dyadic $\delta$-squares intersecting $\{(z,h(z)):z\in I\}$.

With these notions in place, we recall the
following theorem of Orponen, Py\"or\"al\"a,
and Yi \cite[Theorem~1.11]{OPY25}.
\begin{samepage}
\begin{theorem}\label{thm:OPY}
Let $s\in(0,1]$, $t\in[0,2]$, $T\geq1$, and
$\eta>0$. There exist constants
$\varepsilon,\delta_0>0$, depending only on
$s,t,T,\eta$, such that the following holds for every
$\delta\in2^{-\mathbb N}\cap(0,\delta_0]$.
\\ Let $\mathcal H\subset B_{C^2([-2,2])}(1)$ be a
nonempty transversal family over $[-2,2]$ with
transversality constant $T$. Suppose that
$\mathcal H$ is a
$(\delta,t,\delta^{-\varepsilon})$-set. For each $h\in\mathcal H$,
let $\mathcal P(h)$ be a nonempty shading contained in
$[-1,1]^2$ that is a $(\delta,s,\delta^{-\varepsilon})$-set.
Suppose that there exists an integer $M\geq1$,
independent of $h$, such that
\[
|\mathcal P(h)|\geq M
\qquad\text{for every }h\in\mathcal H.
\]
Then
\begin{equation*}
\left|\bigcup_{h\in\mathcal H}\mathcal P(h)\right|
\geq
\delta^\eta
\min\left\{
\delta^{-t},
\delta^{-(s+t)/2},
\delta^{-1}
\right\}M.
\end{equation*}
\end{theorem}
\end{samepage}
Here $s$ controls nonconcentration of each shading, while
$t$ controls nonconcentration of the family $\mathcal H$
in the $C^2$ metric. The factor $M$, which gives a
uniform lower bound for the cardinality of each shading,
plays an important role in Subsection~\ref{sec:families}.
Although the constants $\varepsilon$ and $\delta_0$
are not explicit, in our applications the nonconcentration constants
grow at most as fixed powers of $\log(1/\delta)$,
and hence are bounded by $\delta^{-\varepsilon}$
for all sufficiently small $\delta$.

\subsection{Incidence theorems for transversal families}

The preceding Furstenberg theorem bounds the size of the
union of the shadings from below. The bivariate argument
instead requires an upper bound for the total number of
incidences, counting each square once for every curve
whose shading contains it.
To pass to Frostman measures later, we retain the dependence on the nonconcentration constants in two estimates of O'Regan, Wu, and Yi
\cite{OWY26}.

Let $0<\delta\le1$ be the scale; in the estimates below it is dyadic.
Throughout this subsection, the curve families are
transversal with a fixed constant $T$ and lie in a
fixed bounded ball of $C^2([-2,2])$. Their shadings
are contained in a fixed compact subset of $(-1,1)^2$.
Fixed
translations and coordinatewise rescalings,
with dyadic squares replaced by covers at a comparable scale,
ensure the domain
and norm conditions required by the cited theorems.

A $\delta$-separated family
$\mathcal F\subset C^2([-2,2])$ is called a
\emph{$(\delta,t,C)$-Katz--Tao set} if
\[
|\mathcal F\cap(h+B_{C^2([-2,2])}(r))|
\leq C(r/\delta)^t
\qquad\text{for all }h\in C^2([-2,2])
\text{and }\delta\leq r\leq1.
\]
Here separation is measured in the $C^2$ norm.
For collections of dyadic $\delta$-squares, we use
the analogous counting condition on their centers.
In particular, $|\mathcal P\cap B(z,r)|$ counts
the squares in $\mathcal P$ whose centers lie in $B(z,r)$.
Unlike \eqref{eq:delta-set}, the Katz--Tao condition
uses $\delta^{-t}$ rather than the covering number $|\mathcal F|_\delta$ as its normalization.

To describe nonconcentration in both the values and
slopes of the curves, we use the first jet evaluation map
\[
A_z:C^2(I)\to\R^2,\qquad
A_z(h)=(h(z),h'(z)),\quad z\in I,
\]
for any compact interval $I$. In the incidence estimates
below, $I=[-2,2]$.
Fix $s\ge0$ and $C_1\ge1$.
A transversal, $\delta$-separated family $\mathcal F$
is called \emph{rectangular $(\delta,2s,C_1)$-Katz--Tao}
if
\begin{equation}\label{bi:rectangular-definition}
|\{h\in\mathcal F:A_x(h)\in J_1\times J_2\}|
\leq C_1\left(\frac{|J_1||J_2|}{\delta^2}\right)^s
\end{equation}
for every $x\in[-2,2]$ and all intervals $J_1,J_2\subset\R$
with $|J_1|,|J_2|\geq\delta$.

Transversality ensures that $A_x$ is uniformly
bi-Lipschitz on $\mathcal F$
\cite[Lemma~2.6]{OWY26}.
Thus, for $\delta$-separated
families, these counting conditions agree, up to
fixed constants, with the corresponding
covering-number formulations.

We write $\pi_x(Q)$ for the horizontal projection
of a square $Q$. Given shadings $\mathcal P(h)$,
$h\in\mathcal F$, set
\[
\mathcal P=\bigcup_{h\in\mathcal F}\mathcal P(h),
\qquad
\mathcal I(\mathcal F,\mathcal P(\cdot))
=\sum_{h\in\mathcal F}|\mathcal P(h)|.
\]
Both $\mathcal I(\mathcal F,\mathcal P(\cdot))$ and the quantity
$\gamma_{\mathcal P(\cdot),t}$ defined below depend on the chosen assignment
$h\mapsto\mathcal P(h)$, not just on the union $\mathcal P$.
The two-ends condition prevents a shading from
concentrating in a small region. Given
$0<\e_2<\e_1<1$, a shading $\mathcal P(h)$ is
called \emph{$(\e_1,\e_2)$-two-ends} if
\begin{equation}\label{bi:two-ends-definition}
|\{Q\in\mathcal P(h):\operatorname{center}(Q)
\in B(z,\delta^{\e_1})\}|
\leq\delta^{\e_2}|\mathcal P(h)|
\qquad\text{for all }z\in\R^2.
\end{equation}
This center formulation implies the corresponding
covering-number condition for the union of the squares
with a fixed multiplicative constant. If necessary,
this constant can be absorbed by replacing $\e_2$
with a smaller positive exponent, leaving $\e_1$ unchanged.
A shading is called \emph{$\upsilon$-dense} if
$|\mathcal P(h)|\geq\upsilon\delta^{-1}$.

Given $0\leq t\leq1$, let $\gamma_{\mathcal P(\cdot),t}$
denote the smallest constant $C\geq1$ such that
\[
|\{Q\in\mathcal P(h):
\operatorname{center}(Q)\in B(z,r)\}|
\leq C(r/\delta)^t
\]
for all $h\in\mathcal F$, $z\in\R^2$, and
$\delta\leq r\leq1$.

We recall the two incidence estimates of \cite{OWY26}.
Only the second requires a two-ends hypothesis.

\begin{theorem}[\cite{OWY26}, Theorem 4.4]
\label{bi:owy-rectangular}
Let $0<s<1$, $\sigma=\min\{2s,2-2s\}$,
and $C_1,C_2,T\geq1$.
Suppose that
$\mathcal F\subset B_{C^2([-2,2])}(1)$
is a finite, $\delta$-separated transversal family
over $[-2,2]$ with transversality constant $T$.
Assume that $\mathcal F$ satisfies
\eqref{bi:rectangular-definition} and that every
shading $\mathcal P(h)$ is
$(\delta,\sigma,C_2)$-Katz--Tao.

Then, for every $\eta>0$,
\begin{equation}\label{bi:owy-incidence}
\mathcal I(\mathcal F,\mathcal P(\cdot))
\lesssim_\eta
C_1^{2/3}C_2^{1/3}
\delta^{-2s/3-\eta}
|\mathcal P|^{2/3}|\mathcal F|^{1/3}.
\end{equation}
The estimate holds for all sufficiently small
dyadic $\delta>0$, with the implicit constant
and scale threshold independent of $C_1,C_2$.
\end{theorem}

\begin{theorem}[\cite{OWY26}, Theorem 4.6]
\label{bi:owy-two-ends}
Let $0<s<1$, $\sigma=\min\{2s,2-2s\}$,
and $C,T\geq1$.
Suppose that
$\mathcal F\subset B_{C^2([-2,2])}(1)$
is a finite, $\delta$-separated transversal family
over $[-2,2]$ with transversality constant $T$,
and that $\mathcal F$ is a
$(\delta,2s,C)$-Katz--Tao set.

Assume that every shading $\mathcal P(h)$ is
$(\e_1,\e_2)$-two-ends and $\upsilon$-dense,
where $0<\e_2<\e_1<1$ and
$0\leq\upsilon\leq1$.
Then, for every $\eta>0$,
\begin{equation}\label{bi:owy-two-ends-bound}
|\mathcal P|
\gtrsim_\eta
\delta^{\eta+s\e_1}
C^{-1}\delta^{(2s-1)/2}
\gamma_{\mathcal P(\cdot),\sigma}^{-1/2}
\upsilon^{1/2}
\sum_{h\in\mathcal F}|\mathcal P(h)|.
\end{equation}
The estimate holds for all sufficiently small
dyadic $\delta>0$. For fixed exponents and
transversality constant, the implicit constant
and scale threshold are independent of
$C$, $\upsilon$, and $\gamma_{\mathcal P(\cdot),\sigma}$.
\end{theorem}

\subsection{An \texorpdfstring{$L^2$}{L2} estimate for nondegenerate FIOs}

We next recall a classical $L^2$ estimate for
nondegenerate FIOs, which we use in Section~\ref{sec:fourier-shell}. While incidence estimates
control the number of incidences between curves
and squares, the FIO approach exploits oscillatory
cancellation associated with the same geometry.
The key hypothesis is a nondegeneracy condition
on the phase, ensuring that the associated canonical
relation is locally a graph. We state the estimate
in the form needed later, following
\cite[Section~1.1, (1.5)--(1.6)]{GS02}.

\begin{theorem}[\cite{GS02}]
\label{thm:nondeg-fio-L2}
Let $d\ge2$ be an integer and write $(x,y,\tau)\in\R^d\times\R^d\times\R$.
Suppose that
\[
\Phi\in C^\infty(\R^d\times\R^d\times\R;\R),
\qquad
a\in C_c^\infty(\R^d\times\R^d\times\R),
\]
and that
\[
\left|\det
\begin{pmatrix}
D_xD_y\Phi & \nabla_x\partial_\tau\Phi\\
(\nabla_y\partial_\tau\Phi)^T & \partial_\tau^2\Phi
\end{pmatrix}
\right|\geq c>0
\qquad\text{on }\supp a.
\]
For $\lambda\ge1$ and $g\in C_c^\infty(\R^d)$, define the oscillatory integral operator
\[
(S_{\lambda}g)(x)
=\int_{\R^d}\int_\R
e^{i\lambda\Phi(x,y,\tau)}
a(x,y,\tau)g(y)\,d\tau\,dy.
\]
This operator extends to $L^2(\R^d)$ and satisfies
\[
\|S_{\lambda}\|_{L^2(\R^d)\to L^2(\R^d)}
\lesssim \lambda^{-(d+1)/2},
\qquad \lambda\geq1.
\]
\end{theorem}

\begin{remark}
Let $\{(\Phi_\gamma,a_\gamma)\}_{\gamma\in\Gamma}$ be a family
of phases and amplitudes as in Theorem~\ref{thm:nondeg-fio-L2}.
Write $S_{\lambda,\gamma}$ for the operator defined there with
$(\Phi,a)$ replaced by $(\Phi_\gamma,a_\gamma)$.

Suppose all the amplitudes are supported in a common compact set
$K\subset\R^{2d+1}$. Assume that the determinant in the theorem,
evaluated at $\Phi_\gamma$, has absolute value at least $c_0>0$
on $\supp a_\gamma$ for every $\gamma\in\Gamma$.
On a fixed open neighborhood $U$ of $K$, assume also that
\[
\sup_{\gamma\in\Gamma}
\left(\|\Phi_\gamma\|_{C^N(U)}+\|a_\gamma\|_{C^N(U)}\right)<\infty
\qquad\text{for every integer }N\ge0.
\]
Then
\[
\|S_{\lambda,\gamma}\|_{L^2(\R^d)\to L^2(\R^d)}
\lesssim\lambda^{-(d+1)/2},
\qquad \lambda\ge1,
\]
with a constant independent of $\lambda$ and $\gamma$.
Only finitely many of the derivative bounds enter this constant;
see \cite[Section~1.1, (1.6)]{GS02}.
Section~\ref{sec:fourier-shell} uses this uniform estimate.
\end{remark}

\subsection{A trivariate positive-measure theorem}
For comparison, we record the trivariate positive-measure theorem
of \cite[Theorem~5.1(ii)]{Pha26}. The proof of
Theorem~\ref{thm:multivariate-dimension-expansion} below obtains
the same threshold on a localized trivariate slice using
Proposition~\ref{shell:balanced-local}.

Let $V\subset\R^3$ be open and $F\in C^\infty(V;\R)$. Write
$F_j=\partial_jF$ and $F_{ij}=\partial_i\partial_jF$. Define the
\emph{degeneracy set}
\begin{equation}\label{eq:trivariate-zero-set}
\mathcal Z_{3,F}:=
\bigcap_{\substack{i,j,\ell\in\{1,2,3\}\\
i,j,\ell\text{pairwise distinct}}}
\left\{
x\in V:
F_j(x)F_{i\ell}(x)-F_\ell(x)F_{ij}(x)=0
\right\}.
\end{equation}
Observe that outside $\mathcal Z_{3,F}$, at least one of
these expressions is nonzero. The following
theorem shows that this nondegeneracy condition,
together with a sufficient sum of input dimensions,
ensures positive Lebesgue measure.

\begin{theorem}\label{thm:trivariate-positive-measure}
Let $V\subset\R^3$ be open and
$F\in C^\infty(V)$. Suppose that
$B_1,B_2,B_3\subset[0,1]$ are compact sets
of positive Hausdorff dimension satisfying
\[
B_1\times B_2\times B_3
\subset V\setminus\mathcal Z_{3,F}.
\]
If $\dimH B_1+\dimH B_2+\dimH B_3>2$, then
\[
\Leb\bigl(F(B_1\times B_2\times B_3)\bigr)>0.
\]
\end{theorem}

\begin{remark}
The localization in the cited proof can be justified
using Frostman measures. Choose
\[
0<s_i<\dimH B_i\quad (1\le i\le3),\qquad s_1+s_2+s_3>2,
\]
and $s_i$-Frostman probability measures $\mu_i$
supported on $B_i$. Choose a point in the product
of their supports. Since this point lies outside
$\mathcal Z_{3,F}$, there are closed intervals $J_i$
containing its coordinates in their interiors such that
$J_1\times J_2\times J_3\subset V$ and one of the
expressions in \eqref{eq:trivariate-zero-set} is bounded
away from zero throughout this product.
Each $\mu_i(J_i)>0$, so the normalized restrictions
$\mu_i|_{J_i}/\mu_i(J_i)$ remain $s_i$-Frostman.
Their supports are compact subsets of $B_i\cap J_i$
whose Hausdorff dimensions have sum at least
$s_1+s_2+s_3>2$.
The local nondegenerate FIO argument in
\cite[proof of Theorem~5.1]{Pha26} therefore applies
to these compact subsets and gives the stated conclusion.
\end{remark}

    \subsection{Analytic localization on products}
\label{sec:analytic-localization-products}
The first localization lemma finds a nonvanishing point in a product of compact sets (compare
\cite[Lemma~9.4]{RZ24}). The second lets us discard the analytic zero set
when the coordinate measures are nonatomic.

\begin{lemma}\label{lem:analytic-product-uniqueness}
Let $k\geq1$, and let $\Omega\subset\R^k$ be a
connected open neighborhood of $[0,1]^k$.
Suppose that $g\in C^\omega(\Omega)$ is not
identically zero. If $K_1,\ldots,K_k\subset[0,1]$
are infinite compact sets, then there exists
$a\in K_1\times\cdots\times K_k$ such that
$g(a)\ne0$.
\end{lemma}

\begin{proof}
By compactness, choose $r>0$ such that $I^k\subset\Omega$, where
$I=(-r,1+r)$. Suppose that $g$ vanishes on $K_1\times\cdots\times K_k$.
We replace the factors $K_j$ by $I$ one at a time. At step $j$, for
$1\le j\le k$, fix the preceding coordinates in $I$ and the subsequent
ones in their respective sets $K_i$. The resulting real-analytic
function of $x_j$ on $I$ vanishes on $K_j$, which has an accumulation
point in $I$. By the one-variable identity theorem it vanishes on all
of $I$. After $k$ steps, $g$ vanishes on $I^k$ and hence on the
connected domain $\Omega$, a contradiction.
This completes the proof of the lemma.
\end{proof}

\begin{lemma}
\label{lem:analytic-product-zero-mass}
Let $V\subset\R^k$ be connected and open, and let
$h\in C^\omega(V)$ be not identically zero. If $\nu_i$ is a
nonatomic finite Borel measure on $\R$ for every $1\le i\le k$, then
\begin{equation}\label{analytic:product-zero-mass}
(\nu_1\otimes\cdots\otimes\nu_k)\bigl(V\cap\{h=0\}\bigr)=0.
\end{equation}
\end{lemma}

\begin{proof}
We argue by induction on $k$. When $k=1$, the zeros of a
nonzero real-analytic function are discrete and hence countable.
Since $\nu_1$ is nonatomic, the conclusion follows.

Suppose $k\ge2$ and the assertion holds in dimension $k-1$.
Cover $V$ by countably many open rectangles
$I=I_1\times\cdots\times I_k$ contained in $V$.
By analyticity and connectedness, $h$ cannot vanish identically
on $I$. Thus there exists
$(\mathbf a',a_k)\in I$ such that
$h(\mathbf a',a_k)\ne0$.

To apply the induction hypothesis, write $\mathbf x'=(x_1,\ldots,x_{k-1})$.
The function $\mathbf x'\longmapsto h(\mathbf x',a_k)$
is real analytic and not identically zero on
$I_1\times\cdots\times I_{k-1}$.
By the induction hypothesis,
\[
(\nu_1\otimes\cdots\otimes\nu_{k-1})
\bigl(\{\mathbf x'\in I_1\times\cdots\times I_{k-1}:h(\mathbf x',a_k)=0\}\bigr)=0.
\]
Consequently, for $(\nu_1\otimes\cdots\otimes\nu_{k-1})$-almost every $\mathbf x'\in I_1\times\cdots\times I_{k-1}$, the function
$x_k\mapsto h(\mathbf x',x_k)$ is not identically zero
on $I_k$. Its zeros are therefore countable and have
zero $\nu_k$-measure. By Fubini's theorem,
\[
\begin{aligned}
&(\nu_1\otimes\cdots\otimes\nu_k)
\bigl(I\cap\{h=0\}\bigr)\\
&\qquad=
\int_{I_1\times\cdots\times I_{k-1}}
\nu_k\bigl(\{x_k\in I_k:h(\mathbf x',x_k)=0\}\bigr)
\,\dd(\nu_1\otimes\cdots\otimes\nu_{k-1})(\mathbf x')
=0.
\end{aligned}
\]
Since countably many such rectangles cover $V$,
this proves \eqref{analytic:product-zero-mass}.
\end{proof}

\section{\texorpdfstring{Fourier shell estimates for image measures}{Fourier shell estimates for image measures}}
\label{sec:fourier-shell}

In this section, we establish Fourier shell estimates for
images of products of Frostman measures under real analytic
maps. We use these estimates to bound the dimension of the images and to prove absolute continuity in the spectral applications.

\subsection{Pairing an oscillatory integral operator with Frostman measures}
We first prove a lemma estimating the pairing of an oscillatory
integral operator with Frostman measures. Since Frostman measures need not have $L^2$ densities, we use
the oscillatory phase to identify the relevant spatial frequencies.
Once these frequencies are isolated, we can apply the classical $L^2$ estimate to the localized inputs.

Let $d\ge2$,
$\Phi\in C^\infty(\R^d\times\R^d;\R)$, and
$\{a_\lambda\}_{\lambda\ge1}\subset
C_c^\infty(\R^d\times\R^d\times\R)$.
Consider the oscillatory integral operators
$T_\lambda:L^2(\R^d)\to L^2(\R^d)$ defined by
\[
T_\lambda g(u)
=\lambda\int_{\R^d}\int_\R
e^{-2\pi i\lambda\tau\Phi(u,v)}
a_\lambda(u,v,\tau)g(v)\,\dd\tau\dd v.
\]
The integral kernel of $T_\lambda$ is
\[
K_{T_\lambda}(u,v)=\lambda\int_{\R}
e^{-2\pi i\lambda\tau\Phi(u,v)}a_\lambda(u,v,\tau)\,\dd\tau.
\]
We assume that the amplitudes are supported in a common
compact set $K\subset\R^{2d+1}$, where $|\tau|\approx1$,
and satisfy
\[
\sup_{\lambda\ge1}
\|\partial^\beta a_\lambda\|_{L^\infty}<\infty
\]
for every multi-index $\beta$.
The derivatives of $\Phi$ are assumed bounded on a fixed
neighborhood of the projection of $K$ onto $\R^{2d}$.

We further assume that there exist constants
$\varepsilon_0,c_0>0$, independent of $\lambda$, such that
\begin{equation}\label{eq:pairing-bordered-condition}
\left|\det
\begin{pmatrix}
\tau D_uD_v\Phi & \nabla_u\Phi\\
(\nabla_v\Phi)^T & 0
\end{pmatrix}\right|
\ge c_0
\end{equation}
whenever $(u,v,\tau)\in\supp a_\lambda$ and
$|\Phi(u,v)|\le2\varepsilon_0$.
In addition, we assume that there exists $c_1>0$ such that
\begin{equation}\label{eq:pairing-gradient-condition}
|\nabla_u\Phi(u,v)|,\,
|\nabla_v\Phi(u,v)|\ge c_1
\end{equation}
for every $\lambda\ge1$ and
$(u,v,\tau)\in\supp a_\lambda$.

\begin{lemma}\label{bi:fractal-fio-pairing}
Let $T_\lambda$ be as above, and let
$0<\alpha_i\le d$ for $i=1,2$.
Suppose that $\rho_1,\rho_2$ are positive finite
$\alpha_i$-Frostman measures on $\R^d$, respectively,
satisfying
\[
\rho_i(B(w,\delta))\le A_i\delta^{\alpha_i},
\qquad B_i=\|\rho_i\|,
\]
for every $i=1,2$, $w\in\R^d$, and $\delta>0$.
Then, for every $N\ge1$,
\begin{equation}\label{bi:fractal-fio-bound}
\begin{aligned}
\left|\iint K_{T_\lambda}(u,v)
\,\dd\rho_1(u)\dd\rho_2(v)\right|
&\lesssim
(A_1B_1A_2B_2)^{1/2}
\lambda^{(d+1-\alpha_1-\alpha_2)/2}+C_NB_1B_2\lambda^{-N},
\end{aligned}
\end{equation}
where the constants are independent of $\lambda\ge1$.
\end{lemma}

\begin{proof}
We first localize the measures in frequency, then estimate $T_\lambda$ on $L^2$. For the error outside the selected frequencies, fix $N\ge1$ and an integer $J\ge N+2d+1$.

We begin with the frequencies selected by the oscillatory phase. Let $\xi,\eta\in\R^d$
be the Fourier variables dual to $u,v\in\R^d$, respectively.
The stationary equations for the spatial Fourier transform
of $K_{T_\lambda}$ are
\[
\xi=-\lambda\tau\nabla_u\Phi(u,v),
\qquad
\eta=-\lambda\tau\nabla_v\Phi(u,v).
\]
Since $|\tau|\approx1$ on the amplitude supports,
\eqref{eq:pairing-gradient-condition} and the uniform upper
bounds on the derivatives of $\Phi$ imply that
\[
c\lambda\le|\xi|,|\eta|\le C\lambda
\]
at every stationary point, for some constants $0<c<C$
independent of $\lambda$.

To isolate these frequencies, define
\[
\Gamma=
\overline{\bigcup_{\lambda\ge1}
\left\{
-\tau\nabla_u\Phi(u,v),\,
-\tau\nabla_v\Phi(u,v):
(u,v,\tau)\in\supp a_\lambda
\right\}}.
\]
Then $\Gamma$ is a compact subset of
$\R^d\setminus\{0\}$.
Choose a real, even function
$\varphi\in C_c^\infty(\R^d\setminus\{0\})$
that equals one on a neighborhood of $\Gamma$.
For $\lambda\ge1$, define the Fourier multiplier
$P_\lambda:L^2(\R^d)\to L^2(\R^d)$ by
\[
\widehat{P_\lambda g}(\xi)
=\varphi(\xi/\lambda)\widehat g(\xi).
\]
For a finite Borel measure $\rho$ on $\R^d$, the same multiplier
defines $P_\lambda\rho$. Since $\varphi(\cdot/\lambda)\widehat\rho$
belongs to $L^2(\R^d)$, so does $P_\lambda\rho$ by Plancherel.
Since $\varphi$ is compactly supported, the Fourier
estimate for ball averages for Frostman measures
\cite[Section~3.8, pp.~44--45]{Mat15} gives
\[
\begin{aligned}
\|P_\lambda\rho_i\|_2^2
&=\int_{\R^d}
|\varphi(\xi/\lambda)|^2
|\widehat{\rho_i}(\xi)|^2\,\dd\xi\lesssim
\int_{|\xi|\le C\lambda}
|\widehat{\rho_i}(\xi)|^2\,\dd\xi\lesssim A_iB_i\lambda^{d-\alpha_i}.
\end{aligned}
\]
Taking square roots gives
\begin{equation}\label{bi:fractal-fio-frostman}
\|P_{\lambda}\rho_i\|_2
\lesssim (A_iB_i)^{1/2}\lambda^{(d-\alpha_i)/2}.
\end{equation}
We next show that inserting the Fourier projections $P_\lambda$
introduces a rapidly decaying error. Set $E_\lambda=T_\lambda-P_\lambda T_\lambda P_\lambda$.
Since $\varphi$ is real and even, the Fourier transform
of the kernel of $E_\lambda$ is
\[
\widehat K_{E_\lambda}(\xi,\eta)
=\bigl(1-\varphi(\xi/\lambda)\varphi(\eta/\lambda)\bigr)
\widehat K_{T_\lambda}(\xi,\eta).
\]
Thus, it suffices to estimate $\widehat K_{T_\lambda}$
on the discarded frequency region
\[
\mathcal D_\lambda
=\supp\bigl(1-\varphi(\xi/\lambda)
                    \varphi(\eta/\lambda)\bigr)
\subset\R^{2d}.
\]
By the choice of $\varphi$, at least one of $\xi/\lambda$
and $\eta/\lambda$ stays a fixed positive distance from
$\Gamma$ whenever $(\xi,\eta)\in\mathcal D_\lambda$.
Consequently, the spatial phase
\[
\Psi_{\lambda,\xi,\eta}(u,v,\tau)
=\lambda\tau\Phi(u,v)+u\cdot\xi+v\cdot\eta
\]
satisfies $|\nabla_{u,v}\Psi_{\lambda,\xi,\eta}|
\gtrsim \lambda+|\xi|+|\eta|$
uniformly on the amplitude supports.
For bounded normalized frequencies, this follows from
their separation from $\Gamma$; for large frequencies,
it follows from the uniform bounds on the derivatives
of $\Phi$. Integrating by parts $J$ times in $(u,v)$ gives
\[
|\widehat K_{T_\lambda}(\xi,\eta)|
\lesssim_J
\lambda(\lambda+|\xi|+|\eta|)^{-J},
\qquad (\xi,\eta)\in\mathcal D_\lambda.
\]
Since the multiplier defining $\widehat K_{E_\lambda}$
is uniformly bounded, the same estimate holds for
$\widehat K_{E_\lambda}$. Since $J>2d$, Fourier inversion
and a change of variables yield
\[
\begin{aligned}
\|K_{E_\lambda}\|_\infty
&\le \|\widehat K_{E_\lambda}\|_{L^1(\R^{2d})}\\
&\lesssim_J
\lambda\iint_{\R^{2d}}
(\lambda+|\xi|+|\eta|)^{-J}
\,\dd\xi\dd\eta\\
&\lesssim_J \lambda^{2d+1-J}.
\end{aligned}
\]
Our choice $J\ge N+2d+1$ therefore gives
\begin{equation}\label{bi:fractal-fio-remainder}
T_\lambda=P_\lambda T_\lambda P_\lambda+E_\lambda,
\qquad
\|K_{E_\lambda}\|_\infty\lesssim_N\lambda^{-N}.
\end{equation}
It remains to bound $T_\lambda$ on $L^2$. For this we separate
the contribution near $\{\Phi=0\}$ from its complement. The augmented mixed Hessian
of the phase $\tau\Phi$ is
\[
\mathcal M_\Phi(u,v,\tau)=
\begin{pmatrix}
\tau D_uD_v\Phi&\nabla_u\Phi\\
(\nabla_v\Phi)^T&0
\end{pmatrix}.
\]
For fixed $u$, this is the Jacobian with respect
to $(v,\tau)$ of the map
\[
(v,\tau)\longmapsto
\bigl(\tau\nabla_u\Phi(u,v),\Phi(u,v)\bigr).
\]
Expanding along the last row and then along the last column of each minor gives
\begin{equation}\label{eq:pairing-monge-ampere-scaling}
\det\mathcal M_\Phi(u,v,\tau)=\tau^{d-1}\det\mathcal M_\Phi(u,v,1).
\end{equation}
On a common compact set, a uniform lower bound for the absolute determinant on $\{\Phi=0\}$ extends to a sufficiently small band $|\Phi|\leq2\varepsilon_0$ by continuity. For a continuous family, the same argument applies on a compact set that includes the phase parameters. This does not replace the gradient bounds in \eqref{eq:pairing-gradient-condition} on the full amplitude support.

On $\{\Phi=0\}$ its invertibility is precisely the local
canonical-graph condition. We use the uniform lower bound near this
set to apply the FIO estimate.

Choose $\chi\in C_c^\infty(\R)$ with
$\chi=1$ on $[-\varepsilon_0,\varepsilon_0]$
and $\supp\chi\subset(-2\varepsilon_0,2\varepsilon_0)$.
Decompose $T_\lambda=T_\lambda^{\mathrm{near}} +T_\lambda^{\mathrm{off}}$,
corresponding to the amplitudes
$\chi(\Phi)a_\lambda$ and
$(1-\chi(\Phi))a_\lambda$, respectively.
By \eqref{eq:pairing-bordered-condition}, the augmented determinant
is uniformly nonzero on the support of $\chi(\Phi)a_\lambda$.
Theorem~\ref{thm:nondeg-fio-L2} therefore gives
\[
\|T_\lambda^{\mathrm{near}}\|_{L^2\to L^2}
\lesssim\lambda\lambda^{-(d+1)/2}
=\lambda^{-(d-1)/2}.
\]
On the complementary support,
$|\Phi|\geq\varepsilon_0$. Since
\[
e^{-2\pi i\lambda\tau\Phi}
=-\frac{1}{2\pi i\lambda\Phi}
\partial_\tau e^{-2\pi i\lambda\tau\Phi},
\]
integration by parts $J$ times in $\tau$ gives the following identity
for $|\Phi(u,v)|\ge\varepsilon_0$. The off kernel vanishes when
$|\Phi(u,v)|<\varepsilon_0$.
\[
\begin{aligned}
K_{T_\lambda^{\mathrm{off}}}(u,v)
={}&\lambda(2\pi i\lambda\Phi(u,v))^{-J}\int_\R
e^{-2\pi i\lambda\tau\Phi(u,v)}
\partial_\tau^J
\bigl[(1-\chi(\Phi(u,v)))
a_\lambda(u,v,\tau)\bigr]
\,\dd\tau.
\end{aligned}
\]
The uniform derivative bounds and
$|\Phi|\geq\varepsilon_0$ imply $\|K_{T_\lambda^{\mathrm{off}}}\|_\infty
\lesssim_J\lambda^{1-J}$.
By the uniform compact spatial support, the same bound
holds for the Hilbert--Schmidt norm. In particular,
\[
\|T_\lambda^{\mathrm{off}}\|_{L^2\to L^2}
\lesssim_J\lambda^{1-J}\le\lambda^{-(d-1)/2}.
\]
Combining this with the estimate for
$T_\lambda^{\mathrm{near}}$, we obtain
\begin{equation}\label{bi:fractal-fio-l2}
\|T_\lambda\|_{L^2\to L^2}
\lesssim\lambda^{-(d-1)/2}.
\end{equation}
Finally, $P_\lambda$ is convolution with a real, even Schwartz
kernel. Fubini moves the two copies of $P_\lambda$ onto the measures, giving
\[
\iint K_{P_\lambda T_\lambda P_\lambda}(u,v)\,
\dd\rho_1(u)\dd\rho_2(v)
=\bigl\langle T_\lambda(P_\lambda\rho_2),P_\lambda\rho_1\bigr\rangle_{L^2(\R^d)}.
\]
Combining this identity with \eqref{bi:fractal-fio-remainder},
\eqref{bi:fractal-fio-l2}, and \eqref{bi:fractal-fio-frostman} gives
\[
\begin{aligned}
\left|\iint K_{T_\lambda}(u,v)
\,\dd\rho_1(u)\dd\rho_2(v)\right|
&\quad\lesssim
\lambda^{-(d-1)/2}
\|P_\lambda\rho_1\|_2
\|P_\lambda\rho_2\|_2
+C_NB_1B_2\lambda^{-N}\\
&\quad\lesssim
(A_1B_1A_2B_2)^{1/2}
\lambda^{(d+1-\alpha_1-\alpha_2)/2}
+C_NB_1B_2\lambda^{-N},
\end{aligned}
\]
This finishes the proof of the lemma.
\end{proof}

\subsection{Local collision geometry in two variables}
\label{bi:localization-section}

Let $f=f(x,y)$ be a real analytic function of two variables. We first show
that nonspeciality allows us to
localize the collision geometry to a rectangle where
the required nondegeneracy conditions hold, while
retaining positive Frostman mass. Recall that,
on the regular set, $R_f=f_1/f_2$ and $\mathcal K=\partial_{12}\log|R_f|$.
Raz and Zahl
\cite[Definition~3.1 and Lemmas~3.2--3.3]{RZ24} introduced the auxiliary function below, defined on the set where
$f_1f_2f_{12}\ne0$:
\[
\kappa_f
=f_1\partial_2\left(\frac{f_1f_2}{f_{12}}\right)
-f_2\partial_1\left(\frac{f_1f_2}{f_{12}}\right)
=\left(\frac{f_1f_2}{f_{12}}\right)^2\mathcal K
\]
to characterize analytic special forms
and establish dimension expansion results.
The same quantity arises in the FIO approach
of \cite[Section~4]{Pha26}, where its nonvanishing controls
the fold singularities of the associated canonical relation.
Wherever $f_1f_2f_{12}\neq0$,
the nonvanishing of $\kappa_f$ is equivalent
to that of $\mathcal K$.
We work directly with $\mathcal K$, which remains
well defined when $f_{12}=0$, and use its nonvanishing
to establish transversality of the normalized collision
curves.

We seek a rectangle where $f_1$, $f_2$,
$\mathcal K$, and $\partial_2R_f$ are all bounded
away from zero. To find such a rectangle that retains positive Frostman mass, we use Lemma~\ref{lem:analytic-product-uniqueness} from
Subsection~\ref{sec:analytic-localization-products}.

\begin{lemma}\label{bi:localization}
Let $\Omega\subset\R^2$ be a connected open neighborhood of $[0,1]^2$,
and let $f\in C^\omega(\Omega)$ be nonspecial with $f_1,f_2\not\equiv0$.
Let $0<s<1$, and let $\nu_0,\mu_0$ be $s$-Frostman probability
measures with compact supports in $[0,1]$. There are compact intervals
$I_1,I_2$ of positive $\nu_0$- and $\mu_0$-mass, respectively,
such that, on a neighborhood of $I_1\times I_2$,
\begin{equation}\label{bi:local-geometry}
|f_1|,\ |f_2|,\ |\mathcal K|,\ |\partial_2R_f|
\geq c>0.
\end{equation}
Moreover, the normalized restrictions
\[
\nu=\frac{\nu_0|_{I_1}}{\nu_0(I_1)},
\qquad
\mu=\frac{\mu_0|_{I_2}}{\mu_0(I_2)}
\]
remain $s$-Frostman probability measures.
The intervals may be chosen arbitrarily short about a fixed point
of the product support, within a fixed rectangle on which
\eqref{bi:local-geometry} holds.
\end{lemma}

\begin{proof}
The regular set is nonempty because $f_1$ and $f_2$ are nonzero
analytic functions on connected $\Omega$. If $\mathcal K$ vanished on
that set, the Elekes--R\'onyai criterion \cite[Lemma~3.2]{RZ24} would
give analytic special form on each connected regular component, and
hence on a regular product box. This contradicts nonspeciality.

To locate the nonvanishing quantities on the product support, we clear the denominators. Define
\begin{equation}\label{bi:numerator}
\mathcal N
:=f_1f_2^2f_{112}-f_2^2f_{12}f_{11}
-f_1^2f_2f_{122}+f_1^2f_{22}f_{12}.
\end{equation}
Put $D_f=f_{12}f_2-f_1f_{22}$. On $\Omega^{\reg}$,
\begin{equation}\label{bi:regular-numerators}
\partial_2 R_f=\frac{D_f}{f_2^2},\qquad
\mathcal K=\frac{\mathcal N}{(f_1f_2)^2}.
\end{equation}
Thus $\mathcal N$ is not identically zero, and neither is $D_f$, since its identical vanishing
would imply $\mathcal K\equiv0$ on the regular set.

Since $\Omega$ is connected, the product $\mathcal N f_1f_2D_f$
is a nonzero real analytic function on $\Omega$.
By Lemma~\ref{lem:analytic-product-uniqueness}, there exists $(x_0,y_0)\in\supp\nu_0\times\supp\mu_0$
where this product is nonzero.
Continuity then gives a rectangle containing
$(x_0,y_0)$ on which \eqref{bi:local-geometry}
holds. We may choose arbitrarily short compact
intervals $I_1,I_2$ about $x_0,y_0$, respectively,
such that $\nu_0(I_1),\mu_0(I_2)>0$ and
\eqref{bi:local-geometry} holds on a neighborhood
of $I_1\times I_2$.

Finally, define
\[
\nu=\frac{\nu_0|_{I_1}}{\nu_0(I_1)},
\qquad
\mu=\frac{\mu_0|_{I_2}}{\mu_0(I_2)}.
\]
If $C_{\nu_0},C_{\mu_0}$ are the Frostman constants
of $\nu_0,\mu_0$, then for every interval
$J\subset\R$,
\[
\nu(J)\leq
\frac{C_{\nu_0}}{\nu_0(I_1)}|J|^s,
\qquad
\mu(J)\leq
\frac{C_{\mu_0}}{\mu_0(I_2)}|J|^s.
\]
Thus, $\nu$ and $\mu$ remain $s$-Frostman
probability measures, completing the proof.
\end{proof}

\medskip
\noindent\textit{Normalized collision curves.}
\label{bi:collision-geometry}
The construction here applies to smooth $f$ on a neighborhood
of a rectangle satisfying \eqref{bi:local-geometry}. In particular, it
does not require analyticity.
Choose a compact rectangle contained in the regular set, slightly larger than the
input rectangle, on which
\eqref{bi:local-geometry} holds. Since $f_2$ does not vanish,
the smooth implicit function theorem gives a local inverse
$Y:U_{\mathrm{inv}}\to\R$, where $U_{\mathrm{inv}}\subset\R^2$
is a product neighborhood in the $(X,u)$-coordinates, with
\[
f(X,Y(X,u))=u,\qquad Y(x,f(x,y))=y.
\]
Choose the input intervals so that the image of $I_1\times I_2$
under $(x,y)\mapsto(x,f(x,y))$ is compactly contained in the domain
of $Y$. The domain then also contains $(x+rq,f(x,y))$ for all
$(x,y)\in I_1\times I_2$, all $q$ with $1/2\leq |q|\leq2$, and all sufficiently
small $r>0$. For such $r$ and $q$, define
\begin{equation}\label{bi:normalized-collision-curve}
H_r(x,q;y)
=\frac{Y(x+rq,f(x,y))-y}{r}.
\end{equation}
The variables $(x,q)$ parametrize the curves and $y$ is the
curve variable. For an exact collision with $x'=x+rq$, the
second coordinate is
$y'=Y(x+rq,f(x,y))$, so that
$(y'-y)/r=H_r(x,q;y)$.

The singularity at $r=0$ is removable. Indeed, fixing
$u=f(x,y)$, implicit differentiation gives
\[
\partial_XY(X,u)
=-\frac{f_1(X,Y(X,u))}{f_2(X,Y(X,u))}
=-R_f(X,Y(X,u)).
\]
Since $Y(x,u)=y$, the fundamental theorem of calculus yields
\[
\begin{aligned}
H_r(x,q;y)
&=\frac{Y(x+rq,u)-Y(x,u)}{r}=q\int_0^1\partial_XY(x+trq,f(x,y))\,\dd t.
\end{aligned}
\]
Thus $H_r$ extends smoothly to $r=0$, with
$H_0(x,q;y)=-qR_f(x,y)$. The convergence is uniform together
with the derivatives on fixed compact parameter regions.

For a differentiable curve $h$, write $A_y(h)=(h(y),h'(y))$.
To see the source of transversality, consider the first jet map
\[
(x,q)\longmapsto A_y(H_0(x,q;\cdot))
=\bigl(H_0(x,q;y),\partial_yH_0(x,q;y)\bigr).
\]
Its parameter Jacobian is
\[
D_{(x,q)}A_y(H_0)
=\begin{pmatrix}
-q\partial_1R_f&-R_f\\
-q\partial_{12}R_f&-\partial_2R_f
\end{pmatrix},
\]
and hence
\begin{equation}\label{bi:collision-jet}
\begin{aligned}
\det D_{(x,q)}A_y(H_0)
&=q\bigl((\partial_1R_f)(\partial_2R_f)-R_f\partial_{12}R_f\bigr)=-qR_f^2\mathcal K.
\end{aligned}
\end{equation}
By \eqref{bi:local-geometry}, this determinant, as well as
$\partial_qH_0=-R_f$ and
$\partial_q\partial_yH_0=-\partial_2R_f$, is bounded
away from zero in absolute value for $1/2\leq |q|\leq2$.

\subsection{Fourier shell estimates in two variables}
\label{bi:radon}
We now establish a Fourier shell estimate for the image measure
$\omega=f_*(\nu\otimes\mu)$. The key is to analyze the collision
relation $f(x,y)=f(x',y')$ using the rescaled curves introduced
in Subsection~\ref{bi:localization-section}.
Their values and slopes provide the nondegeneracy needed to
apply FIO estimates uniformly as the
horizontal separation $r\to0$.
For the spectral applications, the two Frostman exponents may differ. We therefore keep them separate in the estimate.

\begin{proposition}[]
\label{bi:radon-shell}
Let $f\in C^\infty(U)$, where $U\subset\R^2$ is an open
rectangle containing $(x_0,y_0)$, and suppose that
\eqref{bi:local-geometry} holds on $U$.
Then there exist compact intervals $I_1,I_2$ containing
$x_0,y_0$ in their interiors, respectively, such that
the following holds.

Let $\nu$ and $\mu$ be $s_1$- and $s_2$-Frostman
probability measures supported on $I_1$ and $I_2$,
respectively, where $0<s_1,s_2\le1$.
Set $\omega=f_*(\nu\otimes\mu)$.
Then, for every nonnegative
$\zeta\in C_c^\infty(\R\setminus\{0\})$,
\begin{equation}\label{bi:radon-shell-bound}
\int_\R \zeta(t/R)|\widehat\omega(t)|^2\,\dd t
\lesssim_\zeta R^{3/2-s_1-s_2},
\qquad R\ge2.
\end{equation}
In particular,
\begin{equation}\label{spectral:unequal-two-site-shell}
\int_{R\le|t|\le2R}
|\widehat\omega(t)|^2\,\dd t
\lesssim R^{3/2-s_1-s_2}.
\end{equation}
The implicit constants may depend on the local geometry,
the Frostman exponents and constants, and $\zeta$,
but are independent of $R$.

Consequently, if $s_1+s_2>1/2$, then
\[
\dimH f(\supp\nu\times\supp\mu)
\ge \min\left\{s_1+s_2-\frac12,1\right\}.
\]
Moreover, if $s_1+s_2>3/2$, then $\omega$ has an
$L^2(\R)$ density.
\end{proposition}

\begin{proof}
Since $\omega=f_*(\nu\otimes\mu)$, its Fourier transform is
\[
\widehat\omega(t)
=\iint e^{-2\pi itf(x,y)}\,\dd\nu(x)\dd\mu(y).
\]
To estimate the Fourier shell integral, set
\[
I_R=\int_\R\zeta(t/R)|\widehat\omega(t)|^2\,\dd t.
\]
Expanding the square and applying Fubini's theorem gives
\[
I_R=\iiiint K_R\bigl(f(x,y)-f(x',y')\bigr)
\,\dd\nu(x)\dd\nu(x')\dd\mu(y)\dd\mu(y'),
\]
where
\[
K_R(v)=\int_\R e^{-2\pi itv}\zeta(t/R)\,\dd t
=R\widehat\zeta(Rv).
\]
Since $\zeta$ is smooth and compactly supported,
its Fourier transform decays rapidly. Thus, for every $N\ge1$,
\begin{equation}\label{bi:radon-kernel-decay}
|K_R(v)|\lesssim_{N,\zeta}R(1+R|v|)^{-N}.
\end{equation}
The decay of $K_R$ localizes the image difference. To compare this with the input differences, we split the integral according to $|x'-x|$.

Choose a nonnegative $\psi\in C_c^\infty(\R)$ satisfying
\[
\supp\psi\subset\{q:1/2<|q|<2\},\qquad
\sum_{r\in2^{\mathbb Z}}\psi(a/r)=1\quad(a\ne0).
\]
Here $r$ ranges over positive dyadic numbers.
Since $\nu$ has no atoms, the diagonal $x=x'$ has
zero $\nu\otimes\nu$-measure. Since $\psi\ge0$ and
$|K_R|\le R\|\zeta\|_1$, dominated convergence gives
\begin{equation*}
I_R=\sum_{r\in2^{\mathbb Z}} I_{R,r},
\end{equation*}
where for each $r\in 2^{\mathbb{Z}}$
\[
\begin{aligned}
I_{R,r}=\iiiint&
K_R\bigl(f(x,y)-f(x',y')\bigr)
\psi\left(\frac{x'-x}{r}\right)\dd\nu(x)\dd\nu(x')\dd\mu(y)\dd\mu(y').
\end{aligned}
\]
Each term is localized to the first-coordinate separation
$|x'-x|\approx r$.

On the support of the dyadic cutoff write $x'=x+rq$, with
$1/2<|q|<2$, and put $z=(y'-y)/r$.
Fix a compact rectangle $J_1\times J_2\subset U$ about $(x_0,y_0)$.
Choose a smaller compact rectangle $J'\Subset\operatorname{int}(J_1\times J_2)$ about $(x_0,y_0)$. The construction in
Subsection~\ref{bi:localization-section} identifies the zero set of
$f(x,y)-f(x+rq,y+rz)$ for $r>0$ with the graph $z=H_r(x,q;y)$ in
\eqref{bi:normalized-collision-curve}. The smooth extension at $r=0$
bounds $H_r$ on the compact set $\mathcal C_*=\{(x,q,y,r):(x,y)\in J',\ 1/2\le|q|\le2,\ 0\le r\le r_*\}$ for some $r_*>0$. Set
\[
C_f=\sup_{J_1\times J_2}|f_1|,\qquad
c_f=\inf_{J_1\times J_2}|f_2|>0,
\]
and choose $D_0>\max\{\sup_{\mathcal C_*}|H_r|,4C_f/c_f\}$. Choose
$\chi_0\in C_c^\infty(\R)$ with $0\leq\chi_0\leq1$, equal to one
on $[-D_0,D_0]$ and supported in $(-2D_0,2D_0)$.
\par
Fix $0<r_0\le r_*$ small enough that $(x+rq,y+rz)$ stays in
$J_1\times J_2$ for $(x,y)$ in $J'$,
$|q|\leq2$, $|z|\leq2D_0$, and $0\leq r\leq r_0$.
By \eqref{bi:collision-jet} and smooth convergence to $H_0$, we
may also require
\begin{equation}\label{bi:radon-curvature}
\left|\det D_{(x,q)}(H_r,\partial_yH_r)\right|
\geq\kappa_0>0
\end{equation}
uniformly on $J'$ for $1/2\leq|q|\leq2$ and
$0\leq r\leq r_0$. Take $I_1\times I_2$ in the interior of $J'$ with
$2\operatorname{diam}(I_1)<r_0$. Every contributing scale
then satisfies $r<r_0$.

Split $I_{R,r}=I_{R,r}^{\mathrm{near}}+
I_{R,r}^{\mathrm{off}}$, where
\[
\begin{aligned}
I_{R,r}^{\mathrm{near}}=\iiiint&
K_R\bigl(f(x,y)-f(x',y')\bigr)
\psi\left(\frac{x'-x}{r}\right)\chi_0\left(\frac{y'-y}{r}\right)
\,\dd\nu(x)\dd\nu(x')\dd\mu(y)\dd\mu(y'),
\end{aligned}
\]
and $I_{R,r}^{\mathrm{off}}=I_{R,r}-I_{R,r}^{\mathrm{near}}$.
All sums in $r$ run over dyadic scales $r\in2^{\mathbb Z}$.
Further dividing the near contribution at the scale $R^{-1}$,
we obtain
\begin{equation}\label{bi:radon-three-pieces}
I_R=I_R^{\mathrm{near,small}}+I_R^{\mathrm{off}}
+I_R^{\mathrm{near,large}},
\end{equation}
where
\[
\begin{aligned}
I_R^{\mathrm{near,small}}
&=\sum_{r<R^{-1}}I_{R,r}^{\mathrm{near}},\qquad
I_R^{\mathrm{off}}
=\sum_{r\in2^{\mathbb Z}}I_{R,r}^{\mathrm{off}},\qquad
I_R^{\mathrm{near,large}}
=\sum_{r\geq R^{-1}}I_{R,r}^{\mathrm{near}}.
\end{aligned}
\]
We estimate the first two terms directly and use
the FIO estimate for the third.
\\
For the first term $I_R^{\mathrm{near,small}}$, we observe that on the support of the integrand defining $I_{R,r}^{\mathrm{near}}$,
\[
|x'-x|\leq2r,\qquad |y'-y|\leq2D_0r.
\]
The product Frostman bounds show that these
configurations have measure $O(r^{s_1+s_2})$.
Using $|K_R|\lesssim R$ from \eqref{bi:radon-kernel-decay}, we obtain
$|I_{R,r}^{\mathrm{near}}|\lesssim Rr^{s_1+s_2}$.
Summing over $r<R^{-1}$ gives
\begin{equation}\label{bi:radon-small-bound}
|I_R^{\mathrm{near,small}}|
\lesssim R\sum_{r<R^{-1}}r^{s_1+s_2}
\lesssim R^{1-s_1-s_2}.
\end{equation}
For the off-incidence term, set $V=|y'-y|$. Its cutoffs require
$V\geq D_0r$ and $|x'-x|\leq2r$. Since $f_2$ has fixed sign on the
regular rectangle, the mean value theorem gives
\[
|f(x,y)-f(x',y')|\geq c_fV-C_f|x'-x|\geq c_fV/2.
\]
Sum the nonnegative horizontal partition weights first; their sum is
at most one. On each dyadic vertical annulus $a\leq V<2a$, we have
$|x'-x|\lesssim a$, so the product Frostman mass is $O(a^{s_1+s_2})$.
By \eqref{bi:radon-kernel-decay}, for $N>s_1+s_2$,
\begin{equation}\label{bi:radon-off-bound}
|I_R^{\mathrm{off}}|
\lesssim_N R\sum_{a\in 2^{\mathbb Z}}(1+Ra)^{-N}a^{s_1+s_2}
\lesssim R^{1-s_1-s_2}.
\end{equation}
The last bound follows by geometric summation on either side of $a=R^{-1}$.

It remains to estimate $I_R^{\mathrm{near,large}}
=\sum_{r\geq R^{-1}}I_{R,r}^{\mathrm{near}}$.
Put $u=(x,q)$, $v=(y,z)$, and define $\phi_r(u,v) =\frac{f(x,y)-f(x+rq,y+rz)}{r}$.
By the fundamental theorem of calculus, one has
\[
\phi_r
=-\int_0^1\bigl(
 qf_1(x+\theta rq,y+\theta rz)
 +zf_2(x+\theta rq,y+\theta rz)
 \bigr)\,\dd\theta.
\]
Hence $\phi_r$ extends smoothly to $r=0$, with
$\phi_0=-qf_1(x,y)-zf_2(x,y)$. Moreover,
\begin{equation}\label{bi:radon-gradients}
\partial_q\phi_r=-f_1(x+rq,y+rz),\qquad
\partial_z\phi_r=-f_2(x+rq,y+rz).
\end{equation}
We next localize the product measures to pairs satisfying
$|x'-x|\approx r$ and $|y'-y|\lesssim r$,
respectively, and then rescale their displacements.
For $0<r\le r_0$, define the linear map $\pi_r:\R^2\to\R^2$ by $\pi_r(a,b)=\left(a,\frac{b-a}{r}\right)$.
Define the pushforward measures $m_r,n_r$ on $\R^2$ by
\[
\begin{aligned}
m_r&=(\pi_r)_*
\left[\psi\left(\frac{x'-x}{r}\right)
\,\dd\nu(x)\dd\nu(x')\right],\\
n_r&=(\pi_r)_*
\left[\chi_0\left(\frac{y'-y}{r}\right)
\,\dd\mu(y)\dd\mu(y')\right].
\end{aligned}
\]
The Frostman bounds give
\[
m_r(\R^2)\lesssim r^{s_1},
\qquad
n_r(\R^2)\lesssim r^{s_2}.
\]
For an arbitrary center $u_0=(\bar x,\bar q)\in\R^2$ and radius
$0<\rho\leq1$, the condition $(x,(x'-x)/r)\in B(u_0,\rho)$ implies
\[
 |x-\bar x|<\rho,\qquad |x'-(x+r\bar q)|<r\rho .
\]
By Fubini and the Frostman bound in each coordinate, using the
bounded weights in $m_r$, the outer interval contributes
$O(\rho^{s_1})$ and each fiber contributes $O((r\rho)^{s_1})$. Thus
$m_r(B(u_0,\rho))\lesssim r^{s_1}\rho^{2s_1}$.

The same argument in $y,y'$ gives the bound for $n_r$, and
the total masses handle $\rho>1$. Uniformly for centers
$u_0,v_0\in\R^2$ and radii $\rho>0$, we have
\begin{equation}\label{bi:radon-pair-measures}
\begin{aligned}
m_r(\R^2)&\lesssim r^{s_1},&
m_r(B(u_0,\rho))&\lesssim r^{s_1}\rho^{2s_1},\\
n_r(\R^2)&\lesssim r^{s_2},&
n_r(B(v_0,\rho))&\lesssim r^{s_2}\rho^{2s_2}.
\end{aligned}
\end{equation}
Set $\lambda=Rr\geq1$. Choose a fixed smooth compactly supported
cutoff $\chi(u,v)$ equal to one on $\supp(m_r\otimes n_r)$ for every
contributing $r$, with support in the regular localization region.
On a common open neighborhood of $\supp\chi$, the points
$(x,y)$ and $(x+rq,y+rz)$ lie in $U$ for every $0\le r\le r_0$.
Multiply $\phi_r$ by a fixed smooth cutoff compactly supported in that
neighborhood and equal to one near $\supp\chi$, then extend by zero. The integral formula for $\phi_r$
gives uniform derivative bounds on this fixed support. Using the same
notation for the extension, define $T_{\lambda,r}:L^2(\R^2)\to L^2(\R^2)$ by
\[
(T_{\lambda,r}g)(u)
=\lambda\int_{\R^2}\int_\R
e^{-2\pi i\lambda\tau\phi_r(u,v)}
\zeta(\tau)\chi(u,v)g(v)
\,\dd\tau\dd v.
\]
Its Schwartz kernel is
\[
K_{T_{\lambda,r}}(u,v)
=\lambda\chi(u,v)\int_\R
e^{-2\pi i\lambda\tau\phi_r(u,v)}
\zeta(\tau)\,\dd\tau.
\]
Where the pairing cutoff is nonzero, we have
$f(x,y)-f(x',y')=r\phi_r(u,v)$. Since $\lambda=Rr$, the change of variables
$t=R\tau$ gives
\begin{equation}\label{bi:radon-fio-pairing}
I_{R,r}^{\mathrm{near}}
=r^{-1}\iint K_{T_{\lambda,r}}(u,v)
\,\dd m_r(u)\dd n_r(v).
\end{equation}
We now verify the hypotheses of
Lemma~\ref{bi:fractal-fio-pairing}.
On a neighborhood of $\supp\chi$, the derivative $\partial_z\phi_r$
is bounded away from zero and the zero set of $\phi_r$ is given by
$z=H_r(x,q;y)$. On this neighborhood we may write
\[
\phi_r(u,v)
=\widetilde h_r(u,v)(z-H_r(x,q;y)),
\]
where $\widetilde h_r$ is smooth and uniformly
nonvanishing. On this local zero set, the augmented mixed Hessian satisfies
\[
\begin{aligned}
&\det
\begin{pmatrix}
\tau D_uD_v\phi_r&\nabla_u\phi_r\\
(\nabla_v\phi_r)^T&0
\end{pmatrix}=
\tau\widetilde h_r^3
\det D_{(x,q)}
(H_r,\partial_yH_r).
\end{aligned}
\]
By \eqref{bi:radon-curvature}, the uniform nonvanishing of
$\widetilde h_r$, and the support of $\zeta$ away from zero,
$|\det\mathcal M_{\phi_r}|$ is bounded below on the collision set.
Compactness of $[0,r_0]\times\supp(\chi\zeta)$ and smoothness at
$r=0$ extend this bound, with a possibly smaller constant, to
$|\phi_r|\leq 2\varepsilon_0$ on the amplitude support for some
$\varepsilon_0>0$ independent of $r$. This verifies
\eqref{eq:pairing-bordered-condition}.
Equations \eqref{bi:radon-gradients} and \eqref{bi:local-geometry}
give \eqref{eq:pairing-gradient-condition}, and the smooth extension
to $r=0$ gives uniform bounds on the derivatives of $\phi_r$.
Apply \eqref{bi:fractal-fio-bound} with $d=2$ and $\alpha_i=2s_i$ for
$i=1,2$. By \eqref{bi:radon-pair-measures},
\[
\left|\iint K_{T_{\lambda,r}}(u,v)\,\dd m_r(u)\dd n_r(v)\right|
\lesssim r^{s_1+s_2}\bigl(\lambda^{3/2-s_1-s_2}+\lambda^{-N}\bigr).
\]
Since $\lambda\geq1$ and $s_1+s_2\leq2$, the remainder is absorbed for
$N\geq2$. Multiplying by $r^{-1}$ in \eqref{bi:radon-fio-pairing} gives
\[
|I_{R,r}^{\mathrm{near}}|\lesssim
r^{s_1+s_2-1}(Rr)^{3/2-s_1-s_2}=R^{3/2-s_1-s_2}r^{1/2}.
\]
The half-power of $r$ makes the dyadic sum geometric:
\begin{equation}\label{bi:radon-fio-bound}
|I_R^{\mathrm{near,large}}|
\lesssim R^{3/2-s_1-s_2}\sum_{R^{-1}\leq r<r_0}r^{1/2}
\lesssim R^{3/2-s_1-s_2}.
\end{equation}
Combining this with
\eqref{bi:radon-three-pieces},
\eqref{bi:radon-small-bound}, and
\eqref{bi:radon-off-bound} proves
\eqref{bi:radon-shell-bound}.

We now pass from the shell estimate to a Hausdorff dimension bound.
Choose a nonnegative
even function
$\zeta_0\in C_c^\infty(\R\setminus\{0\})$
with $\zeta_0\geq1$ on $1\leq|t|\leq2$. This choice in
\eqref{bi:radon-shell-bound} gives \eqref{spectral:unequal-two-site-shell}.
When $s_1+s_2>1/2$, apply Lemma~\ref{lem:fourier-shell-consequences} to \eqref{spectral:unequal-two-site-shell} with $\beta=s_1+s_2-1/2$. Since $\omega$ is supported on $f(\supp\nu\times\supp\mu)$, the lemma gives the asserted dimension bound. It also gives an $L^2$ density when $s_1+s_2>3/2$.
This finishes the proof of the proposition.
\end{proof}

We record the following consequence of Proposition~\ref{bi:radon-shell}
and analytic localization.

\begin{lemma}[]
\label{spectral:bivariate-frostman-ac}
Let $f:\R^2\to\R$ be real analytic, and suppose that
\begin{equation}\label{spectral:nondegeneracy}
f_1f_2\,\partial_yR_f\,
\partial_{xy}\log|R_f|\ne0,
\qquad R_f=\frac{f_1}{f_2},
\end{equation}
at some point. Let $\nu,\mu$ be compactly supported
$s_1$- and $s_2$-Frostman probability measures on $\R$,
respectively, where $0<s_1,s_2\le1$.
If $s_1+s_2>3/2$, then $f_*(\nu\otimes\mu)\ll\Leb$.
\end{lemma}

\begin{proof}
Let $\mathcal N$ be as in \eqref{bi:numerator} and $D_f$ as in
\eqref{bi:regular-numerators}. Set $h=f_1f_2D_f\mathcal N$.
By \eqref{bi:regular-numerators}, $h$ is nonzero at the point
in \eqref{spectral:nondegeneracy}. The marginals are nonatomic,
so Lemma~\ref{lem:analytic-product-zero-mass} gives
$(\nu\otimes\mu)(\{h=0\})=0$.

Around every point where $h\ne0$, the conditions in
\eqref{bi:local-geometry} hold on a sufficiently small rectangle.
Choose countably many compact product boxes $Q_j$ whose interiors
cover $\{h\ne0\}$ and on which Proposition~\ref{bi:radon-shell}
applies. For every box of positive product mass, the normalized
marginal restrictions retain the Frostman exponents. Since
$s_1+s_2>3/2$, their pushforward has an $L^2$ density.
Apply Lemma~\ref{lem:local-to-global-ac} to $\nu\otimes\mu$ with
$w_j=\mathbf{1}_{Q_j}$: the boxes cover $(\nu\otimes\mu)$-almost
every point, and each restriction of positive mass is a multiple of
its normalization. Thus $f_*(\nu\otimes\mu)\ll\Leb$.
\end{proof}

\subsection{Fourier shell estimates in three or more variables}
\label{shell:higher}

We now extend the Fourier shell estimates to smooth functions
$f:U\to\R$, where $U\subset\R^k$ is open and $k\ge3$.
For every $1\le i\le k$, let $\nu_i$ be an $s_i$-Frostman
probability measure on $\R$, with $0<s_i\le1$, supported
on a compact interval $I_i$ such that
$I_1\times\cdots\times I_k\subset U$.
Set
\[
\boldsymbol\nu=\bigotimes_{i=1}^k\nu_i,
\qquad
\omega=f_*\boldsymbol\nu.
\]
As in the bivariate case, expanding the squared Fourier
transform gives
\[
|\widehat\omega(\tau)|^2
=\iint_{U\times U}
e^{-2\pi i\tau\Phi(\mathbf x,\mathbf y)}
\,\dd\boldsymbol\nu(\mathbf x)
\,\dd\boldsymbol\nu(\mathbf y),
\]
where the collision function $\Phi:U\times U\to\R$ is defined by $\Phi(\mathbf x,\mathbf y)=f(\mathbf x)-f(\mathbf y)$.
 The key observation is that the $2k$ scalar coordinates
of $(\mathbf x,\mathbf y)$ can be partitioned into two
groups of size $k$. In two variables the Monge--Amp\`ere determinant of the
unrescaled collision vanishes at every diagonal point for every balanced
partition, which is why the preceding proof rescales the collision. For
$k\ge3$ a balanced partition can instead be nondegenerate on the diagonal. The partition approach to generalized Radon transforms was
developed by Greenleaf, Iosevich, and Taylor
\cite{GIT22,GIT24} to study the nonempty interior of
configuration sets. Pham \cite[Theorem~3.4]{Pha26}
extended this approach to obtain positive Lebesgue measure
for general $k$-point configuration sets.

In our setting, the two copies of the input supply the $2k$ scalar coordinates of the collision relation. Write
\[
\mathbf w=(x_1,\ldots,x_k,y_1,\ldots,y_k)\in\R^{2k},
\]
and choose a partition $\sigma=(\sigma_L,\sigma_R)$ of
$\{1,\ldots,2k\}$ with $|\sigma_L|=|\sigma_R|=k$.
The corresponding coordinate vectors are
\[
x_{\sigma_L}=(w_j)_{j\in\sigma_L}\in\R^k,
\qquad
x_{\sigma_R}=(w_j)_{j\in\sigma_R}\in\R^k,
\]
where the indices in each group are arranged in increasing order.
Each measure $\nu_i$ occurs twice in
$\boldsymbol\nu\otimes\boldsymbol\nu$. Accordingly, set
\[
\widetilde\nu_j=
\begin{cases}
\nu_j,&1\le j\le k,\\
\nu_{j-k},&k+1\le j\le2k,
\end{cases}
\]
The product measures $\rho_{\sigma_L}=\bigotimes_{j\in\sigma_L}\widetilde\nu_j$
and $\rho_{\sigma_R}=\bigotimes_{j\in\sigma_R}\widetilde\nu_j$
are both on $\R^k$. Let $\pi_\sigma:\R^{2k}\to\R^k\times\R^k$ be the
coordinate permutation $(\mathbf x,\mathbf y)\mapsto(x_{\sigma_L},x_{\sigma_R})$.
Write $\Phi_\sigma=\Phi\circ\pi_\sigma^{-1}:\pi_\sigma(U\times U)\to\R$
for the phase in the reordered coordinates.

We assume the collision phase
$\Phi(\mathbf x,\mathbf y)=f(\mathbf x)-f(\mathbf y)$ satisfies
the Monge--Amp\`ere condition
\eqref{shell:bordered-condition}, that is, for some
partition $\sigma=(\sigma_L,\sigma_R)$ and some
$\mathbf z\in U$,
\begin{equation}\label{shell:local-bordered-condition}
\left.\det\begin{pmatrix}
D_{x_{\sigma_L}}D_{x_{\sigma_R}}\Phi_\sigma
& \nabla_{x_{\sigma_L}}\Phi_\sigma\\
(\nabla_{x_{\sigma_R}}\Phi_\sigma)^T & 0
\end{pmatrix}\right|_{(x_{\sigma_L},x_{\sigma_R})=\pi_\sigma(\mathbf z,\mathbf z)}
\ne0.
\end{equation}
The following proposition applies
Lemma~\ref{bi:fractal-fio-pairing} after localization
to a product neighborhood of $\mathbf z$.

\begin{samepage}
\begin{proposition}[]
\label{shell:balanced-local}
Let $k\ge3$ and let $f\in C^\infty(U;\R)$, where
$U\subset\R^k$ is open. Define $\Phi(\mathbf x,\mathbf y) =f(\mathbf x)-f(\mathbf y)$.
Suppose that \eqref{shell:local-bordered-condition} holds
for some $\mathbf z\in U$ and some partition
$\sigma=(\sigma_L,\sigma_R)$ of the $2k$ coordinates into two groups
of size $k$.
\\
Then there exists a product neighborhood
$Q=I_1\times\cdots\times I_k$ of $\mathbf z$
with the following property.
Let $\nu_i$ be an $s_i$-Frostman probability measure
on $\R$, supported in $I_i$, where
$0<s_i\le1$ for every $1\le i\le k$, and set $\omega=f_*(\nu_1\otimes\cdots\otimes\nu_k)$.
\\
Then
\begin{equation}\label{shell:balanced-bound}
\int_{R\le|\tau|\le2R}
|\widehat\omega(\tau)|^2\,\dd\tau
\lesssim R^{(k+1)/2-\sum_{i=1}^k s_i},
\qquad R\ge2.
\end{equation}
The implicit constant may depend on $f$, $Q$, and
the Frostman constants, but is independent of $R$.
In particular, if $\sum_{i=1}^k s_i>\frac{k+1}{2}$,
then $\omega$ has an $L^2(\R)$ density.
\end{proposition}
\end{samepage}

\begin{proof}
We first localize the phase near the nondegenerate point
$(\mathbf z,\mathbf z)$. Choose a nonnegative function
$\zeta\in C_c^\infty(\R\setminus\{0\})$
supported where $1/2<|\theta|<3$ and equal to one
on $1\le|\theta|\le2$.
By \eqref{eq:pairing-monge-ampere-scaling} with $d=k$,
the augmented mixed Hessian determinant of the joint phase
$\Psi(u,v,\theta)=\theta\Phi_\sigma(u,v)$, where
$u=x_{\sigma_L}$ and $v=x_{\sigma_R}$, is $\theta^{k-1}$
times its value at $(u,v,1)$.

Since $|\theta|$ is bounded away from zero on
$\supp\zeta$, condition  \eqref{shell:local-bordered-condition}
and continuity allow us to choose compact intervals
$I_i\subset\R$ containing $z_i$ in their interiors such that
$Q=I_1\times\cdots\times I_k\subset U$ and the determinant
is uniformly bounded away from zero whenever
$\mathbf x,\mathbf y\in Q$ and $\theta\in\supp\zeta$.
The same bound holds on a slightly larger neighborhood,
where both gradients $\nabla_{x_{\sigma_L}}\Phi_\sigma$
and $\nabla_{x_{\sigma_R}}\Phi_\sigma$ are also bounded away from zero.

To apply Lemma~\ref{bi:fractal-fio-pairing}, we restrict
the oscillatory integral operator to this nondegenerate region.
Choose $a\in C_c^\infty(\R^k\times\R^k)$ supported
in this neighborhood and equal to one on $\pi_\sigma(Q\times Q)$.
Multiply $\Phi_\sigma$ by a smooth cutoff compactly supported in
$\pi_\sigma(U\times U)$ and equal to one near $\supp a$, then extend
it by zero to $\R^k\times\R^k$. We retain $\Phi_\sigma$ for the extension.

Define $T_{R}:L^2(\R^k)\to L^2(\R^k)$ by
\[
T_{R}g(x_{\sigma_L})
=\int_{\R^k}
K_{T_{R}}(x_{\sigma_L},x_{\sigma_R})
g(x_{\sigma_R})\,\dd x_{\sigma_R},
\]
where
\[
K_{T_{R}}(x_{\sigma_L},x_{\sigma_R})
=R\,
\int_\R
e^{-2\pi iR\theta\Phi_\sigma(x_{\sigma_L},x_{\sigma_R})}
a(x_{\sigma_L},x_{\sigma_R})\zeta(\theta)\,\dd\theta.
\]
Since the measures $\nu_i$ are supported in $I_i$,
the cutoff $a$ equals one on $\supp(\rho_{\sigma_L}\otimes\rho_{\sigma_R})$.
Hence
\[
\begin{aligned}
&\int_\R\zeta(\tau/R)
|\widehat\omega(\tau)|^2\,\dd\tau=
\iint_{\R^k\times\R^k}
K_{T_{R}}(x_{\sigma_L},x_{\sigma_R})
\,\dd\rho_{\sigma_L}(x_{\sigma_L})
\,\dd\rho_{\sigma_R}(x_{\sigma_R}).
\end{aligned}
\]
On the amplitude support, $|\theta|\approx1$, the two group gradients of $\Phi_\sigma$ are bounded away from zero, and the collision determinant, with mixed derivative block $\theta D_{x_{\sigma_L}}D_{x_{\sigma_R}}\Phi_\sigma$, is bounded away from zero. These bounds verify
\eqref{eq:pairing-bordered-condition} and
\eqref{eq:pairing-gradient-condition}.
Moreover, the amplitude
$a(x_{\sigma_L},x_{\sigma_R})\zeta(\theta)$
is smooth and compactly supported, independently of $R$.
Thus, $T_{R}$ satisfies all the hypotheses of
Lemma~\ref{bi:fractal-fio-pairing}, with uniform constants.

It remains to record which Frostman exponents enter the pairing estimate. Define $\widetilde s_j=s_j$ for $1\le j\le k$ and
$\widetilde s_j=s_{j-k}$ for $k<j\le2k$. For $B=L,R$, let
$\alpha_{\sigma_B}=\sum_{j\in\sigma_B}\widetilde s_j$.
The product measures $\rho_{\sigma_L}$ and $\rho_{\sigma_R}$ are then
$\alpha_{\sigma_L}$- and $\alpha_{\sigma_R}$-Frostman, with
\[
\alpha_{\sigma_L}+\alpha_{\sigma_R}=2\sum_{i=1}^k s_i,
\]
since each measure $\nu_i$ occurs twice.
Applying Lemma~\ref{bi:fractal-fio-pairing}
with $d=k$ and $N>k$, we obtain
\[
\begin{aligned}
\left|
\iint K_{T_{R}}(x_{\sigma_L},x_{\sigma_R})
\,\dd\rho_{\sigma_L}(x_{\sigma_L})
\,\dd\rho_{\sigma_R}(x_{\sigma_R})
\right|\,
&\,\lesssim
R^{(k+1-\alpha_{\sigma_L}-\alpha_{\sigma_R})/2}
+O_N(R^{-N})\\
&\,\lesssim
R^{(k+1)/2-\sum_{i=1}^k s_i}.
\end{aligned}
\]
Since $\zeta\ge1$ on $1\le|\theta|\le2$,
this proves \eqref{shell:balanced-bound}.

If $\sum_{i=1}^k s_i>(k+1)/2$, then Lemma~\ref{lem:fourier-shell-consequences} applied to \eqref{shell:balanced-bound} with $\beta=\sum_{i=1}^k s_i-(k-1)/2>1$ gives an $L^2(\R)$ density. This completes the proof.
\end{proof}

    \section{\texorpdfstring{Dimension expansion for functions of three or more variables}{Dimension expansion for functions of three or more variables}}
\label{sec:analytic}

We prove Theorems~\ref{thm:multivariate-dimension-expansion}
and~\ref{thm:balanced-collision-image}. For the first, we relate images
of transversal families to incidence counts, then localize to a
three-variable slice in Subsection~\ref{analytic:proof-multivariate}.
The second uses the balanced Fourier shell estimate directly;
its proof is in Subsection~\ref{analytic:proof-balanced}.

\subsection{Dimension of images of transversal curve families}
\label{sec:families}

The discretized Furstenberg theorem controls the number
of squares in a union of shadings, whereas we seek the
Hausdorff dimension of the values attained by a family
of curves. To connect the two, we use one Hausdorff cover of the attained values for the whole family of curves. Pairing each covering
interval with horizontal intervals of the same length produces
squares meeting the graphs. These squares form the shadings
to which we apply the Furstenberg theorem. Comparing the
resulting lower bound with the upper count from the
common cover gives the estimate needed for
Theorem~\ref{thm:multivariate-dimension-expansion}.

\begin{samepage}
    \begin{lemma}\label{lem:common-image}
Let $0<t\leq2$ and $0<s\leq1$. Suppose that
$P\subset \R^2$ and $B\subset (-1,1)$ are compact
sets that support Frostman probability measures
$\mu$ and $\nu$, respectively, satisfying
\begin{equation*}
\mu(B(x,r))\leq C_\mu r^t,
\qquad
\nu(B(z,r))\leq C_\nu r^s
\quad\text{for all }x\in\R^2,\ z\in\R,\ r>0.
\end{equation*}
Suppose that $\mathcal H=\{h_p:p\in P\}\subset B_{C^2([-2,2])}(1)$
is a transversal family over $[-2,2]$, and that
the map $p\mapsto h_p$ is bi-Lipschitz with respect
to the $C^2$ metric. Assume also that $|h_p(z)|<1$ for $(p,z)\in P\times B$.
Then
\begin{equation*}
\dimH\{h_p(z):p\in P,\ z\in B\}
\geq
\min\left\{t,\frac{s+t}{2},1\right\}.
\end{equation*}
\end{lemma}
\end{samepage}
The key is to retain the cardinality of each shading in
Theorem~\ref{thm:OPY}. At a selected scale $\delta$, let
$\mathcal J_n$ be the intervals of length $\delta$ in a single
Hausdorff cover of the image,
and select horizontal intervals of comparable $\nu$-mass $w$.
If each selected curve carries $\nu$-mass at least $q$ over these
horizontal intervals with values in $\bigcup\mathcal J_n$, it meets
at least $q/(2w)$ squares. For a loss $\eta>0$, writing
$b=\min\{t,(s+t)/2,1\}$, the two counts give
\[
 \frac{q}{2w}\delta^{\eta-b}
 \leq\text{number of selected graph squares}
 \leq\frac{|\mathcal J_n|}{w}.
\]
The factor $w^{-1}$ cancels. The proof below arranges that $q$
loses only a power of $\log(1/\delta)$, so this comparison
rules out a Hausdorff cover with exponent below $b$.
No upper bound on packing dimension for the horizontal set is needed. Related scale-selection arguments appear in
\cite[Lemma~3.3]{HSY22}, \cite[Section~1.3]{OPY25}, and
\cite[Section~9]{RZ24}. For affine functions, the same
conclusion also follows from the sharp exceptional projection
theorem of Ren and Wang~\cite[Theorem~1.2]{RW25}.

\begin{proof}
Let
\[
D=\{h_p(z):p\in P,\ z\in B\},
\qquad
b=\min\left\{t,\frac{s+t}{2},1\right\}.
\]
Suppose, to the contrary, that $\dimH D<b$,
and choose $u,\eta>0$ such that $\dimH D<u<u+2\eta<b$.
Let $T\geq1$ be a transversality constant for
$\mathcal H$, and let $\varepsilon,\delta_0>0$
be given by Theorem~\ref{thm:OPY} with parameters
$s,t,T,\eta$.

Since $\dimH D<u$, for every sufficiently large integer
$n_0$, there exists a countable family $\mathcal J$
of dyadic intervals covering $D$, each of length at most
$2^{-n_0}$, such that
\begin{equation}\label{eq:common-image-cover-cost}
\sum_{J\in\mathcal J}|J|^u\leq1.
\end{equation}
Indeed, passing from ordinary intervals to dyadic intervals
changes the cost by at most a fixed factor, which can
be absorbed by choosing the original cover with sufficiently
small cost. We choose $n_0$ sufficiently large
to meet all the small-scale requirements below and fix
such a cover.
For each $n\geq n_0$, let $\mathcal J_n$
be the collection of distinct intervals in $\mathcal J$
of length $2^{-n}$, and set
\begin{equation*}
D_n=\bigcup_{J\in\mathcal J_n}J.
\end{equation*}
Since every interval in $\mathcal J_n$ has length
$2^{-n}$, \eqref{eq:common-image-cover-cost} gives $|\mathcal J_n|\leq 2^{nu}$.

For each $n\geq n_0$, let $E_n=\{(p,z)\in P\times B:h_p(z)\in D_n\}$.
The evaluation map $(p,z)\mapsto h_p(z)$ is continuous,
so each $E_n$ is Borel.
Integrating its indicator in $z$ gives a Borel function of
$p$, which also ensures measurability of the parameter
sets selected below.

Since the sets $D_n$ cover $D$, the sets $E_n$
cover $P\times B$. As $\mu$ and $\nu$ are
probability measures, choosing $n_0$ sufficiently
large that $\sum_{n\geq n_0}n^{-2}<1$, we find
some $n\geq n_0$ such that
\begin{equation}\label{eq:common-image-scale-mass}
(\mu\otimes\nu)(E_n)\geq n^{-2}.
\end{equation}
Fix such an $n$ and write $\delta=2^{-n}$.

Partition the $z$-axis into dyadic intervals of length
$\delta$. For each $j\geq0$, let $\mathcal I_j$
be the family of such intervals $I$ satisfying $2^{-j}\leq\nu(I)<2^{-j+1}$.
By Lemma~\ref{lem:dyadic-cell-masses} with $d=1$ and $\gamma=1$, the intervals meeting $B$ with $\nu(I)<\delta^2$ carry $O(\delta)$ total $\nu$-mass. The corresponding product mass of $E_n$ is no larger, since $\mu$ is a probability measure. For sufficiently large $n$, \eqref{eq:common-image-scale-mass} therefore leaves at least $n^{-2}/2$ of the product mass in the classes $\mathcal I_0,\ldots,\mathcal I_{2n}$.

By the pigeonhole principle, there exists
$j\in\{0,\ldots,2n\}$ such that, setting
$\mathcal I=\mathcal I_j$ and $w=2^{-j}$,
\begin{gather}
w\leq\nu(I)<2w
\quad\text{for all }I\in\mathcal I,
\\
(\mu\otimes\nu)
\left(E_n\cap\left(P\times\bigcup_{I\in\mathcal I}I\right)\right)
\geq c_1n^{-3}=:2q.
\label{eq:common-image-bin-mass}
\end{gather}
The products $I\times J$, with $I\in\mathcal I$
and $J\in\mathcal J_n$, form a family of dyadic
$\delta$-squares covering
$(\bigcup_{I\in\mathcal I}I)\times D_n$.
We use this same family of squares for every graph
in the subsequent argument.

The mass estimate \eqref{eq:common-image-bin-mass} averages over all curves, while
the Furstenberg theorem requires a substantial shading
on each selected curve. We therefore set
\[
P_*=
\left\{p\in P:
\nu\left\{z\in B\cap\bigcup_{I\in\mathcal I}I:
h_p(z)\in D_n\right\}\geq q
\right\}.
\]
By Fubini's theorem and \eqref{eq:common-image-bin-mass},
we have
\[
\begin{aligned}
2q
&\leq\int_P
\nu\left\{z\in B\cap\bigcup_{I\in\mathcal I}I:
h_p(z)\in D_n\right\}\,d\mu(p)\\
&\leq\mu(P_*)+q\mu(P\setminus P_*)\\
&\leq\mu(P_*)+q.
\end{aligned}
\]
Consequently,
\begin{equation}\label{eq:common-image-good-parameters}
\mu(P_*)\geq q\geq c_2n^{-3}.
\end{equation}

We now discretize the retained parameters. The restricted measure still satisfies the original nonconcentration bound.
Indeed, by \eqref{eq:common-image-good-parameters}, $\sigma=\mu|_{P_*}$ satisfies
\[
\sigma(\R^2)\geq q\geq c_2n^{-3},
\qquad
\sigma(B(x,r))\leq C_\mu r^t
\quad\text{for all }x\in\R^2,\ r>0.
\]
Let $Q_0$ be a fixed bounded rectangle containing $P$.
Applying Lemma~\ref{lem:discretization} with $E=P_*$,
ambient dimension $2$, Frostman exponent $t$, and $\kappa=3$,
we obtain a finite set $P_\delta\subset P_*$
such that
\[
|P_\delta|\geq c_3n^{-4}\delta^{-t},
\]
and $P_\delta$ is a $(\delta,t,Cn^4)$-set.
Since $p\mapsto h_p$ is bi-Lipschitz,
after enlarging $C$ by a factor depending only
on the bi-Lipschitz constants,
the family
$\mathcal H_\delta=\{h_p:p\in P_\delta\}$
is also a $(\delta,t,Cn^4)$-set in the $C^2$
metric and inherits the transversality constant $T$.

The comparable-mass horizontal intervals allow us
to construct a shading of controlled size and
nonconcentration on each selected curve. Indeed,
every $p\in P_\delta\subset P_*$
satisfies
\[
\nu\left\{z\in B\cap\bigcup_{I\in\mathcal I}I:
h_p(z)\in D_n\right\}\geq q.
\]
Since $\nu(I)<2w$ for every $I\in\mathcal I$,
there are at least $q/(2w)$ distinct intervals
$I\in\mathcal I$ containing some $z\in B$
with $h_p(z)\in D_n$.
For each such interval $I$, choose one such point $z$ and let
$J\in\mathcal J_n$ be the interval containing $h_p(z)$. Let
$\mathcal P(h_p)$ be the shading formed by the resulting squares $I\times J$.
These squares have distinct horizontal intervals,
so
\begin{equation}\label{eq:common-image-row-lower}
|\mathcal P(h_p)|\geq\frac{q}{2w}.
\end{equation}
Fix $\xi\in\R^2$ and $\delta\le r\le1$. The horizontal projections of the squares in $\mathcal P(h_p)$ meeting $B(\xi,r)$ are distinct dyadic intervals, each of $\nu$-mass at least $w$, and lie in an interval of length at most $2r+2\delta\le4r$. Thus the Frostman bound and \eqref{eq:common-image-row-lower} give
\[
\#\{Q\in\mathcal P(h_p):Q\cap B(\xi,r)\ne\varnothing\}
\lesssim \frac{r^s}{w}\lesssim q^{-1}r^s|\mathcal P(h_p)|
\lesssim n^3r^s|\mathcal P(h_p)|.
\]
Since the $\delta$-covering number of a union of dyadic $\delta$-squares is comparable to their number, $\mathcal P(h_p)$ is a $(\delta,s,Cn^3)$-set.

Set $M=\max\{1,\lceil q/(2w)\rceil\}$.
By \eqref{eq:common-image-row-lower},
\[
|\mathcal P(h_p)|\geq M,
\qquad\text{for all }p\in P_\delta.
\]
Since $P$ and $B$ are compact, and
$|h_p(z)|<1$ for all $(p,z)\in P\times B$,
the selected graph points lie a fixed positive
distance from the boundary of $[-1,1]^2$.
Thus, all selected squares are contained in
$[-1,1]^2$ for sufficiently large $n_0$.

By our choice of $n_0$, we also have
$Cn^4\leq\delta^{-\varepsilon}$ and
$\delta\leq\delta_0$ for every $n\geq n_0$.
We have constructed a finite transversal family
$\mathcal H_\delta$ with a shading $\mathcal P(h)$ for each
$h\in\mathcal H_\delta$, satisfying the required
nonconcentration estimates.
Theorem~\ref{thm:OPY} therefore gives
\begin{equation*}
\left|\bigcup_{p\in P_\delta}\mathcal P(h_p)\right|
\geq\delta^{\eta-b}M
\geq\frac{q}{2w}\delta^{\eta-b}.
\end{equation*}

For the upper bound, every selected square has the form $I\times J$ with $I\in\mathcal I$ and $J\in\mathcal J_n$. Since
$|\mathcal I|\leq w^{-1}$, we have
\begin{equation*}
\left|\bigcup_{p\in P_\delta}\mathcal P(h_p)\right|
\leq|\mathcal I||\mathcal J_n|
\leq\frac{|\mathcal J_n|}{w}.
\end{equation*}
Comparing the two counts cancels the factor $w^{-1}$ and gives
\[
|\mathcal J_n|
\geq\frac{q}{2}\delta^{\eta-b}
\geq c_4n^{-3}\delta^{\eta-b}.
\]
Since $\delta=2^{-n}$, for sufficiently large
$n_0$ we have $c_4n^{-3}>\delta^\eta$.
Consequently,
\[
|\mathcal J_n|
>\delta^{2\eta-b}
>\delta^{-u},
\]
where the last inequality follows from
$u+2\eta<b$. This contradicts
$|\mathcal J_n|\leq\delta^{-u}$ from \eqref{eq:common-image-cover-cost}. This completes the proof of the lemma.
\end{proof}

We now apply the common image estimate to prove
Theorem~\ref{thm:multivariate-dimension-expansion}. The key is to find a
transversal family of curves arising from the
analytic function $f$. We first establish a
differential criterion for analytic special form.
When this structure fails, analyticity allows us
to locate a nondegenerate point in the Cartesian
product of the input sets. The resulting local
family of curves satisfies the hypotheses of
Lemma~\ref{lem:common-image}, giving the desired
dimension bound.

\subsection{Characterizing analytic special form}

If $f=G(H_1+\cdots+H_k)$ on a regular box, then
$f_j/f_\ell=H_j'/H_\ell'$ is independent of every third
coordinate $x_i$. Thus the numerator of its $x_i$ derivative must vanish. We use its negative as a test for this structure. We use
the sign convention below because it also gives the determinant
of the corresponding value-and-slope map.

\begin{definition}
Let $k\geq3$, let $\Omega\subset\R^k$ be open, and let
$f\in C^2(\Omega)$. Write $f_i=\partial_{x_i}f$ and
$f_{ij}=\partial_{x_i}\partial_{x_j}f$, for $1\le i,j\le k$. For pairwise
distinct $i,j,\ell\in\{1,\ldots,k\}$, define the
\emph{auxiliary function}
\begin{equation}\label{eq:G-def}
 \Gaux_{i;j\ell}[f]:=f_jf_{i\ell}-f_\ell f_{ij}.
\end{equation}
When the function $f$ is fixed, we suppress $[f]$ from the notation.
\end{definition}
\begin{remark}
When $k=3$, the auxiliary functions \eqref{eq:G-def} agree,
up to sign and relabeling, with those of \cite[(5.1)]{Pha26}.
Note that
$\Gaux_{i;j\ell}=-\Gaux_{i;\ell j}$, and
\begin{equation*}
 \Gaux_{i;j\ell}
 =\Gaux_{j;i\ell}-\Gaux_{\ell;ij}.
\end{equation*}
Consequently, if two of these three auxiliary functions vanish
at a point, then so does the third.
\end{remark}
The next theorem shows that the simultaneous vanishing of
these auxiliary functions characterizes analytic special form
on a regular box. It
gives a $k$-variable version, on regular boxes, of
\cite[Lemma~5.3]{Pha26}.

\begin{samepage}
\begin{theorem}\label{thm:analytic-special-form}
Let $k\geq3$, let
$Q=I_1\times\cdots\times I_k\subset\R^k$
be an open box, and let $f\in C^\omega(Q)$
satisfy $f_1\cdots f_k\neq0$ on $Q$.
Then the following are equivalent:
\begin{itemize}
\item[(i)] $\Gaux_{i;j\ell}\equiv0$ on $Q$
for all pairwise distinct $i,j,\ell$ in $\{1,\ldots,k\}$.
\item[(ii)] There exist univariate real analytic
functions $H_1,\ldots,H_k,G$, where $H_j:I_j\to\R$ for $1\le j\le k$
and $G:I\to\R$ on an open interval $I$ containing $H_1(I_1)+\cdots+H_k(I_k)$, such that
\[
f(x_1,\ldots,x_k)
=G\bigl(H_1(x_1)+\cdots+H_k(x_k)\bigr),
\qquad x\in Q.
\]
\end{itemize}
\end{theorem}
\end{samepage}

\begin{proof}
Assume (i). We first show that the ratios of
first derivatives separate into one-variable
factors. These factors will then provide the
coordinate changes needed to recover the
special form.

Let $i,j,\ell$ be pairwise distinct indices in $\{1,\ldots,k\}$. Since
$f_j$ and $f_\ell$ do not vanish on $Q$,
the function
\[
L_{j\ell}:=\log\left|\frac{f_j}{f_\ell}\right|
\]
is real analytic on $Q$. Differentiating
with respect to $x_i$, we obtain
\[
\partial_iL_{j\ell}
=\frac{f_{ij}}{f_j}-\frac{f_{i\ell}}{f_\ell}
=-\frac{\Gaux_{i;j\ell}}{f_jf_\ell}=0.
\]
Thus, $L_{j\ell}$ depends only on $x_j$ and $x_\ell$.
We also write $L_{j\ell}(x_j,x_\ell)$ for the resulting function on $I_j\times I_\ell$.

Choose $m\notin\{j,\ell\}$ with $1\le m\le k$ and observe that $L_{j\ell}=L_{jm}-L_{\ell m}$.
Since $L_{jm}$ depends only on $(x_j,x_m)$
and $L_{\ell m}$ only on $(x_\ell,x_m)$, we have
$\partial_j\partial_\ell L_{j\ell}=0$ on $Q$.
As $Q$ is a product of intervals, integration
gives univariate real analytic functions
$u_{j\ell}$ and $w_{j\ell}$ such that
\begin{equation}\label{eq:R-separate}
L_{j\ell}(x_j,x_\ell)
=u_{j\ell}(x_j)+w_{j\ell}(x_\ell),
\qquad x\in Q.
\end{equation}

We next choose the one-variable functions consistently
across all pairs by fixing the first coordinate as a
reference. By \eqref{eq:R-separate}, we have $L_{j1}=u_{j1}(x_j)+w_{j1}(x_1), \quad j\ne1$.
If $j,\ell\ne1$ are distinct, then
\[
L_{j\ell}
=u_{j1}(x_j)-u_{\ell1}(x_\ell)
+w_{j1}(x_1)-w_{\ell1}(x_1).
\]
Since $L_{j\ell}$ is independent of $x_1$,
the difference $w_{j1}-w_{\ell1}$ must be
constant. Thus, there exist constants $c_j$,
with $c_2=0$ and $2\le j\le k$, such that
\[
w_{j1}=w_{21}+c_j,\qquad j\ne1.
\]
Setting $u_1=-w_{21}$ and
$u_j=u_{j1}+c_j$ for $j\ne1$, we obtain
\begin{equation}\label{eq:R-u}
L_{j\ell}(x_j,x_\ell)
=u_j(x_j)-u_\ell(x_\ell),
\qquad j\ne\ell.
\end{equation}

These one-variable functions now determine the
desired change of coordinates. Since each $f_j$
is nonvanishing on the connected box $Q$, its sign
$\varepsilon_j\in\{-1,1\}$ is constant. For each $1\le j\le k$,
choose a univariate real analytic function $H_j$ such that
$H_j'(x_j)=\varepsilon_j e^{u_j(x_j)}$.

Exponentiating \eqref{eq:R-u} and restoring signs gives
\begin{equation*}
\frac{f_j}{f_\ell}
=\frac{H_j'(x_j)}{H_\ell'(x_\ell)}
\quad\text{for all }j\ne\ell.
\end{equation*}
Set $a(x)=\frac{f_1(x)}{H_1'(x_1)}$, $x\in Q$. Then
\begin{equation}\label{eq:fj-lambda}
f_j=a H_j'(x_j)
\quad\text{for all }j=1,\ldots,k.
\end{equation}
Set
$S(x):=H_1(x_1)+\cdots+H_k(x_k)$.
Then \eqref{eq:fj-lambda} yields
\begin{equation*}
df=\sum_{j=1}^k f_j\,dx_j
=a\sum_{j=1}^k H_j'(x_j)\,dx_j
=a\,dS.
\end{equation*}

Since each $H_j'$ is nonvanishing on $I_j$,
the function $H_j$ is an analytic diffeomorphism
from $I_j$ onto the interval $J_j=H_j(I_j)$.
Consequently, the map
\[
\Psi:Q\longrightarrow J_1\times\cdots\times J_k,
\qquad
\Psi(x)=\bigl(H_1(x_1),\ldots,H_k(x_k)\bigr),
\]
is an analytic diffeomorphism. In the new
coordinates $z_j=H_j(x_j)$ for $1\le j\le k$, define
\[
F=f\circ\Psi^{-1},
\qquad
\widetilde a=a\circ\Psi^{-1}.
\]
By \eqref{eq:fj-lambda},
\[
\partial_{z_j}F
=\left(\frac{f_j}{H_j'}\right)\circ\Psi^{-1}
=\widetilde a,
\qquad j=1,\ldots,k.
\]

The equality of these partial derivatives shows
that $F$ depends only on the sum of its coordinates.
Indeed, let
\[
T(z)=z_1+\cdots+z_k,
\qquad
I=T(J_1\times\cdots\times J_k).
\]
Then $dF=\widetilde a\,dT$, so $F$ is
constant along each level set of $T$ inside
the box. To see this, note that these intersections
are convex. Hence, whenever $T(z)=T(z')$,
\[
\begin{aligned}
F(z')-F(z)
&=\int_0^1
\widetilde a(z+v(z'-z))
\sum_{j=1}^k(z_j'-z_j)\,dv=0.
\end{aligned}
\]
Thus, there is a well-defined function
$G:I\to\R$ such that $F=G\circ T$.

It remains to verify that $G$ is analytic.
Fix $t_0\in I$ and choose
$z^0\in J_1\times\cdots\times J_k$
with $T(z^0)=t_0$. For $t$ sufficiently
close to $t_0$, we have $G(t)=F(z_1^0+t-t_0,z_2^0,\ldots,z_k^0)$.
The right-hand side is analytic in $t$.
Since these local expressions agree on
overlapping intervals, $G$ is analytic on $I$.
Returning to the original coordinates, we obtain
\[
f(x)
=G\bigl(H_1(x_1)+\cdots+H_k(x_k)\bigr),
\qquad x\in Q,
\]
which proves (ii).

Conversely, suppose that (ii) holds, and write
$S=H_1(x_1)+\cdots+H_k(x_k)$. Then
\[
f_j=G'(S)H_j'(x_j),
\qquad
f_{ij}=G''(S)H_i'(x_i)H_j'(x_j),
\qquad i\ne j.
\]
Consequently, whenever $i,j,\ell$ are
pairwise distinct indices in $\{1,\ldots,k\}$,
\[
\begin{aligned}
\Gaux_{i;j\ell}
&=G'(S)G''(S)H_i'
\bigl(H_j'H_\ell'-H_\ell'H_j'\bigr)=0.
\end{aligned}
\]
This proves (i) and completes the proof.
\end{proof}

\subsection{Local dimension expansion}

The preceding theorem characterizes analytic special form
through the vanishing of the auxiliary functions.
In three variables, after relabeling so that $\Gaux_{3;12}(x^0)\ne0$,
the nonvanishing has a useful geometric consequence near $x^0$:
the values and slopes of the curves
$z\mapsto f(x_1,x_2,z)$ vary independently with
the parameters $(x_1,x_2)$. The following lemma
shows that this nondegeneracy yields a transversal
family of curves, allowing us to apply
Lemma~\ref{lem:common-image}.

\begin{lemma}\label{lem:local-transverse}
Let $U\subset\R^3$ be open and $f\in C^3(U)$.
Suppose that there exists
$x^0=(x_1^0,x_2^0,x_3^0)\in U$ such that
\begin{equation*}
\Gaux_{3;12}(x^0)\neq0.
\end{equation*}
Then there exist a compact rectangle
$Q\subset\R^2$ containing $(x_1^0,x_2^0)$
in its interior and a compact interval
$I\subset\R$ containing $x_3^0$ in its
interior, with $Q\times I\subset U$, such that
the family
\[
\mathcal H=\{h_p:p\in Q\},
\qquad
h_p(z)=f(p,z),
\]
is transversal over $I$.
Moreover, the map $p\mapsto h_p$ is
bi-Lipschitz from $Q$ into $C^2(I)$.
\end{lemma}

\begin{proof}
Write $p^0=(x_1^0,x_2^0)$ and $z_0=x_3^0$.
To quantify the variation of the curves'
values and slopes with their parameters, consider
the map from parameters to first jets $\mathcal J_z:U_z\to\R^2$, where
$U_z=\{p\in\R^2:(p,z)\in U\}$, given by $\mathcal J_z(p)=\bigl(f(p,z),f_3(p,z)\bigr)$.
By symmetry of mixed second derivatives, the parameter Jacobian satisfies
\[
\det D_p\mathcal J_{z_0}(p^0)
=(f_1f_{23}-f_2f_{13})(x^0)
=\Gaux_{3;12}(x^0)\ne0.
\]
Thus, $D_p\mathcal J_{z_0}(p^0)$ is invertible,
so there exists $m>0$ such that
\[
|D_p\mathcal J_{z_0}(p^0)(v)|\geq m|v|,
\qquad\text{for all }v\in\R^2.
\]
Since $(p,z)\mapsto D_p\mathcal J_z(p)$ is continuous,
we may choose a sufficiently small compact rectangle
$Q$ centered at $p^0$ and a compact interval $I$
centered at $z_0$, with $Q\times I\subset U$, such that
\[
\|D_p\mathcal J_z(p)-D_p\mathcal J_{z_0}(p^0)\|_{\mathrm{op}}
\leq\frac{m}{2}
\qquad\text{for all }(p,z)\in Q\times I.
\]
For $z\in I$, let $A_z:C^2(I)\to\R^2$ be given by
$A_z(h)=(h(z),h'(z))$. Then $\mathcal J_z(p)=A_z(h_p)$ for $p\in Q$.
Since $Q$ is convex, for any $p,p'\in Q$ and $z\in I$,
\[
\begin{aligned}
\mathcal J_z(p)-\mathcal J_z(p')
&=D_p\mathcal J_{z_0}(p^0)(p-p')\\
&\quad+\int_0^1
\bigl(D_p\mathcal J_z(p'+t(p-p'))
-D_p\mathcal J_{z_0}(p^0)\bigr)(p-p')\,dt.
\end{aligned}
\]
Consequently,
\begin{equation}\label{eq:local-jet-lower}
|\mathcal J_z(p)-\mathcal J_z(p')|
\geq\frac{m}{2}|p-p'|,
\qquad z\in I.
\end{equation}
On the other hand, since $f\in C^3(U)$ and
$Q\times I$ is compact, there exists $L>0$ such that
\begin{equation}\label{eq:local-c2-upper}
\|h_p-h_{p'}\|_{C^2(I)}
\leq L|p-p'|
\qquad\text{for all }p,p'\in Q.
\end{equation}
By \eqref{eq:local-jet-lower} and \eqref{eq:local-c2-upper}, for every $p,p'\in Q$ and $z\in I$,
\[
\frac m2|p-p'|
\le |h_p(z)-h_{p'}(z)|+|h_p'(z)-h_{p'}'(z)|
\le \|h_p-h_{p'}\|_{C^2(I)}\le L|p-p'|.
\]
Thus $p\mapsto h_p$ is bi-Lipschitz, and \eqref{eq:transversality} holds with $T=\max\{1,2L/m\}$. This completes the proof of the lemma.
\end{proof}

Once Lemma~\ref{lem:local-transverse} supplies the transversal family, Lemma~\ref{lem:common-image} gives the following local dimension estimate.

\begin{corollary}\label{cor:trivariate-local-dimension-bound}
Under the hypotheses of Lemma~\ref{lem:local-transverse},
there exist open intervals $V_i$ containing $x_i^0$ for $1\le i\le3$,
with $V_1\times V_2\times V_3\subset U$, such that the following holds.
Suppose that $A_i\subset V_i$ are compact sets
supporting $\alpha_i$-Frostman probability measures,
where $\alpha_i\in(0,1]$, for $1\le i\le3$. Then
\begin{equation}\label{eq:local-product-image}
\dimH f(A_1\times A_2\times A_3)
\geq
\min\left\{
\alpha_1+\alpha_2,
\frac{\alpha_1+\alpha_2+\alpha_3}{2},
1
\right\}.
\end{equation}
\end{corollary}

\begin{proof}
Lemma~\ref{lem:local-transverse} provides a transversal
family of curves near $x^0$. To apply
Lemma~\ref{lem:common-image}, we normalize the
domain and the $C^2$ norm of this family while
preserving its transversality.

Let $Q$ and $I$ be given by
Lemma~\ref{lem:local-transverse}.
Set $z_0=x_3^0$. Shrink $I$ about $z_0$ so that $I=[z_0-2r,z_0+2r]$
for some $r>0$. The upper $C^2$ bound is preserved by restriction,
and \eqref{eq:local-jet-lower} still holds for every $z$ in the smaller interval.
Choose sufficiently small open intervals
$V_1,V_2,V_3$ containing $x_1^0,x_2^0,x_3^0$,
respectively, such that
\[
V_1\times V_2\subset Q,
\qquad
V_3\subset(z_0-r,z_0+r).
\]

Rescale $I$ to $[-2,2]$ by setting
$z=z_0+r\theta$. Since $f\in C^3(U)$
and $Q\times I$ is compact, there exists
$M<\infty$ such that
\[
\sup_{p\in Q}
\|f(p,z_0+r\,\cdot)\|_{C^2([-2,2])}\leq M.
\]
Choose $\rho>0$ sufficiently small that
$\rho M<1$, and define
\[
\widetilde h_p(\theta)
=\rho f(p,z_0+r\theta),
\qquad p\in Q,\quad\theta\in[-2,2].
\]
Then
\[
\|\widetilde h_p\|_{C^2([-2,2])}<1,
\qquad p\in Q.
\]
For $p,p'\in Q$, put $\widetilde g=\widetilde h_p-\widetilde h_{p'}$.
The two components of the first jet are multiplied by $\rho$ and
$\rho r$, respectively. Thus \eqref{eq:local-jet-lower} and
\eqref{eq:local-c2-upper} give, for every $\theta\in[-2,2]$,
\[
\begin{aligned}
\frac{\rho m}{2}\min\{1,r\}|p-p'|
&\le |\widetilde g(\theta)|+|\widetilde g'(\theta)| \\
&\le \|\widetilde g\|_{C^2([-2,2])}
\le \rho L\max\{1,r,r^2\}|p-p'|,
\end{aligned}
\]
where $m,L$ are the constants in those estimates. The outer bounds
make $p\mapsto\widetilde h_p$ bi-Lipschitz. The pointwise lower bound
and the norm upper bound give transversality.

Now set
\[
P=A_1\times A_2,
\qquad
\widetilde A_3
=\left\{\frac{z-z_0}{r}:z\in A_3\right\}
\subset(-1,1).
\]
Let $\mu_i$ be an $\alpha_i$-Frostman
probability measure supported on $A_i$,
for $i=1,2,3$. Then $\mu_1\otimes\mu_2$
is an $(\alpha_1+\alpha_2)$-Frostman
probability measure on $P$.
Since $z\mapsto(z-z_0)/r$ is bi-Lipschitz,
the pushforward of $\mu_3$ under this map
is an $\alpha_3$-Frostman probability measure
on $\widetilde A_3$.

Applying Lemma~\ref{lem:common-image} with
$t=\alpha_1+\alpha_2$ and $s=\alpha_3$,
we obtain
\[
\dimH\{\widetilde h_p(\theta):
p\in P,\ \theta\in\widetilde A_3\}
\geq
\min\left\{
\alpha_1+\alpha_2,
\frac{\alpha_1+\alpha_2+\alpha_3}{2},
1
\right\}.
\]
Finally,
\[
\{\widetilde h_p(\theta):
p\in P,\ \theta\in\widetilde A_3\}
=\rho f(A_1\times A_2\times A_3).
\]
Since multiplication by $\rho>0$ preserves
Hausdorff dimension, this proves
\eqref{eq:local-product-image}.
{This finishes the proof of the corollary.}
\end{proof}

\subsection{Proof of Theorem~\ref{thm:multivariate-dimension-expansion}}
\label{analytic:proof-multivariate}

The proof of Theorem~\ref{thm:multivariate-dimension-expansion}
reduces the problem to three variables. The key is to use
analyticity to find a point in the Cartesian product of the
input sets where the required differential condition is nonzero.
We then fix the remaining variables and apply the trivariate
estimates.

\begin{proof}[Proof of Theorem~\ref{thm:multivariate-dimension-expansion}]
If $f$ is of analytic special form, alternative (i) holds. Otherwise,
since no $f_i$ vanishes identically on the connected domain $\Omega$,
the sets $\{f_i\ne0\}$ are open and dense in $\Omega$. Their intersection
contains an open product box. By Theorem~\ref{thm:analytic-special-form},
not all the auxiliary functions in \eqref{eq:G-def} vanish identically
on this box. Relabel the coordinates so that $\Gaux_{3;12}\not\equiv0$.

Let $A\subset[0,1]$ be Borel with $\dimH A=\alpha>0$, and fix
$0<\beta<\alpha$. By Frostman's lemma, choose a $\beta$-Frostman
probability $\mu_A$ with compact support $K\subset A$. The support $K$
is infinite, so Lemma~\ref{lem:analytic-product-uniqueness} gives
$\mathbf a=(a_1,\ldots,a_k)\in K^k$ with $\Gaux_{3;12}(\mathbf a)\ne0$.
Fix the last $k-3$ coordinates at $a_4,\ldots,a_k$ and set
\[
 V=\{\mathbf x\in\R^3:(\mathbf x,a_4,\ldots,a_k)\in\Omega\},
 \qquad F(\mathbf x)=f(\mathbf x,a_4,\ldots,a_k)\quad(\mathbf x\in V).
\]
For $k=3$, the fixed tuple is empty and $F=f$. At
$\mathbf a'=(a_1,a_2,a_3)$, equality of mixed partials gives
$(F_1F_{23}-F_2F_{13})(\mathbf a')=\Gaux_{3;12}[F](\mathbf a')
=\Gaux_{3;12}[f](\mathbf a)\ne0$,
so Corollary~\ref{cor:trivariate-local-dimension-bound} applies.

For $\Phi_F(\mathbf x,\mathbf y)=F(\mathbf x)-F(\mathbf y)$, use the
partition $\sigma_L=\{1,2,6\}$, $\sigma_R=\{3,4,5\}$, with coordinates
ordered as $(x_1,x_2,y_3)$ and $(y_1,y_2,x_3)$, respectively. The right
order is an even permutation of increasing index order. Let
$\mathcal D_F:V\times V\to\R$ denote the Monge--Amp\`ere determinant
for $\Phi_F$ in this ordered partition. At $(\mathbf z,\mathbf z)$
with $\mathbf z\in V$, its matrix, with rows ordered as $(x_1,x_2,y_3,*)$ and columns ordered as $(y_1,y_2,x_3,*)$, where $*$ denotes the border, is
\[
 \begin{pmatrix}0_{2\times2}&B\\-B^T&C\end{pmatrix},\qquad
 B=\begin{pmatrix}F_{13}&F_1\\F_{23}&F_2\end{pmatrix},\qquad
 C=\begin{pmatrix}0&-F_3\\F_3&0\end{pmatrix},
\]
with derivatives evaluated at $\mathbf z$. Moving the last two columns
to the front is even and gives diagonal blocks $B$ and $-B^T$. Thus the
Monge--Amp\`ere determinant satisfies
\begin{equation}\label{analytic:trivariate-balanced-determinant}
 \mathcal D_F(\mathbf z,\mathbf z)
 =(F_1F_{23}-F_2F_{13})(\mathbf z)^2.
\end{equation}
It is nonzero at $(\mathbf a',\mathbf a')$, so
Proposition~\ref{shell:balanced-local} applies there as well.

For every $1\le i\le3$, choose a compact interval $J_i$ with $a_i$
in its interior so that $J_1\times J_2\times J_3$ lies in both local
neighborhoods. Since $a_i\in K=\supp\mu_A$, the normalized restriction
$\mu_{A,i}=\mu_A|_{J_i}/\mu_A(J_i)$ is a $\beta$-Frostman probability
with compact support $B_i\subset A$. The fixed coordinates belong to
$A$, so Corollary~\ref{cor:trivariate-local-dimension-bound} gives
\[
 \dimH f(A^k)\ge\dimH F(B_1\times B_2\times B_3)
 \ge\min\{2\beta,3\beta/2,1\}=\min\{3\beta/2,1\}.
\]
Letting $\beta\uparrow\alpha$ proves \eqref{eq:main-bound}.
If $\alpha>2/3$, choose $2/3<\beta<\alpha$ in the same construction.
Since $3\beta>2$, Proposition~\ref{shell:balanced-local} gives an
$L^2(\R)$ density for $F_*(\mu_{A,1}\otimes\mu_{A,2}\otimes\mu_{A,3})$.
This probability is carried by the compact set
$F(B_1\times B_2\times B_3)\subset f(A^k)$, so $\Leb(f(A^k))>0$.
This recovers the positive-measure threshold of
\cite[Theorem~5.1(ii)]{Pha26} and completes the proof.
\end{proof}

\subsection{Proof of Theorem~\ref{thm:balanced-collision-image}}
\label{analytic:proof-balanced}

The local shell estimate in Proposition~\ref{shell:balanced-local} gives the following product-set and product-measure version of Theorem~\ref{thm:balanced-collision-image}.
\begin{samepage}
\begin{theorem}
\label{thm:balanced-collision-products}
Let $k\ge3$, let $\Omega\subset\R^k$ be connected and open,
and let $f\in C^\omega(\Omega)$.
Fix compact intervals $J_i\subset\R$ for $1\le i\le k$,
with $Q_0=J_1\times\cdots\times J_k\subset\Omega$.
Define $\Phi(\mathbf x,\mathbf y)=f(\mathbf x)-f(\mathbf y)$.
Suppose that there exist a partition
$\sigma=(\sigma_L,\sigma_R)$ of the $2k$ scalar coordinates
into two groups of size $k$ and a point $\mathbf z_0\in\Omega$
such that the Monge--Amp\`ere condition
\begin{equation}\label{analytic:balanced-product-condition}
\left.
\det\begin{pmatrix}
D_{x_{\sigma_L}}D_{x_{\sigma_R}}\Phi_\sigma
& \nabla_{x_{\sigma_L}}\Phi_\sigma\\
(\nabla_{x_{\sigma_R}}\Phi_\sigma)^T & 0
\end{pmatrix}
\right|_{(x_{\sigma_L},x_{\sigma_R})=\pi_\sigma(\mathbf z_0,\mathbf z_0)}
\ne0
\end{equation}
holds. Then the following conclusions hold.
\\
\begin{enumerate}[label=(\roman*)]
\item Let $\nu_i$ be an $s_i$-Frostman probability measure
supported on $J_i$, where $0<s_i\le1$ for every $1\le i\le k$.
If $\sum_{i=1}^k s_i>\frac{k+1}{2}$,
then $f_*\left(\bigotimes_{i=1}^k\nu_i\right)\ll\Leb$.

\item Let $A_i\subset J_i$ be Borel sets with
$\alpha_i=\dimH A_i>0$ for every $1\le i\le k$.
Then
\[
\dimH f(A_1\times\cdots\times A_k)
\ge
\max\left\{0,
\min\left\{1,\sum_{i=1}^k\alpha_i-\frac{k-1}{2}\right\}
\right\}.
\]
Moreover, if $\sum_{i=1}^k\alpha_i>\frac{k+1}{2}$,
then $\Leb(f(A_1\times\cdots\times A_k))>0$.
\end{enumerate}
\end{theorem}
\end{samepage}

\begin{proof}
Fix the partition $\sigma=(\sigma_L,\sigma_R)$ from the theorem and
let $h:\Omega\to\R$ be the determinant in
\eqref{analytic:balanced-product-condition} evaluated at
$(x_{\sigma_L},x_{\sigma_R})=\pi_\sigma(\mathbf z,\mathbf z)$, i.e.
\[
h(\mathbf z)=
\left.
\det\begin{pmatrix}
D_{x_{\sigma_L}}D_{x_{\sigma_R}}\Phi_\sigma
&\nabla_{x_{\sigma_L}}\Phi_\sigma\\
(\nabla_{x_{\sigma_R}}\Phi_\sigma)^T&0
\end{pmatrix}
\right|_{(x_{\sigma_L},x_{\sigma_R})=\pi_\sigma(\mathbf z,\mathbf z)}.
\]
Then $h$ is real analytic
and $h(\mathbf z_0)\ne0$. For $s_i$-Frostman probability measures $\nu_i$
supported on $J_i$, with $s_i>0$ for every $1\le i\le k$, set
$\boldsymbol\nu=\bigotimes_{i=1}^k\nu_i$. The measures $\nu_i$ are nonatomic,
so Lemma~\ref{lem:analytic-product-zero-mass} gives
\begin{equation}\label{analytic:balanced-zero-set}
\boldsymbol\nu\bigl(Q_0\cap\{h=0\}\bigr)=0.
\end{equation}

Cover $Q_0\cap\{h\ne0\}$ by the interiors of countably many compact
product rectangles $Q_\ell=\prod_{i=1}^k I_{i,\ell}$, each contained in a
product neighborhood supplied by Proposition~\ref{shell:balanced-local}.
Whenever $\boldsymbol\nu(Q_\ell)>0$, the normalized restriction
$\boldsymbol\nu|_{Q_\ell}/\boldsymbol\nu(Q_\ell)$ is the product of the
normalized marginal restrictions $\nu_i|_{I_{i,\ell}}/\nu_i(I_{i,\ell})$.
They retain their exponents $s_i$. Thus its pushforward $\omega_\ell$
under $f$ satisfies
\begin{equation}\label{analytic:balanced-local-shell}
\int_{R\le|\tau|\le2R}|\widehat{\omega_\ell}(\tau)|^2\,\dd\tau
\lesssim_\ell R^{(k+1)/2-\sum_{i=1}^k s_i},\qquad R\ge2.
\end{equation}

For part (i), $\sum_{i=1}^k s_i>(k+1)/2$, so
Lemma~\ref{lem:fourier-shell-consequences} gives an $L^2(\R)$ density
for each $\omega_\ell$. The positive-mass rectangles cover
$\boldsymbol\nu$-almost every point by \eqref{analytic:balanced-zero-set}.
Apply Lemma~\ref{lem:local-to-global-ac} to $\boldsymbol\nu$ and the
rectangle indicators: $f_*(\boldsymbol\nu|_{Q_\ell})=
\boldsymbol\nu(Q_\ell)\omega_\ell\ll\Leb$, so $f_*\boldsymbol\nu\ll\Leb$.

For part (ii), let $\alpha_i=\dimH A_i>0$ for every $1\le i\le k$.
If $\sum_{i=1}^k\alpha_i\le(k-1)/2$, the dimension bound is trivial.
Otherwise choose $0<s_i<\alpha_i$ with $\sum_{i=1}^k s_i>(k-1)/2$.
If $\sum_{i=1}^k\alpha_i>(k+1)/2$, take the $s_i$ close enough to
$\alpha_i$ that $\sum_{i=1}^k s_i>(k+1)/2$ as well.
By Frostman's lemma, choose $s_i$-Frostman probabilities supported on
compact sets $K_i\subset A_i$. For at least one rectangle $Q_\ell$ the
product mass is positive. Its $\omega_\ell$ is carried by the compact
set $f(K_1\times\cdots\times K_k)$. From
\eqref{analytic:balanced-local-shell} and
Lemma~\ref{lem:fourier-shell-consequences},
\[
\dimH f(A_1\times\cdots\times A_k)
\ge \dimH\supp\omega_\ell
\ge \min\left\{1,\sum_{i=1}^k s_i-\frac{k-1}{2}\right\}.
\]
Taking the supremum over the admissible exponents gives the
dimension bound. If $\sum_{i=1}^k\alpha_i>(k+1)/2$, fix a choice with
$\sum_{i=1}^k s_i>(k+1)/2$. The corresponding $\omega_\ell$ has an
$L^2(\R)$ density, so its compact carrier has positive Lebesgue measure.
This completes the proof.
\end{proof}

\section{\texorpdfstring{Dimension expansion for functions of two variables}{Dimension expansion for functions of two variables}}
\label{sec:bivariate}

In this section, we prove Theorem~\ref{thm:bivariate-dimension-expansion}
by studying pairs of input points whose images under
$f$ are nearly equal. After localization and rescaling,
the collision relation gives rise to transversal families
of curves. We use incidence estimates to control the
corresponding curve-square incidences, while FIO
estimates exploit oscillatory cancellation in the
same geometry. Together, these approaches yield
\eqref{eq:bivariate-bound}.
Subsection~\ref{bi:proof-completion} combines the two bounds and
proves the positive-measure conclusion.

Throughout this section, we assume that $f$ satisfies
the hypotheses of Theorem~\ref{thm:bivariate-dimension-expansion} and is
nonspecial, since otherwise alternative (i) holds.
Implicit constants may depend on the fixed local
geometry, smooth cutoffs, Frostman constants, and
exponents, but not on the separation scales or
frequencies. The discrete incidence estimates are
applied at sufficiently small dyadic scales, while
the measure estimates hold at arbitrary scales
between zero and one.

We use the localization and normalized collision curves from
Subsection~\ref{bi:localization-section}. On the localized intervals, the
next proposition bounds the probability that two independent samples from
the image measure are within distance $\delta$.

\begin{proposition}\label{bi:collision-bound}
Let $0<s<1$. The intervals $I_1,I_2$ in
Lemma~\ref{bi:localization} may be chosen sufficiently small
so that the following holds. For any $s$-Frostman probability
measures $\nu,\mu$ supported on $I_1,I_2$, respectively,
every $\eta>0$, and every $0<\delta\leq1$,
\begin{equation}\label{bi:collision-estimate}
(\nu\otimes\mu)^{\otimes2}
\bigl\{|f(x,y)-f(x',y')|\leq\delta\bigr\}
\lesssim_\eta
\delta^{\min\{4s/3,\,2s-s^2\}-\eta}.
\end{equation}
The implicit constant may depend on $s,\eta$, the fixed local
geometry, and the Frostman constants of $\nu,\mu$, but is
independent of $\delta$.

Consequently, the image measure $\omega=f_*(\nu\otimes\mu)$
has finite $u$-energy for every
$0<u<\min\{4s/3,2s-s^2\}$.
\end{proposition}

We prove this proposition in Subsections~\ref{bi:incidences}--\ref{bi:collisions} and use it in the proof of
Theorem~\ref{thm:bivariate-dimension-expansion} in
Subsection~\ref{bi:proof-completion}.

\subsection{Incidence bounds for collision counting}
\label{bi:incidences}

We derive from the rectangular incidence estimate of
O'Regan, Wu, and Yi
\cite[Theorem~4.4]{OWY26} a bound adapted to
the collision problem. The key point is to retain
the dependence on the nonconcentration constants,
since these constants later depend on the masses of the selected cells.
For $s>2/3$, we obtain a stronger estimate by keeping
track of the cardinality of each shading.

\begin{lemma}\label{bi:density-incidence}
Let $0<s<1$ and $\vartheta=\max\{2/3,s\}$.
Suppose that $\mathcal F$ is a $\delta$-separated
transversal family that is rectangular
$(\delta,2s,C_1)$-Katz--Tao, and that every
shading $\mathcal P(h)$ is $(\delta,s,K)$-Katz--Tao,
where $C_1,K\geq1$.
Assume that
$\mathcal F\subset B_{C^2([-2,2])}(1)$ is finite,
with a fixed transversality constant, and that all
shading squares are contained in $[-1,1]^2$.
Then, for every $\eta>0$,
\begin{equation}\label{bi:density-incidence-bound}
\mathcal I(\mathcal F,\mathcal P(\cdot))
\lesssim_\eta
\delta^{-s\vartheta-\eta}
C_1^\vartheta K^{1-\vartheta}
|\mathcal F|^{1-\vartheta}|\mathcal P|^\vartheta.
\end{equation}
The estimate holds for all sufficiently small
dyadic $\delta>0$. The implicit constant and
scale threshold are independent of $C_1,K$.
\end{lemma}

\begin{proof}
Set $\sigma=\min\{2s,2-2s\}$.
If $s\leq2/3$, then $\sigma\geq s$.
Since $r/\delta\geq1$ for $\delta\leq r\leq1$,
the $(\delta,s,K)$-Katz--Tao condition implies
the $(\delta,\sigma,K)$-Katz--Tao condition.
Applying Theorem~\ref{bi:owy-rectangular}
with $C_2=K$, we obtain
\[
\mathcal I(\mathcal F,\mathcal P(\cdot))
\lesssim_\eta
\delta^{-2s/3-\eta}
C_1^{2/3}K^{1/3}
|\mathcal F|^{1/3}|\mathcal P|^{2/3},
\]
which is \eqref{bi:density-incidence-bound}
with $\vartheta=2/3$.

Suppose now that $s>2/3$, so that
$\sigma=2-2s<s$ and $\vartheta=s$.
The improvement comes from retaining the
cardinality of each shading when applying
the rectangular incidence estimate.

We first group the curves according to their
shading cardinalities. Since the curves have
uniformly bounded slopes on a fixed interval,
each graph meets $O(\delta^{-1})$ dyadic
$\delta$-squares. Thus, every nonempty shading
satisfies $1\leq|\mathcal P(h)|\lesssim\delta^{-1}$.
Partition the curves with nonempty shadings
into subfamilies $\mathcal F_N$ such that
\[
N\leq|\mathcal P(h)|<2N,
\qquad h\in\mathcal F_N,
\]
where $N$ ranges over the dyadic integers for which $\mathcal F_N\ne\varnothing$.
There are $O(\log(1/\delta))$ such classes.

Write
\[
\mathcal P_N
=\bigcup_{h\in\mathcal F_N}\mathcal P(h),
\qquad
\mathcal I_N
=\sum_{h\in\mathcal F_N}|\mathcal P(h)|.
\]
By construction,
\[
\mathcal I(\mathcal F,\mathcal P(\cdot))
=\sum_N\mathcal I_N,
\qquad
N|\mathcal F_N|
\leq\mathcal I_N<2N|\mathcal F_N|.
\]

For every $h\in\mathcal F_N$,
$z\in\R^2$, and $\delta\leq r\leq1$,
the cardinality bound and the Katz--Tao
condition give
\[
|\mathcal P(h)\cap B(z,r)|
\leq\min\left\{2N,K(r/\delta)^s\right\}.
\]
Since $0<\sigma<s$, interpolating these
two bounds with exponent $\sigma/s$ yields
\[
|\mathcal P(h)\cap B(z,r)|
\leq(2N)^{1-\sigma/s}
K^{\sigma/s}(r/\delta)^\sigma.
\]
Thus, each shading in this class is
$(\delta,\sigma,C_2)$-Katz--Tao with
$C_2=(2N)^{1-\sigma/s}K^{\sigma/s}$.
Applying Theorem~\ref{bi:owy-rectangular}
to $\mathcal F_N$, we obtain, for every
$\varepsilon>0$,
\[
\begin{aligned}
N|\mathcal F_N|
&\leq\mathcal I_N\\
&\lesssim_\varepsilon
\delta^{-2s/3-\varepsilon}
C_1^{2/3}C_2^{1/3}
|\mathcal F_N|^{1/3}|\mathcal P_N|^{2/3}\\
&\lesssim_\varepsilon
\delta^{-2s/3-\varepsilon}
C_1^{2/3}K^{\sigma/(3s)}
N^{(1-\sigma/s)/3}
|\mathcal F_N|^{1/3}|\mathcal P_N|^{2/3}.
\end{aligned}
\]
Since $\sigma=2-2s$, we have
\[
1-\frac{1-\sigma/s}{3}
=\frac{2s+\sigma}{3s}
=\frac{2}{3s}.
\]
Collecting the powers of $N$ and
$|\mathcal F_N|$ therefore gives
\[
N^{2/(3s)}|\mathcal F_N|^{2/3}
\lesssim_\varepsilon
\delta^{-2s/3-\varepsilon}
C_1^{2/3}K^{\sigma/(3s)}
|\mathcal P_N|^{2/3}.
\]
Raising both sides to the power $3s/2$
and using $\sigma/2=1-s$, we obtain
\[
N|\mathcal F_N|^s
\lesssim_\varepsilon
\delta^{-s^2-3s\varepsilon/2}
C_1^sK^{1-s}|\mathcal P_N|^s.
\]
Since $\mathcal I_N<2N|\mathcal F_N|$,
it follows that
\begin{equation*}
\mathcal I_N
\lesssim_\varepsilon
\delta^{-s^2-3s\varepsilon/2}
C_1^sK^{1-s}
|\mathcal F_N|^{1-s}|\mathcal P_N|^s.
\end{equation*}
Summing the preceding estimate over the cardinality classes and applying H\"older's inequality, we obtain
\[
\begin{aligned}
\mathcal I(\mathcal F,\mathcal P(\cdot))
&\lesssim_\varepsilon
\delta^{-s^2-3s\varepsilon/2}
C_1^sK^{1-s}
\sum_N|\mathcal F_N|^{1-s}|\mathcal P_N|^s\\
&\lesssim_\varepsilon
\delta^{-s^2-3s\varepsilon/2}
C_1^sK^{1-s}
|\mathcal F|^{1-s}
\left(\sum_N|\mathcal P_N|\right)^s.
\end{aligned}
\]
Although the sets $\mathcal P_N$ may overlap,
each is contained in $\mathcal P$, and there
are only $O(\log(1/\delta))$ classes. Hence
\[
\sum_N|\mathcal P_N|
\lesssim\log(1/\delta)|\mathcal P|,
\]
which gives
\[
\mathcal I(\mathcal F,\mathcal P(\cdot))
\lesssim_\varepsilon
\delta^{-s^2-3s\varepsilon/2}
\log(1/\delta)^s
C_1^sK^{1-s}
|\mathcal F|^{1-s}|\mathcal P|^s.
\]
Choosing $\varepsilon=\eta/(3s)$ and
absorbing the logarithmic factor into
$\delta^{-\eta/2}$ proves
\eqref{bi:density-incidence-bound}.
The implicit constant and scale threshold
are independent of $C_1,K$.
This completes the proof of the lemma.
\end{proof}

\subsection{A weighted incidence estimate}
\label{bi:weighted-section}

The preceding incidence estimate concerns finite families
of curves and squares, whereas the collision problem involves
measures. To pass from one setting to the other, we impose
two nonconcentration conditions: a first jet bound controlling
the distribution of curves and a tube bound controlling
the target mass near graph segments. By grouping cells
of comparable mass, we can apply
Lemma~\ref{bi:density-incidence} on each mass class. The cell masses then cancel in the resulting estimate.

\begin{lemma}\label{bi:weighted-incidences}
Let $0<s<1$, and let $\mathcal F\subset B_{C^2([-2,2])}(1/4)$
be a compact family transversal over $[-2,2]$, with constant $T$.
Let $\rho$ be a finite positive Borel measure on $\mathcal F$, and let
$\mu_{\mathrm{tar}}$ be a finite positive Borel measure on $\R^2$
supported in a fixed compact subset of $(-1,1)^2$.
Assume that their total masses are bounded.

Suppose that
\begin{equation}\label{bi:weighted-jets}
\rho\{h\in\mathcal F:A_x(h)\in J_1\times J_2\}
\lesssim (|J_1||J_2|)^s
\end{equation}
for all $x\in[-2,2]$ and bounded intervals $J_1,J_2\subset\R$.
Suppose also that, for every $0<\delta\leq r\leq1$,
every interval $J$ of length at most $r$,
and every $g\in C^1([-2,2])$ with $\|g\|_{C^1([-2,2])}\leq2$,
\begin{equation}\label{bi:weighted-tubes}
\mu_{\mathrm{tar}}\{(x,y)\in[-2,2]\times\R:x\in J,\ |y-g(x)|\leq\delta\}
\lesssim \delta^s r^s.
\end{equation}
Then, for every $\eta>0$ and every fixed $C_0\geq1$,
\begin{equation}\label{bi:weighted-incidence-bound}
\iint
\mathbf1_{\{|y-h(x)|\leq C_0\delta\}}
\,\dd\mu_{\mathrm{tar}}(x,y)\dd\rho(h)
\lesssim_{\eta,C_0}
\delta^{\min\{4s/3,\,2s-s^2\}-\eta}
\end{equation}
for all $0<\delta\leq1$.
The implicit constant may also depend on $s$,
the transversality constant of $\mathcal F$,
the total mass bounds and the constants in
\eqref{bi:weighted-jets}--\eqref{bi:weighted-tubes},
but is independent of $\delta$.
\end{lemma}

\begin{proof}
We discretize the measures $\rho$ and $\mu_{\mathrm{tar}}$
at scale $\delta$ and apply
Lemma~\ref{bi:density-incidence} to cells of
comparable mass. The first jet map provides
a natural discretization of the curve family,
since transversality controls the curves through
their values and slopes at a single point.

Fix a sufficiently small dyadic $\delta>0$.
For each dyadic $\delta$-square $Q$ in the
first jet plane and each dyadic $\delta$-square
$S$ in the target plane, set
\[
E_Q=\{h\in\mathcal F:A_0(h)\in Q\},
\qquad
m_Q=\rho(E_Q),
\qquad
n_S=\mu_{\mathrm{tar}}(S).
\]
By \eqref{eq:transversality} at $z=0$, any $h,\widetilde h\in E_Q$
satisfy
\[
\begin{aligned}
\|h-\widetilde h\|_{C^2([-2,2])}
&\leq T\bigl(
|h(0)-\widetilde h(0)|
+|h'(0)-\widetilde h'(0)|
\bigr)\leq2T\delta.
\end{aligned}
\]
In particular, their first jets differ by
$O(\delta)$ uniformly on $[-2,2]$.

The two partitions have $O(\delta^{-2})$ positive-mass cells, since $A_0(\mathcal F)$ and $\supp\mu_{\mathrm{tar}}$ lie in fixed bounded sets. Lemma~\ref{lem:dyadic-cell-masses} with $d=2$ and $\gamma=2$ bounds the discarded mass by $O(\delta^2)$ in each partition. As both measures have bounded total mass, these cells contribute $O(\delta^2)$ to \eqref{bi:weighted-incidence-bound}. The remaining cells form $O(\log(1/\delta))$ classes in each partition, with $w\leq m_Q<2w$ and $v\leq n_S<2v$.

Fix a pair of mass classes. To replace the
parameter cells by a discrete family of curves,
choose a representative $h_Q\in E_Q$ from
each selected cell. Representatives from
neighboring cells need not be $\delta$-separated.
A finite coloring of the dyadic first jet grid
resolves this difficulty: within each color,
distinct squares are separated by at least
$\delta$, and hence so are the corresponding
curves in the $C^2$ norm.

We work with one color at a time. Let $F$
be its family of representative curves and
$P$ the target squares in the chosen mass class.
Since the total masses are bounded,
\begin{equation}\label{bi:weighted-cell-counts}
|F|\lesssim w^{-1},
\qquad
|P|\lesssim v^{-1}.
\end{equation}

We first verify the rectangular Katz--Tao
condition for $F$. Fix $x\in[-2,2]$ and
intervals $J_1,J_2$ with $|J_1|,|J_2|\geq\delta$.
If $A_x(h_Q)\in J_1\times J_2$, then
transversality gives
$A_x(E_Q)\subset J_1^*\times J_2^*$,
where $J_i^*$ is a fixed enlargement of $J_i$ for $i=1,2$.
Since the parameter cells are disjoint
and each has $\rho$-mass at least $w$,
\eqref{bi:weighted-jets} yields
\[
w\,|\{h\in F:A_x(h)\in J_1\times J_2\}|
\lesssim(|J_1||J_2|)^s.
\]
Thus, $F$ is rectangular
$(\delta,2s,C_1)$-Katz--Tao with
\begin{equation}\label{bi:weighted-C1}
C_1\lesssim\frac{\delta^{2s}}{w}.
\end{equation}

We first select candidate squares.
For each $h_Q\in F$, let $\mathcal S(h_Q)$
consist of those squares $S\in P$ containing
a point $(x,y)$ for which
$|y-\widetilde h(x)|\leq C_0\delta$
for some $\widetilde h\in E_Q$.
At a witnessing point $(x_S,y_S)\in S$,
$|y_S-h_Q(x_S)|\le(C_0+2T)\delta$. Since $S$ has side length $\delta$
and $h_Q$ is uniformly Lipschitz, $|y-h_Q(x)|\le C\delta$ at every
$(x,y)\in S$, with $C$ independent of $S$ and $\delta$.

Fix $z\in\R^2$ and $\delta\leq r\leq1$.
The squares in $\mathcal S(h_Q)$ whose centers
lie in $B(z,r)$ are contained in a bounded
number of graph neighborhoods of width
$O(\delta)$ over intervals of length at most
$r$.
Each such neighborhood can be covered
by finitely many tubes of width $\delta$ about
vertical translates $h_Q+k\delta$, with $k$
ranging over a bounded set of integers.
For sufficiently small $\delta$, these translates
have $C^1$ norm at most $2$, so the tube hypothesis
applies at the original scales $\delta,r$,
including when $r=\delta$.
Applying \eqref{bi:weighted-tubes}
with fixed enlargements and using
$\mu_{\mathrm{tar}}(S)\geq v$, we obtain
\[
v\,|\mathcal S(h_Q)\cap B(z,r)|
\lesssim\delta^sr^s.
\]
Consequently, every candidate collection is
$(\delta,s,K)$-Katz--Tao with
\begin{equation}\label{bi:weighted-K}
K\lesssim\frac{\delta^{2s}}{v}.
\end{equation}
The single-cell cases of
\eqref{bi:weighted-jets} and
\eqref{bi:weighted-tubes} also imply
$w,v\lesssim\delta^{2s}$.
We may therefore choose $C_1,K\geq1$
while retaining \eqref{bi:weighted-C1}
and \eqref{bi:weighted-K}.

A selected square may lie near a representative graph without meeting it. To turn proximity into incidence, we use finitely many vertical translates of each graph. Indeed, for every
$S\in\mathcal S(h)$, there exists $(x,y)\in S$
such that $|y-h(x)|\leq C\delta$.
Since $S$ has side length $\delta$,
we can choose an integer $j$, with $|j|\leq M$
for some fixed $M$, such that
$(x,h(x)+j\delta)\in S$.

For each $|j|\leq M$, define $F_j=\{h+j\delta:h\in F\}$,
and assign each square in $\mathcal S(h)$
to one translate whose graph it meets.
Denote the resulting shading of $h+j\delta$
by $\mathcal P_j(h+j\delta)$ and set
$\mathcal P_j=\bigcup_{h\in F}\mathcal P_j(h+j\delta)$.
By construction,
\[
\sum_{|j|\leq M}
\mathcal I(F_j,\mathcal P_j(\cdot))
=\sum_{h\in F}|\mathcal S(h)|,
\qquad
\mathcal P_j\subset P.
\]

For each fixed $j$, vertical translation
preserves differences between curves,
so $F_j$ has the same transversality constant
and $\delta$-separation as $F$.
Moreover,
$A_x(h+j\delta)=A_x(h)+(j\delta,0)$,
which preserves the rectangular Katz--Tao
condition. Each translated shading is a
subset of the corresponding candidate collection
and therefore remains $(\delta,s,K)$-Katz--Tao.
Discarding any curves with empty translated
shadings does not affect these properties.

Finally, since $|j|\leq M$,
\[
\|h+j\delta\|_{C^2([-2,2])}
\leq\frac14+M\delta<1
\]
for sufficiently small $\delta$.
Thus, each translated family satisfies the
hypotheses of Lemma~\ref{bi:density-incidence}
with the same constants $C_1,K$, and $T$.

Set $\vartheta=\max\{2/3,s\}$.
Each contributing pair of cells $E_Q,S$
has product measure at most
$\rho(E_Q)\mu_{\mathrm{tar}}(S)<4wv$.
The construction above therefore bounds
the contribution of the selected mass classes
and color by
\[
4wv\sum_{|j|\leq M}
\mathcal I(F_j,\mathcal P_j(\cdot)).
\]
For every $|j|\le M$, Lemma~\ref{bi:density-incidence} with loss
$\eta/2$ gives
\[
wv\,\mathcal I(F_j,\mathcal P_j(\cdot))
\lesssim_\eta \delta^{-s\vartheta-\eta/2}
wv\,C_1^\vartheta K^{1-\vartheta}
|F_j|^{1-\vartheta}|\mathcal P_j|^\vartheta.
\]
By \eqref{bi:weighted-cell-counts}, $|F_j|\lesssim w^{-1}$ and
$|\mathcal P_j|\lesssim v^{-1}$. Thus the factor after the power of
$\delta$ is bounded by
\[
(w C_1)^\vartheta(v K)^{1-\vartheta}\lesssim\delta^{2s}
\]
by \eqref{bi:weighted-C1} and \eqref{bi:weighted-K}. Hence
$wv\,\mathcal I(F_j,\mathcal P_j(\cdot))\lesssim_\eta
\delta^{2s-s\vartheta-\eta/2}$, independently of the mass classes.

Set $b=\min\{4s/3,2s-s^2\}$.
Since $2s-s\vartheta=b$, each pair of
mass classes contributes at most
$C_\eta\delta^{b-\eta/2}$.
There are $O(\log^2(1/\delta))$ such pairs
and only finitely many colors and translates.
Including the discarded cells, we obtain
\[
\begin{aligned}
&\iint
\mathbf1_{\{|y-h(x)|\leq C_0\delta\}}
\,\dd\mu_{\mathrm{tar}}(x,y)\dd\rho(h)\\
&\qquad\lesssim_{\eta,C_0}
\delta^{b-\eta/2}\log^2(1/\delta)+\delta^2
\lesssim_{\eta,C_0}\delta^{b-\eta}.
\end{aligned}
\]
Here $b<1$, and we absorbed the logarithmic
factor using
$\log^2(1/\delta)\lesssim_\eta\delta^{-\eta/2}$.

For arbitrary sufficiently small $\delta>0$,
choose a dyadic $\delta'$ with
$\delta\leq\delta'<2\delta$ and use monotonicity.
The remaining scales follow from the bounded
total masses of $\rho$ and $\mu_{\mathrm{tar}}$.
This proves \eqref{bi:weighted-incidence-bound}.
{This finishes the proof of the lemma.}
\end{proof}
\subsection{Near-collision bounds for the image measure}
\label{bi:collisions}

To estimate near-collisions, we first use an exact level curve to relate the horizontal and vertical displacements.
At a horizontal separation scale $r>0$, write
$x'=x+rq$ and $y'=y+rz$, with $1\leq q<2$.
For an exact collision $f(x,y)=f(x',y')$, the nonvanishing
of $f_2$ first gives $|y'-y|\lesssim r$, by the mean value
theorem. Thus $z$ remains bounded, and Taylor's formula yields
\[
0=f(x+rq,y+rz)-f(x,y)
 =r\bigl(qf_1(x,y)+zf_2(x,y)\bigr)+O(r^2).
\]
It follows that $z=-qR_f(x,y)+O(r)$, where $R_f=f_1/f_2$.
The graph of $y\mapsto-qR_f(x,y)$ therefore describes the
limiting normalized displacement along a level curve.

To count near-collisions at arbitrarily small tolerances, we
replace this approximation by an exact family of collision
curves. The nonvanishing of $\mathcal K=\partial_{12}\log|R_f|$
ensures that their values and slopes vary independently with
$(x,q)$. After normalization, near-collisions become incidences
with tubes of width $O(\delta/r)$ around these curves.
The first jet and tube estimates below allow us to apply
Lemma~\ref{bi:weighted-incidences} at each separation scale.

\begin{proof}[Proof of Proposition~\ref{bi:collision-bound}]
Use the family $H_r$ from \eqref{bi:normalized-collision-curve}
on the local chart of Subsection~\ref{bi:localization-section}. We now
rescale it for the weighted incidence estimate.

Let $V$ be a regular rectangle about the product-support point on
which $(x,f(x,y))\in U_{\mathrm{inv}}$ and $Y(x,f(x,y))=y$. Let $y_0$
be the second coordinate of that point. Choose $\ell>0$ and shrink
$I_1,I_2$ about the same point so that $I_2\subset y_0+\ell(-1/2,1/2)$
and $I_1\times(y_0+\ell[-2,2])\Subset V$. Choose $0<r_0\le1$ so that
the shifted inverse arguments $(x+rq,f(x,y_0+\ell\theta))$ stay in
$U_{\mathrm{inv}}$ and their inverse points remain in the larger
regular rectangle. Thus $H_r(x,q;y_0+\ell\theta)$ and its derivatives
extend smoothly on the range
\[
x\in I_1,\quad q\in[1,2],\quad
\theta\in[-2,2],\quad 0\leq r\leq r_0,
\]
and the nondegeneracy just established persists uniformly.
With this choice of $r_0$, shrink $I_1$ once more so that
$\operatorname{diam}(I_1)<r_0$.

For later use, fix constants on the larger regular rectangle,
\[
M_1=\sup|f_1|<\infty,\qquad
m_2=\inf|f_2|>0,
\qquad {C_{\mathrm{vert}}}=\frac{1+2M_1}{m_2}.
\]
Choose $c>0$ once and for all so small that $cC_{\mathrm{vert}}<1/2$
and, with $h_{r,x,q}(\theta)=cH_r(x,q;y_0+\ell\theta)$,
we have
$\|h_{r,x,q}\|_{C^2([-2,2])}<1/4$
uniformly on the above parameter domain. All subsequent
quantitative constants refer to this final normalization.

Write $p=(x,q)$ and
$F_{r,\theta}(p)=A_\theta(h_{r,p})$.
The chain rule and \eqref{bi:collision-jet} give
\[
\det D_pF_{0,\theta}(p)
=-c^2\ell qR_f(x,y)^2\mathcal K(x,y),
\qquad y=y_0+\ell\theta.
\]
By smooth convergence, after reducing $r_0$ and, if needed,
shrinking $I_1$ again to retain
$\operatorname{diam}(I_1)<r_0$, there are constants
$d,M>0$ such that
\[
|\det D_pF_{r,\theta}(p)|\geq d,
\qquad \|D_pF_{r,\theta}(p)\|_{\mathrm{op}}\leq M
\]
for all admissible $(r,\theta,p)$.
For a $2\times2$ matrix, the least singular value is the
absolute determinant divided by the largest singular value.
Consequently, with $m=d/M>0$,
\begin{equation}\label{bi:jet-singular-lower}
|D_pF_{r,\theta}(p)v|\geq m|v|,
\qquad v\in\R^2.
\end{equation}

Pointwise invertibility does not by itself imply an inverse
Lipschitz bound on the full parameter rectangle.
Uniform continuity of $D_pF_{r,\theta}$ allows us to cover
$I_1\times[1,2]$ by finitely many closed convex rectangles
$Q_1,\ldots,Q_N$ such that, for a reference point $p_a\in Q_a$,
\[
\|D_pF_{r,\theta}(p)-D_pF_{r,\theta}(p_a)\|_{\mathrm{op}}
\leq\frac m2,
\qquad p\in Q_a,
\]
uniformly in $r$ and $\theta$.
The segment argument in the proof of
Lemma~\ref{lem:local-transverse}, with reference derivative
$D_pF_{r,\theta}(p_a)$ and \eqref{bi:jet-singular-lower}, gives the
lower bound below. The uniform derivative bound gives the upper bound.
For every $p,p'\in Q_a$, we have
\begin{equation}\label{bi:jet-bilipschitz}
\frac m2|p-p'|
\leq |A_\theta(h_{r,p})-A_\theta(h_{r,p'})|
\leq M|p-p'|.
\end{equation}
Smooth parameter dependence also gives
$\|h_{r,p}-h_{r,p'}\|_{C^2([-2,2])}
\leq L_h|p-p'|$ for a fixed $L_h$.
Thus each family
$\mathcal F_{r,a}=\{h_{r,p}:p\in Q_a\}$
is compact and uniformly transversal over $[-2,2]$,
with a transversality constant independent of $r$ and $a$.

We also need the stronger one-dimensional control in $q$.
Put $y=y_0+\ell\theta$ and
$y_r=Y(x+rq,f(x,y))$.
Differentiating the defining identities gives
\begin{equation}\label{bi:q-derivative-value}
\partial_qh_{r,x,q}(\theta)=-cR_f(x+rq,y_r)
\end{equation}
and
\begin{equation}\label{bi:q-derivative-slope}
\partial_qh_{r,x,q}'(\theta)
=-c\ell\,\partial_2R_f(x+rq,y_r)
\frac{f_2(x,y)}{f_2(x+rq,y_r)}.
\end{equation}
Indeed, the second identity uses
$\partial_y y_r=f_2(x,y)/f_2(x+rq,y_r)$.
By \eqref{bi:local-geometry}, both $q$-derivatives have fixed
signs and absolute values bounded below by positive constants
$a_1,a_2$, uniformly in $r,\theta,x,q$.
These fixed geometric choices are independent of $\nu,\mu$
and $\delta$.

We now fix a dyadic separation
$r\leq x'-x<2r$, with $\delta\leq r\leq r_0$.
The near-collision condition forces the vertical displacement
to be $O(r)$. Indeed, the mean value theorem and the fixed
sign of $f_2$ give
\[
\begin{aligned}
m_2|y'-y|
&\leq |f(x',y')-f(x',y)|\\
&\leq |f(x',y')-f(x,y)|+|f(x,y)-f(x',y)|\\
&\leq \delta+M_1|x'-x|\\
&\leq (1+2M_1)r.
\end{aligned}
\]
Thus $|y'-y|\leq C_{\mathrm{vert}}r$.
With the target coordinates
\[
\theta=\frac{y-y_0}{\ell},\qquad
\widetilde z=\frac{c(y'-y)}{r},
\]
our choices ensure $|\theta|<1/2$ and $|\widetilde z|<1/2$.
Thus $\widetilde z=cz$, where $z=(y'-y)/r$ is the vertical displacement
used in the preliminary expansion.

We construct the curve measures separately on the fixed
parameter charts. Let
$$G_r=\{(x,x')\in I_1^2:r\leq x'-x<2r\}.$$
Define $\Pi_r:G_r\to I_1\times[1,2]$ by $\Pi_r(x,x')=(x,(x'-x)/r)$,
and $\Gamma_r:I_1\times[1,2]\to C^2([-2,2])$ by $\Gamma_r(p)=h_{r,p}$.
For each $1\leq a\leq N$, define
\[
\rho_{r,a}
=r^{-s}(\Gamma_r\circ\Pi_r)_*
\bigl((\nu\otimes\nu)|_{G_r\cap\Pi_r^{-1}(Q_a)}\bigr).
\]
This measure is supported on the compact transversal family
$\mathcal F_{r,a}$. Its mass satisfies
\[
\|\rho_{r,a}\|
\leq r^{-s}\int\nu([x+r,x+2r))\,\dd\nu(x)
\lesssim1.
\]
To verify the rectangular first jet bound, fix
$\theta\in[-2,2]$ and intervals $J_1,J_2$ of lengths
$t_1,t_2>0$. Let $E_a=\{p\in Q_a:F_{r,\theta}(p)\in J_1\times J_2\}$.
By \eqref{bi:jet-bilipschitz}, any two points of $E_a$
have first coordinates at distance
$O(\max\{t_1,t_2\})$. For fixed $x$, the bounds
$|\partial_qh|\geq a_1$ and $|\partial_qh'|\geq a_2$ from \eqref{bi:q-derivative-value} and \eqref{bi:q-derivative-slope},
together with their fixed signs, show that any two
admissible $q$ values satisfy
\[
|q-q'|\leq\min\{t_1/a_1,t_2/a_2\}
\lesssim\min\{t_1,t_2\}.
\]
Thus the $x$-projection of $E_a$ lies in an interval of
length $O(\max\{t_1,t_2\})$, and for each fixed $x$
the corresponding $x'=x+rq$ lie in an interval of length
$O(r\min\{t_1,t_2\})$.
Fubini and the Frostman bounds in $x$ and $x'$ give, with the
factor $r^{-s}$ from the definition of $\rho_{r,a}$,
\begin{equation}\label{bi:collision-rectangular}
\begin{aligned}
\rho_{r,a}\{h:A_\theta(h)\in J_1\times J_2\}
&\lesssim r^{-s}\max\{t_1,t_2\}^s
\bigl(r\min\{t_1,t_2\}\bigr)^s=(t_1t_2)^s.
\end{aligned}
\end{equation}
The bound is uniform in $r,\theta,a$.

The target measure records the normalized second-coordinate pairs. Set $E_r=\{(y,y')\in I_2^2:|y'-y|\leq C_{\mathrm{vert}}r\}$
and define
\begin{align*}
\Psi_r:E_r&\longrightarrow[-1/2,1/2]^2,\\
(y,y')&\longmapsto(\theta,\widetilde z)
=\left(\frac{y-y_0}{\ell},\frac{c(y'-y)}r\right).
\end{align*}
Set
\[
\mu_{\mathrm{tar},r}
=r^{-s}(\Psi_r)_*\bigl((\mu\otimes\mu)|_{E_r}\bigr).
\]
The Frostman bound yields
\[
\|\mu_{\mathrm{tar},r}\|
\leq r^{-s}\int\mu([y-C_{\mathrm{vert}}r,y+C_{\mathrm{vert}}r])\,\dd\mu(y)
\lesssim1.
\]
Moreover, $\supp\mu_{\mathrm{tar},r}\subset[-1/2,1/2]^2$,
independently of $r$.

For the tube condition, fix $0<\e\leq t\leq1$,
an interval $J$ of length at most $t$, and
$g\in C^1([-2,2])$ with $\|g\|_{C^1}\leq2$.
Fix also a constant $C\geq1$;
the implicit constant in the following estimate
may depend on $C$.
The condition $\theta\in J$ restricts $y$ to an interval
of length at most $\ell t$. For each such $y$,
\[
\left|\frac{c(y'-y)}r
-g\left(\frac{y-y_0}{\ell}\right)\right|\leq C\e
\]
restricts $y'$ to an interval of length $2Cr\e/c$.
Fubini and the Frostman bounds in $y$ and $y'$ give, with the
factor $r^{-s}$ from the definition of $\mu_{\mathrm{tar},r}$,
\begin{equation}\label{bi:collision-tube-window}
\begin{aligned}
&\mu_{\mathrm{tar},r}\{(\theta,\widetilde z){\in[-2,2]\times\R}:\theta\in J,
\ |\widetilde z-g(\theta)|\leq C\e\}\lesssim r^{-s}(r\e)^s t^s
=\e^st^s.
\end{aligned}
\end{equation}
Thus $\mu_{\mathrm{tar},r}$ satisfies \eqref{bi:weighted-tubes}
uniformly in $r$.

To turn a near-collision into a tube condition, recall that $y_r=Y(x',f(x,y))$, so that
$f(x',y_r)=f(x,y)$. The two points $(x',y')$ and $(x',y_r)$
lie in the larger regular rectangle. At fixed $x'$, the mean value
theorem therefore gives
\[
m_2|y'-y_r|\le |f(x',y')-f(x',y_r)|\le\delta.
\]
With $q=(x'-x)/r$, the definitions give
$\widetilde z-h_{r,x,q}(\theta)=(c/r)(y'-y_r)$. Set $C_Y=m_2^{-1}$.

Consequently,
\begin{equation}\label{bi:collision-tube-reduction}
|\widetilde z-h_{r,x,q}(\theta)|\leq C_*\frac{\delta}{r},
\qquad C_*:=\max\{1,cC_Y\}.
\end{equation}

Let $\mathcal C_r(\delta)$ be the
$(\nu\otimes\mu)^{\otimes2}$-mass of the $\delta$-near-collisions $|f(x,y)-f(x',y')|\leq\delta$ with
$r\leq x'-x<2r$. Since $|y'-y|\leq C_{\mathrm{vert}}r$ for every
such collision and the rectangles $Q_a$ cover the parameter
domain, \eqref{bi:collision-tube-reduction} and the pushforward definitions give
\[
\begin{aligned}
\mathcal C_r(\delta)
&\leq r^{2s}\sum_{a=1}^N
\iint\mathbf1_{\{|\widetilde z-h(\theta)|\leq C_*\delta/r\}}
\,\dd\mu_{\mathrm{tar},r}(\theta,\widetilde z)\dd\rho_{r,a}(h).
\end{aligned}
\]
The factor $r^{2s}$ restores the two factors $r^{-s}$
in the normalized measures; no Lebesgue Jacobian is involved.

Set $b=\min\{4s/3,2s-s^2\}$.
By \eqref{bi:collision-rectangular} and
\eqref{bi:collision-tube-window}, the hypotheses of
Lemma~\ref{bi:weighted-incidences} hold for each
$(\mathcal F_{r,a},\rho_{r,a},\mu_{\mathrm{tar},r})$, with
constants independent of $r$ and $a$.
Apply the lemma at tube scale $\varepsilon=\delta/r\in(0,1]$,
with width multiplier $C_*$. This gives
\[
\begin{aligned}
\mathcal C_r(\delta)
&\lesssim_\eta r^{2s}
\left(\frac\delta r\right)^{b-\eta}=\delta^{b-\eta}r^{2s-b+\eta}.
\end{aligned}
\]
As $2s-b=\max\{2s/3,s^2\}>0$,
summing over dyadic $\delta\leq r\leq r_0$ yields
\[
\sum_{\substack{r\in2^{\mathbb Z}\\
\delta\leq r\leq r_0}}\mathcal C_r(\delta)
\lesssim_\eta\delta^{b-\eta}.
\]
Every positive separation $x'-x\geq2\delta$
lies in one of the dyadic shells $[r,2r)$ with $r\geq\delta$.
The remaining near-diagonal pairs satisfy
$|x-x'|<2\delta$. Their $\nu\otimes\nu$-mass
is $O(\delta^s)$. If they also form a $\delta$-collision,
the mean value theorem as above gives
\[
m_2|y'-y|
\leq\delta+M_1|x'-x|
\leq(1+2M_1)\delta,
\]
so, for each such $(x,x')$, the admissible $(y,y')$ have
$\mu\otimes\mu$-mass $O(\delta^s)$, uniformly in $(x,x')$.
Integrating in $(x,x')$ gives a contribution of $O(\delta^{2s})$,
which is bounded by $O(\delta^{b-\eta})$ because $2s>b$.
Negative horizontal separations are treated by exchanging
the two input pairs; the diagonal has zero mass because
$\nu$ is nonatomic. Since the final diameter of $I_1$
is less than $r_0$, these cases exhaust the collision
region. The estimate for $\delta$ bounded away from zero follows
by enlarging the implicit constant. This proves
\eqref{bi:collision-estimate}.

Let $\omega=f_*(\nu\otimes\mu)$ and set $b=\min\{4s/3,2s-s^2\}$.
For $0<u<b$, choose $0<\eta<b-u$. By the layer-cake formula
and \eqref{bi:collision-estimate},
\[
I_u(\omega)
=u\int_0^\infty r^{-u-1}
 (\omega\otimes\omega)\{|v-v'|<r\}\,\dd r
\lesssim_{u,\eta} 1+\int_0^1 r^{b-u-\eta-1}\,\dd r<\infty.
\]
This completes the proof of the proposition.
\end{proof}

\subsection{\texorpdfstring{Proof of Theorem~\ref{thm:bivariate-dimension-expansion}}{Proof of Theorem~\ref{thm:bivariate-dimension-expansion}}}
\label{bi:proof-completion}

We first deduce a dimension bound from the Fourier shell estimate when the input dimensions differ. We then combine this estimate with the incidence bound of Proposition~\ref{bi:collision-bound} to prove
Theorem~\ref{thm:bivariate-dimension-expansion}.

\begin{theorem}\label{bi:unequal-dimensions}
    Let $\Omega\subset\R^2$ be a connected open neighborhood of $[0,1]^2$.
Let $f:\Omega\to\R$ be real analytic and nonspecial in the sense of
\eqref{eq:special form}, with neither $\partial_x f$ nor $\partial_y f$
identically zero. Let $A,B\subset[0,1]$ be Borel sets with positive
Hausdorff dimension. If $\dimH A+\dimH B>1/2$, then
\[
\dimH f(A\times B)\ge
\min\left\{1,\dimH A+\dimH B-\frac12\right\}.
\]
If $\dimH A+\dimH B>3/2$, then $\Leb(f(A\times B))>0$.
\end{theorem}

\begin{proof}
Choose $0<s_1<\dimH A$ and $0<s_2<\dimH B$ with $s_1+s_2>1/2$.
By Frostman's lemma, take $s_1$- and $s_2$-Frostman probabilities
supported on compact subsets of $A$ and $B$, respectively.
Apply Lemma~\ref{bi:localization} with $s=\min\{s_1,s_2\}$
and shrink the intervals about the resulting support point so that
Proposition~\ref{bi:radon-shell} applies. The normalized restrictions
$\nu$ and $\mu$ are still $s_1$- and $s_2$-Frostman, respectively.
Proposition~\ref{bi:radon-shell} gives image dimension at least
$\min\{1,s_1+s_2-1/2\}$. Letting $s_1\uparrow\dimH A$ and
$s_2\uparrow\dimH B$ proves the dimension bound.
If $\dimH A+\dimH B>3/2$, choose $s_1+s_2>3/2$ instead.
The same proposition gives an $L^2$ density for $f_*(\nu\otimes\mu)$,
a probability carried by the compact set
$f(\supp\nu\times\supp\mu)\subset f(A\times B)$. This proves the
positive-measure conclusion and completes the proof of the theorem.
\end{proof}

\begin{proof}[Proof of Theorem~\ref{thm:bivariate-dimension-expansion}]
Let $\Omega\subset\R^2$ be a connected open neighborhood of
$[0,1]^2$, and let $f\in C^\omega(\Omega)$, with neither
$f_1$ nor $f_2$ identically zero.
If $f$ has analytic special form \eqref{eq:special form},
then alternative (i) holds. We may therefore assume that
$f$ is nonspecial.

Let $A,B\subset[0,1]$ be Borel sets with
$\dimH A,\dimH B\ge\alpha$, where $0<\alpha\le1$.
Fix $0<s<\alpha$. By Frostman's lemma, choose
$s$-Frostman probability measures $\nu_0,\mu_0$
supported on compact subsets of $A$ and $B$, respectively.

By Lemma~\ref{bi:localization}, choose compact intervals $J_1,J_2$
whose product contains a point of $\supp\nu_0\times\supp\mu_0$ in
its interior. Choose them sufficiently small that
\eqref{bi:local-geometry} holds on a neighborhood of $J_1\times J_2$
and both Propositions~\ref{bi:collision-bound} and~\ref{bi:radon-shell}
apply to measures supported on these intervals. The normalized
restrictions $\nu=\nu_0|_{J_1}/\nu_0(J_1)$ and
$\mu=\mu_0|_{J_2}/\mu_0(J_2)$ are $s$-Frostman probabilities.
Both estimates apply to $\omega=f_*(\nu\otimes\mu)$, a probability
carried by the compact set
$f(\supp\nu\times\supp\mu)\subset f(A\times B)$.

Set $b(s)=\min\{4s/3,2s-s^2\}$. Proposition~\ref{bi:collision-bound}
gives finite $u$-energy for every $0<u<b(s)$, so
$\dimH\supp\omega\ge b(s)$ by \cite[Theorem~2.8]{Mat15}.
Together with Proposition~\ref{bi:radon-shell}, this gives
\[
\dimH f(A\times B)\ge
\max\{b(s),\min\{2s-1/2,1\}\};
\]
if $s\le1/4$, the second term is nonpositive and can be ignored.
The two terms defining $b(s)$ switch at $s=2/3$. Up to that point
$4s/3\ge 2s-1/2$; after it, the difference between $2s-s^2$ and
$2s-1/2$ is $1/2-s^2$. The Fourier term thus takes over at
$s=1/\sqrt 2$ and reaches one at $s=3/4$.
Letting $s\uparrow\alpha$ gives \eqref{eq:bivariate-bound}.

If $\dimH A+\dimH B>3/2$, the positive-measure conclusion
follows from Theorem~\ref{bi:unequal-dimensions}.
This finishes the proof of the theorem.
\end{proof}

\begin{remark}
Under the hypotheses of Theorem~\ref{thm:intro-positive-bivariate},
a local representation as in \eqref{eq:special form} would give
$\partial_{xy}\log|f_1/f_2|=0$ there. By the identity theorem,
the analytic numerator $\mathcal N$ in \eqref{bi:numerator}
would vanish throughout the connected neighborhood. Since
$f_1f_2$ is nowhere zero on the compact square, it remains
nonzero on an open rectangle containing the square. There,
$\log|f_1/f_2|=a(x)+b(y)$ for univariate analytic $a,b$.

Choose analytic coordinates $u=H_1(x)$ and $v=H_2(y)$ by taking
\[
 H_1'(x)=\operatorname{sgn}(f_1/f_2)e^{a(x)},
 \qquad H_2'(y)=e^{-b(y)}.
\]
The two partial derivatives of $f$ in these coordinates agree.
Each fiber of $u+v$ in the coordinate rectangle is a connected
segment, so $f=G(u+v)$. Fixing $v$ locally shows that $G$ is
analytic on the sum interval. This gives the excluded
representation on the square.

Theorem~\ref{thm:intro-positive-bivariate} follows from
the proof of Theorem~\ref{thm:bivariate-dimension-expansion}
by taking $\alpha=\min\{\dimH A,\dimH B\}>3/4$.
\end{remark}

\section{Spectral applications}
\label{sec:spectral-averaging}

We apply the Fourier shell estimates and analytic localization
to spectral problems for Jacobi operators. In the finite-dimensional
setting, we establish absolute continuity of eigenvalue laws
and of linear eigenvalue statistics with coefficients that are not all
equal. To treat linear statistics, we reduce the nondegeneracy test to three consecutive sites.
For finite averaging at the first $k\ge3$ sites, we verify the
Monge--Amp\`ere condition by a balanced partition of the collision
variables. For half-line operators,
estimates for the boundary ratio yield absolute continuity
of averaged spectral measures and applications to boundary
perturbations.

\subsection{Absolute continuity of level spacings and linear combinations of eigenvalues}
\label{spectral:proof-stats}

We prove Theorem~\ref{spectral:stats-all}, which also yields
Theorem~\ref{spectral:stats-three}. Indeed, taking $n=3$ and
$c=e_j-e_i$, where $e_i$ denotes the $i$th standard basis vector,
gives the level spacing $F_c=\lambda_j-\lambda_i$
for every $1\le i<j\le3$.

\begin{proof}[Proof of Theorem~\ref{spectral:stats-all}]

Let $H(\mathbf x)$ be the $n\times n$ Jacobi matrix \eqref{eq:finite-jacobi-matrix} with diagonal
$\mathbf x\in\R^n$, where $n\ge3$, and fixed nonzero off-diagonal
entries. Write $\lambda_1(\mathbf x)<\cdots<\lambda_n(\mathbf x)$
for its ordered eigenvalues. The nonzero off-diagonal
entries make $H(\mathbf x)$ irreducible, so its eigenvalues are
simple for every $\mathbf x\in\R^n$. Since $H$ depends real
analytically on $\mathbf x$, the ordered eigenvalue functions are
real analytic on $\R^n$. Fix $c\in\R^n$ whose entries are not all
equal and set $F_c=\sum_{j=1}^n c_j\lambda_j$. Let the diagonal variables $X_m$ be independent, with compactly
supported nonatomic laws $\nu_m$ for $1\le m\le n$. Assume that for the consecutive triple
$\mathcal J=\{r,r+1,r+2\}$, where $1\le r\le n-2$, the laws $\nu_m$ are $s_m$-Frostman
with $0<s_m\le1$ for $m\in\mathcal J$ and
$\sum_{m\in\mathcal J}s_m>2$.

To apply Theorem~\ref{thm:balanced-collision-products}(i),
we need to verify the Monge--Amp\`ere condition \eqref{analytic:balanced-product-condition} for
the three diagonal variables indexed by $\mathcal J$.
By \eqref{analytic:trivariate-balanced-determinant},
it suffices to find an ordering $(i,j,\ell)$ of
$\mathcal J$ such that
\[
\Gaux_{i;j\ell}[F_c]
=(\partial_jF_c)(\partial_{i\ell}F_c)
-(\partial_\ell F_c)(\partial_{ij}F_c)
\not\equiv0.
\]
We first rule out analytic special form for the three-site statistic.
We then transfer its nondegeneracy to the consecutive block in the full
matrix and apply the balanced product theorem on almost every slice.

\begin{claim}
If $n=3$ and the coordinates of $c$ are not all equal,
then $F_c$ is nonspecial, and none of its first
partial derivatives vanishes identically.
\end{claim}

Since $H(\mathbf x)$ depends real analytically
on $\mathbf x$ and has simple eigenvalues, each ordered eigenvalue
$\lambda_j(\mathbf x)$ is real analytic on $\R^3$.
Let $\tau\in\R$ denote the spectral parameter, and define
\[
P(\tau,\mathbf x)
=\det(H(\mathbf x)-\tau\mathrm{Id}_3),
\]
where $\mathrm{Id}_3$ is the $3\times3$ identity matrix.
On the one hand, differentiating
$P(\lambda_j(\mathbf x),\mathbf x)=0$
with respect to $x_i$ gives
\[
\partial_i\lambda_j(\mathbf x)
=-\frac{P_{x_i}(\lambda_j(\mathbf x),\mathbf x)}
        {P_\tau(\lambda_j(\mathbf x),\mathbf x)}, \quad \text{ for }1\leq i\leq 3.
\]
Here $P_{x_i}$ and $P_\tau$ denote the partial derivatives
with respect to $x_i$ and $\tau$, respectively.
Since $\lambda_j(\mathbf x)$ is simple,
$P_\tau(\lambda_j(\mathbf x),\mathbf x)\ne0$.

We will also use the eigenvector formula for the same derivative.
Varying only the $i$th diagonal entry gives the
one-parameter perturbation
\[
T_i(t)=H(\mathbf x+te_i)
      =H(\mathbf x)+te_ie_i^*, \quad t\in \R,
\]
where $e_i$ is the $i$th standard basis vector.
Since $T_i'(0)=e_ie_i^*$, Kato's perturbation theorem
\cite[Chapter~II, Theorem~5.4]{Kato95} gives
\[
\partial_i\lambda_j(\mathbf x)
=\langle T_i'(0)u_j,u_j\rangle
=|\langle u_j,e_i\rangle|^2,
\]
where $u_j$ is a normalized eigenvector associated
with the simple eigenvalue $\lambda_j(\mathbf x)$.
Consequently, the two formulas give, for $1\leq i \leq 3$,
\begin{equation}\label{spectral:eigenvalue-statistic-derivatives}
\partial_iF_c(\mathbf x)
=-\sum_{j=1}^3c_j
\frac{P_{x_i}(\lambda_j(\mathbf x),\mathbf x)}
     {P_\tau(\lambda_j(\mathbf x),\mathbf x)}
=\sum_{j=1}^3c_j
|\langle u_j(\mathbf x),e_i\rangle|^2.
\end{equation}
Since $P(\tau,\mathbf x)$ is a polynomial, each eigenvalue
$\lambda_j(\mathbf x)$ is algebraic in any one diagonal
variable when the others are fixed. The same holds for
$\partial_iF_c$ and the ratios
$\partial_iF_c/\partial_jF_c$ wherever they are defined.
We use this algebraic dependence to rule out special
representations of $F_c$.

To identify the limiting gradient of $F_c$, we let the diagonal entries become large while keeping the off-diagonal entries fixed.
Choose $\mathbf v=(v_1,v_2,v_3)\in\R^3$ with distinct
entries, and let $M>0$ be a scaling parameter.
The matrices $H(M\mathbf v)$ and
$M^{-1}H(M\mathbf v)$ have the same eigenspaces in the same order. The fixed
off-diagonal entries of the latter matrix tend to zero, so
\[
M^{-1}H(M\mathbf v)
\longrightarrow\operatorname{diag}(\mathbf v)
\qquad\text{as }M\to\infty.
\]
Let $\pi$ be the permutation of $\{1,2,3\}$ such that $v_{\pi(1)}<v_{\pi(2)}<v_{\pi(3)}$.
The limiting matrix has simple eigenvalues
$v_{\pi(j)}$ with corresponding eigenvectors
$e_{\pi(j)}$. By continuity of the eigenspaces
associated with simple eigenvalues, the normalized
eigenvectors $u_j(M\mathbf v)$ converge, up to sign,
to $e_{\pi(j)}$. Consequently,
\[
|\langle u_j(M\mathbf v),e_i\rangle|^2
\longrightarrow
|\langle e_{\pi(j)},e_i\rangle|^2
=\delta_{i,\pi(j)},
\]
for every $1\le i,j\le3$, where $\delta_{i,\pi(j)}$
is the Kronecker delta.
Thus each ordered eigenvector becomes concentrated
in the coordinate corresponding to the ordering
of the diagonal entries.

The formula \eqref{spectral:eigenvalue-statistic-derivatives} now gives
\[
\partial_iF_c(M\mathbf v)
=\sum_{j=1}^3c_j
|\langle u_j(M\mathbf v),e_i\rangle|^2
\longrightarrow c_{\pi^{-1}(i)}.
\]
Hence
\[
\nabla F_c(M\mathbf v)
\longrightarrow
(c_{\pi^{-1}(1)},c_{\pi^{-1}(2)},c_{\pi^{-1}(3)}).
\]
By varying the ordering of $\mathbf v$, every permutation
of $c$ occurs as a limit of $\nabla F_c(M\mathbf v)$.
Since the coordinates of $c$ are not all equal,
there are at least three distinct permutations,
all with the same positive Euclidean norm.
A one-dimensional subspace contains at most two vectors
of any given nonzero norm. Hence the gradients of $F_c$
cannot all lie in a fixed one-dimensional subspace. Also, for each
coordinate $i$, one can place a nonzero coefficient of $c$ in that
coordinate of a limiting gradient. Thus no $\partial_iF_c$ vanishes
identically.
We use this fact to rule out analytic special
representations of $F_c$.

Suppose, toward a contradiction, that $F_c$ has the
analytic special form \eqref{eq:special form} on a
nonempty product box $Q=I_1\times I_2\times I_3$:
\begin{equation}\label{spectral:statistic-special-form}
F_c(\mathbf x)=g(S(\mathbf x)),
\qquad
S(\mathbf x)=\sum_{i=1}^3h_i(x_i),
\end{equation}
where each $h_i:I_i\to\R$ is real analytic and
$g$ is real analytic on an open interval containing $S(Q)$.
Since no first partial derivative of $F_c$ vanishes
identically, analyticity allows us to shrink $Q$
so that
\begin{equation}\label{spectral:statistic-special-derivatives}
\partial_iF_c(\mathbf x)=g'(S)h_i'(x_i)\ne0
\end{equation}
for every $i\in\{1,2,3\}$ throughout $Q$.
The special representation is strongly constrained by
the behavior of eigenvalues under diagonal shifts.
Indeed, adding the same constant $t$ to all diagonal
entries gives
\[
H(\mathbf x+t(1,1,1))
=H(\mathbf x)+t\mathrm{Id}_3,
\]
so every eigenvalue increases by $t$. Consequently,
\[
F_c(\mathbf x+t(1,1,1))
=F_c(\mathbf x)+t\sum_{j=1}^3c_j.
\]
Differentiating with respect to $t$ at $t=0$ and using
\eqref{spectral:statistic-special-form}, we obtain $g'(S)\sum_{i=1}^3h_i'(x_i)
=\sum_{j=1}^3c_j$.
Differentiating with respect to $x_i$ and applying
the product and chain rules gives
\[
g''(S)h_i'(x_i)\sum_{j=1}^3h_j'(x_j)
+g'(S)h_i''(x_i)=0.
\]
Substituting the preceding identity and dividing
by $g'(S)h_i'(x_i)\ne0$, we obtain
\[
\frac{h_i''(x_i)}{h_i'(x_i)}
=-\frac{g''(S)}{[g'(S)]^2}
\sum_{j=1}^3c_j,
\qquad 1\le i\le3.
\]
Thus the three ratios $h_i''/h_i'$ agree on $Q$,
although each depends on a different variable.
They must therefore equal a common constant $\gamma$.
Solving $h_i''=\gamma h_i'$ gives $h_i'(x_i)=a_i e^{\gamma x_i}$, for $a_i\ne0$,
for every $i\in\{1,2,3\}$.
The formula \eqref{spectral:statistic-special-derivatives} then gives
\[
\frac{\partial_iF_c}{\partial_jF_c}
=\frac{h_i'(x_i)}{h_j'(x_j)}
=\frac{a_i}{a_j}e^{\gamma(x_i-x_j)},
\qquad i\ne j.
\]
With the other variables fixed, by the discussion after \eqref{spectral:eigenvalue-statistic-derivatives}, the left-hand side
is algebraic in $x_i$, whereas $e^{\gamma x_i}$
is transcendental over $\R(x_i)$ whenever
$\gamma\ne0$. Indeed, a polynomial relation
for $e^{\gamma x_i}$ would extend analytically
to all real $x_i$, contradicting exponential
growth in the appropriate direction.
Hence $\gamma=0$, and every $h_i'=a_i$ is constant.

Consequently,
\[
\nabla F_c(\mathbf x)
=g'(S(\mathbf x))(a_1,a_2,a_3)
\]
on $Q$. Hence $a_j\partial_iF_c-a_i\partial_jF_c=0$ on $Q$ for all
$1\le i,j\le3$. These functions are real analytic on the connected
space $\R^3$, so the identities hold everywhere. Thus every gradient
of $F_c$ lies in the fixed line spanned by $(a_1,a_2,a_3)$.
This contradicts the permutation limits of
$\nabla F_c(M\mathbf v)$ obtained above.
Therefore $F_c$ has no analytic special
representation on any regular product box.
By Theorem~\ref{thm:analytic-special-form},
there exists an ordering $(i,j,\ell)$ of
$\{1,2,3\}$ such that
\[
\Gaux_{i;j\ell}[F_c]\not\equiv0.
\]
The three-site case establishes the required nondegeneracy
only for $n=3$. When $n>3$, we must verify the same condition
for the prescribed consecutive triple $\mathcal J$, even though
$F_c$ depends on all $n$ diagonal variables. The following
claim achieves this by reducing the general matrix to a
three-site block.

\begin{claim}[]
Let $n>3$, and let $\mathcal J=\{r,r+1,r+2\}$
be a consecutive triple. If the coordinates of $c$
are not all equal, then there exists an ordering
$(i,j,\ell)$ of $\mathcal J$ such that
\[
\Gaux_{i;j\ell}[F_c]\not\equiv0.
\]
In particular, $F_c$ is nonspecial.
\end{claim}
To prove the claim, we isolate the three-site block
indexed by $\mathcal J$ by making all other diagonal
entries large. This produces three bounded eigenvalues
converging to those of the block. We will show that
their first and second derivatives also converge,
while the derivatives of the remaining eigenvalues
vanish.

Since the coordinates of $c$ are not all equal, there
exists $1\le j_0\le n-2$ such that
$(c_{j_0},c_{j_0+1},c_{j_0+2})$ is not constant.
Keep $\mathbf y=(x_r,x_{r+1},x_{r+2})$ free and set
\[
x_m=Mv_m,\qquad m\notin\mathcal J,
\]
where $M>0$ is a scaling parameter and the numbers
$v_m$ are distinct and nonzero, with exactly $j_0-1$
of them negative.
After reordering the coordinates according to
$\mathcal J$ and its complement, the matrix becomes
\[
H_M(\mathbf y)=
\begin{pmatrix}
A(\mathbf y)&B\\
B^*&MV+C
\end{pmatrix},
\qquad
V=\operatorname{diag}(v_m:m\notin\mathcal J),
\]
where $A(\mathbf y)$ is the three-site Jacobi block
and $B,C$ are fixed matrices, with $C$ having zero
diagonal.
Define the corresponding three-site statistic
\[
g(\mathbf y)
=\sum_{i=1}^3c_{j_0+i-1}\lambda_i(A(\mathbf y)).
\]
Its coefficients are not all equal. The three-site
case therefore provides a point $\mathbf y_0$
and an ordering $(i,j,\ell)$ of $\mathcal J$ such that
\[
\Gaux_{i;j\ell}[g](\mathbf y_0)\ne0,
\]
where derivatives are taken with respect to the
original coordinates indexed by $\mathcal J$.

We first establish convergence of the bounded
eigenvalues through second derivatives.
Let $z\in\mathbb C$ denote the spectral parameter
and set $\epsilon=M^{-1}$. Since $V$ is invertible,
\[
\begin{aligned}
&(z\mathrm{Id}_{n-3}-MV-C)^{-1}=
\epsilon(\epsilon z\mathrm{Id}_{n-3}
-V-\epsilon C)^{-1}
\longrightarrow0
\end{aligned}
\]
as $\epsilon\to0$, locally uniformly for bounded $z$.

For bounded $z$ and sufficiently large $M$, the exterior block $z\mathrm{Id}_{n-3}-MV-C$ is invertible. We can therefore eliminate the exterior coordinates by taking the Schur complement:
\[
\begin{aligned}
\det(z\mathrm{Id}_n-H_M(\mathbf y))
&=\det(z\mathrm{Id}_{n-3}-MV-C)\\
&\quad\times
\det\left(
z\mathrm{Id}_3-A(\mathbf y)
-B(z\mathrm{Id}_{n-3}-MV-C)^{-1}B^*
\right).
\end{aligned}
\]
Thus the bounded eigenvalues are determined by
the second determinant.
To control their derivatives, define
\[
\mathcal F(\epsilon,z,\mathbf y)
=
\det\left(
z\mathrm{Id}_3-A(\mathbf y)
-\epsilon B
(\epsilon z\mathrm{Id}_{n-3}-V-\epsilon C)^{-1}B^*
\right).
\]
Since $V$ is invertible, $\mathcal F$ extends
holomorphically to $\epsilon=0$, where
\[
\mathcal F(0,z,\mathbf y)
=\det(z\mathrm{Id}_3-A(\mathbf y)).
\]
Denote the ordered eigenvalues of $A(\mathbf y)$
by $\lambda_i^A(\mathbf y)$, where $1\le i\le3$.
Their simplicity gives
\[
\begin{aligned}
&\partial_z\mathcal F
(0,\lambda_i^A(\mathbf y_0),\mathbf y_0)=
\prod_{\substack{1\le j\le3\\j\ne i}}
\bigl(\lambda_i^A(\mathbf y_0)
-\lambda_j^A(\mathbf y_0)\bigr)\ne0.
\end{aligned}
\]
The holomorphic implicit function theorem therefore
provides functions $\Lambda_i(\epsilon,\mathbf y)$,
holomorphic near $(0,\mathbf y_0)$, satisfying
\[
\mathcal F(\epsilon,\Lambda_i(\epsilon,\mathbf y),
\mathbf y)=0,
\qquad
\Lambda_i(0,\mathbf y)=\lambda_i^A(\mathbf y).
\]
For real $\mathbf y$ near $\mathbf y_0$ and sufficiently
large $M$, the exterior determinant does not vanish for bounded $z$.
Thus the three simple roots above give distinct bounded eigenvalues
of $H_M(\mathbf y)$. Ordered among themselves, these are
\[
\lambda_i^{(M)}(\mathbf y)
=\Lambda_i(M^{-1},\mathbf y),
\qquad 1\le i\le3.
\]
Their holomorphic dependence on $\epsilon$ implies
\[
\lambda_i^{(M)}(\mathbf y)
\longrightarrow\lambda_i^A(\mathbf y)
\]
locally uniformly on a fixed complex neighborhood
$\mathcal U$ of $\mathbf y_0$.
By Cauchy's integral formula, on every smaller
polydisc $\mathcal U'$ whose closure is contained in $\mathcal U$,
\[
\sup_{\mathbf y\in\mathcal U'}
\left|
\partial_{\mathbf y}^{\beta}
\bigl(
\lambda_i^{(M)}(\mathbf y)
-\lambda_i^A(\mathbf y)
\bigr)
\right|
\longrightarrow0,
\qquad |\beta|\le2.
\]
We next show that the remaining $n-3$ eigenvalues
have vanishing first and second derivatives in
the limit. Consider the rescaled matrix
\[
\epsilon H_M(\mathbf y)=
\begin{pmatrix}
\epsilon A(\mathbf y)&\epsilon B\\
\epsilon B^*&V+\epsilon C
\end{pmatrix}.
\]
At $\epsilon=0$, each $v_m$, with
$m\notin\mathcal J$, is a simple eigenvalue.
Let $\Lambda_m^{\mathrm{ext}}(\epsilon,\mathbf y)$
denote its analytic eigenvalue branch.
By Kato's perturbation theorem
\cite[Chapter~II, \S1.4 and Theorem~5.4]{Kato95},
its derivative at $\epsilon=0$ equals the
corresponding diagonal entry of $C$, which is zero.
Consequently,
\[
\Lambda_m^{\mathrm{ext}}(\epsilon,\mathbf y)
=v_m+O(\epsilon^2),
\]
locally uniformly for complex $\mathbf y$
near $\mathbf y_0$.
Returning to the original matrix, define
\[
\lambda_m^{\mathrm{ext},(M)}(\mathbf y)
=M\Lambda_m^{\mathrm{ext}}(M^{-1},\mathbf y).
\]
Then
\[
\lambda_m^{\mathrm{ext},(M)}(\mathbf y)
=Mv_m+O(M^{-1}).
\]
The error is locally uniform for complex $\mathbf y$
near $\mathbf y_0$. Cauchy's estimates therefore give
\[
\partial_{\mathbf y}^{\beta}
\lambda_m^{\mathrm{ext},(M)}(\mathbf y)
=O(M^{-1}),
\qquad 1\le|\beta|\le2.
\]
The three bounded branches and the $n-3$ distinct
exterior branches exhaust the spectrum for real $\mathbf y$ near
$\mathbf y_0$ and sufficiently large $M$.
By construction, exactly $j_0-1$ of the numbers
$v_m$ are negative. Thus $j_0-1$ exterior eigenvalues
tend to $-\infty$, while the remaining exterior
eigenvalues tend to $+\infty$. Consequently,
the three bounded eigenvalues occupy ranks
$j_0,j_0+1,j_0+2$ for sufficiently large $M$.
Their coefficients in $F_c$ are precisely
$c_{j_0},c_{j_0+1},c_{j_0+2}$.

Let $\mathbf x_M(\mathbf y)$ denote the diagonal vector
whose entries indexed by $\mathcal J$ are $\mathbf y$
and whose remaining entries are $Mv_m$.
The convergence of the bounded eigenvalue branches
and the vanishing derivatives of the exterior branches give
\[
\partial_{\mathbf y}^{\beta}
F_c(\mathbf x_M(\mathbf y))
\longrightarrow
\partial_{\mathbf y}^{\beta}g(\mathbf y),
\qquad 1\le|\beta|\le2,
\]
locally uniformly near $\mathbf y_0$.
Since the auxiliary function depends only on first and
second derivatives, we conclude that
\[
\Gaux_{i;j\ell}[F_c]
\bigl(\mathbf x_M(\mathbf y_0)\bigr)
\longrightarrow
\Gaux_{i;j\ell}[g](\mathbf y_0)\ne0.
\]
For every sufficiently large $M$,
$\Gaux_{i;j\ell}[F_c](\mathbf x_M(\mathbf y_0))\ne0$.
Hence $\Gaux_{i;j\ell}[F_c]\not\equiv0$ for the prescribed consecutive three-site block.

\noindent\emph{Absolute continuity.}
With the nondegeneracy established, we can apply Theorem~\ref{thm:balanced-collision-products}(i)
to the chosen three diagonal variables.

When $n=3$, the identity
\eqref{analytic:trivariate-balanced-determinant}
shows that the square of
$\Gaux_{i;j\ell}[F_c]$ is the diagonal
Monge--Amp\`ere determinant for the ordering
$(j,\ell,i)$. Since this determinant is nonzero
at some point and the three Frostman exponents
sum to more than $2$, the theorem gives
\[
(F_c)_*(\nu_1\otimes\nu_2\otimes\nu_3)\ll\Leb.
\]
Suppose $n>3$. Group the coordinates indexed by $\mathcal J$ into
$\mathbf y$ and the remaining coordinates into $\mathbf w$, in increasing
index order within each group. We write $F_c(\mathbf y,\mathbf w)$
and the auxiliary functions of $F_c$ in these reordered coordinates,
retaining the original site indices for derivatives. Set
\[
\mu=\bigotimes_{m\in\mathcal J}\nu_m,
\qquad
\mu_{\mathrm{out}}
=\bigotimes_{m\notin\mathcal J}\nu_m.
\]
The preceding claim provides a point
$(\mathbf y_*,\mathbf w_*)$ such that $\Gaux_{i;j\ell}[F_c](\mathbf y_*,\mathbf w_*)\ne0$.
Fix $\mathbf y_*$. The function
$\mathbf w\mapsto\Gaux_{i;j\ell}[F_c](\mathbf y_*,\mathbf w)$ is real
analytic on $\R^{n-3}$ and nonzero at $\mathbf w_*$. Since every exterior
law is nonatomic, Lemma~\ref{lem:analytic-product-zero-mass} shows that
this function is nonzero for $\mu_{\mathrm{out}}$-almost every $\mathbf w$.
For each such $\mathbf w$, the three-variable map
\[
F_{\mathbf w}(\mathbf y)=F_c(\mathbf y,\mathbf w)
\]
satisfies the Monge--Amp\`ere condition.
Theorem~\ref{thm:balanced-collision-products}(i)
therefore gives
\[
(F_{\mathbf w})_*\mu\ll\Leb.
\]
Finally, Fubini's theorem yields, for every
Lebesgue-null Borel set $N\subset\R$,
\[
(F_c)_*\left(\bigotimes_{m=1}^n\nu_m\right)(N)
=\int (F_{\mathbf w})_*\mu(N)\,
\dd\mu_{\mathrm{out}}(\mathbf w)=0.
\]
Thus $F_c(\mathbf X)$ has an absolutely continuous
law, completing the proof.
\end{proof}

\subsection{Spectral averaging for finite Jacobi matrices}
\label{spectral:proof-finite}

For Theorem~\ref{spectral:finite}, it suffices to prove absolute continuity of the ordered eigenvalue laws. When $k=2$, we verify the bivariate nondegeneracy
condition by examining the eigenvalue level curves.
When $k\ge3$, we use the three-term eigenvector recurrence
to verify the Monge--Amp\`ere condition for a suitable
balanced partition of the collision variables.

\begin{proof}[Proof of Theorem~\ref{spectral:finite}]
Recall that $2\le k\le n$ and that the first $k$ diagonal
entries are independent random variables $X_1,\ldots,X_k$
with compactly supported $s_i$-Frostman probability laws
$\nu_i$, satisfying
\[
0<s_i\le1\quad\text{for every }1\le i\le k,
\qquad
\sum_{i=1}^k s_i>\frac{k+1}{2}.
\]
The remaining diagonal entries
$\mathbf b=(b_{k+1},\ldots,b_n)\in\R^{n-k}$ and the
nonzero off-diagonal entries $t_1,\ldots,t_{n-1}$ are fixed.

For $\mathbf x=(x_1,\ldots,x_k)\in\R^k$, write
\[
H_{\mathbf b}(\mathbf x)
=
\operatorname{diag}(x_1,\ldots,x_k,b_{k+1},\ldots,b_n)
+\sum_{i=1}^{n-1}t_i
(e_ie_{i+1}^*+e_{i+1}e_i^*).
\]
Let $\boldsymbol\nu=\bigotimes_{i=1}^k\nu_i$ be the joint
law of the random diagonal entries. Denote the ordered
eigenvalues of $H_{\mathbf b}(\mathbf x)$ by
$\lambda_1^{\mathbf b}(\mathbf x)<\cdots<
\lambda_n^{\mathbf b}(\mathbf x)$, and set
\[
\rho_j=(\lambda_j^{\mathbf b})_*\boldsymbol\nu,
\qquad 1\le j\le n.
\]
The connection between these eigenvalue laws and the
averaged scalar spectral measures is given by
\eqref{spectral:scalar-spectral-weights}. Indeed, given
$u\in\mathbb C^n$ and a Borel set $B\subset\R$, we have
\[
\begin{aligned}
\overline\sigma_u(B)
&=\int_{\R^k}\sum_{j=1}^n
\mathbf1_B(\lambda_j^{\mathbf b}(\mathbf x))
|\langle u,v_j^{\mathbf b}(\mathbf x)\rangle|^2
\,\dd\boldsymbol\nu(\mathbf x)\le\|u\|^2\sum_{j=1}^n\rho_j(B),
\end{aligned}
\]
where $v_j^{\mathbf b}(\mathbf x)$ is the unit eigenvector
associated with $\lambda_j^{\mathbf b}(\mathbf x)$,
as in Section~\ref{spectral:random-jacobi}.
Thus it suffices to establish $\rho_j\ll\Leb$ for every
$1\le j\le n$.
The three-term recurrence determines every eigenvector
from its first coordinate, which cannot vanish for a nonzero
eigenvector. Thus each eigenspace is one-dimensional. The analytic
implicit function theorem gives local analytic eigenvalue branches,
and simplicity fixes their order, so each $\lambda_j^{\mathbf b}$
is real analytic on $\R^k$.

\medskip
\noindent\textbf{Case 1: Two random sites ($k=2$).}
Write $\mathbf x=(x,y)$. Since $s_1+s_2>3/2$,
Lemma~\ref{spectral:bivariate-frostman-ac} reduces the
problem to verifying the nondegeneracy condition
\eqref{spectral:nondegeneracy} at some point for each
ordered eigenvalue map of $H_{\mathbf b}(x,y)$.

We study this condition along a level curve
\[
\{(x,y)\in\R^2:
\lambda_j^{\mathbf b}(x,y)=\tau\},
\]
where $\tau\in\R$ is a fixed spectral parameter.
The Schur complement gives an explicit parametrization
of these curves.

\medskip
\noindent\emph{The eigenvalue level curve.}
Suppose first that $n>2$. Let $D$ be the
$(n-2)\times(n-2)$ principal submatrix on sites
$3,\ldots,n$, and let $e$ be the first standard
basis vector of $\R^{n-2}$. Relative to the
decomposition into the first two sites and the tail,
\[
H_{\mathbf b}(x,y)=
\begin{pmatrix}
x&t_1&0\\
t_1&y&t_2e^*\\
0&t_2e&D
\end{pmatrix}.
\]
For $n>2$, define the tail function
$a:\R\setminus\operatorname{spec}(D)\to\R$ by
\begin{equation}\label{spectral:tail-function}
a(\tau)=\tau+t_2^2
\left\langle
(D-\tau\mathrm{Id}_{n-2})^{-1}e,e
\right\rangle,
\qquad \tau\notin\operatorname{spec}(D),
\end{equation}
where $\operatorname{spec}(D)$ denotes the set of
eigenvalues of the tail matrix $D$.
When $n=2$, there is no tail; we take
$\operatorname{spec}(D)=\varnothing$ and set $a(\tau)=\tau$. To verify the nondegeneracy condition, we seek a simpler
equation for the eigenvalue maps from which their derivatives
can be computed explicitly. This motivates the introduction
of the function $a(\tau)$.
\\
We record the positivity of $a'(\tau)$, which will
be used to verify the nondegeneracy condition.
Differentiating \eqref{spectral:tail-function} gives
\begin{equation}\label{spectral:positive-derivative}
g(\tau):=a'(\tau)
=1+t_2^2
\left\|(D-\tau\mathrm{Id}_{n-2})^{-1}e\right\|^2
>0.
\end{equation}
When $n=2$, the tail term is absent.
\\
Set $T=t_1^2>0$. The Schur complement identity gives,
whenever $\tau\notin\operatorname{spec}(D)$,
\[
\begin{aligned}
\det(H_{\mathbf b}(x,y)-\tau\mathrm{Id}_n)
&=\det(D-\tau\mathrm{Id}_{n-2})\times
\det\begin{pmatrix}
x-\tau&t_1\\
t_1&y-a(\tau)
\end{pmatrix}.
\end{aligned}
\]
Thus $\tau$ is an eigenvalue precisely when
\begin{equation}\label{spectral:level-curve}
(x-\tau)(y-a(\tau))=T.
\end{equation}
The same equation follows directly from the
characteristic polynomial when $n=2$.
\\
Fix $\tau\notin\operatorname{spec}(D)$.
The level curve \eqref{spectral:level-curve} is
parametrized by
\[
x=\tau+r,\qquad
y=a(\tau)+\frac{T}{r},
\qquad r\in\R\setminus\{0\}.
\]
This parametrization allows us to let $r$ approach
zero while keeping the eigenvalue $\tau$ fixed.

\medskip
\noindent\emph{Nondegeneracy along the level curve.}
For each $r_0\in\R\setminus\{0\}$, let
$(x_0,y_0)=(\tau+r_0,a(\tau)+T/r_0)$ be the corresponding
point on \eqref{spectral:level-curve}. Since the eigenvalues of
$H_{\mathbf b}(x_0,y_0)$ are simple, there is a unique
index $j\in\{1,\ldots,n\}$ such that
\[
\lambda_j^{\mathbf b}(x_0,y_0)=\tau.
\]
For this point, set $j=j(r_0)$ and define
\[
f:\R^2\to\R,\qquad
f(x,y)=\lambda_j^{\mathbf b}(x,y).
\]
Thus $f(x_0,y_0)=\tau$, and $f$ is real analytic on $\R^2$.
Write $f_1=\partial_x f$ and $f_2=\partial_y f$.
By continuity, $f$ remains outside
$\operatorname{spec}(D)$ on a neighborhood of
$(x_0,y_0)$. Equation~\eqref{spectral:level-curve}
therefore gives
\[
(x-f(x,y))(y-a(f(x,y)))=T.
\]
On this neighborhood, write
\[
r=r(x,y)=x-f(x,y),
\qquad
W=W(x,y)=g(f(x,y))r(x,y)^2+T.
\]
We suppress the dependence on $(x,y)$ below.
Differentiating the eigenvalue equation with respect
to $x$ and $y$ yields
\[
\begin{aligned}
(1-f_1)(y-a(f))-rg(f)f_1&=0,\\
-f_2(y-a(f))+r(1-g(f)f_2)&=0.
\end{aligned}
\]
Substituting $y-a(f)=T/r$, we obtain
\[
(1-f_1)T=g(f)r^2f_1,
\qquad
Tf_2=r^2(1-g(f)f_2).
\]
Solving these identities for $f_1$ and $f_2$ gives
\[
f_1=\frac{T}{W},
\qquad
f_2=\frac{r^2}{W}.
\]
In particular, the ratio of the first derivatives is
\[
R_f=\frac{f_1}{f_2}=\frac{T}{r^2}.
\]
Since $\partial_2r=-f_2$, we also have
\[
\partial_2R_f
=-\frac{2T}{r^3}\partial_2r
=\frac{2T}{rW}.
\]
We have therefore established
\begin{equation}\label{spectral:first-derivatives}
f_1=\frac{T}{W},
\qquad
f_2=\frac{r^2}{W},
\qquad
R_f=\frac{T}{r^2},
\qquad
\partial_2R_f=\frac{2T}{rW}.
\end{equation}
Since $T>0$, $g(f)>0$, and $r\ne0$, the quantities
$f_1$, $f_2$, and $\partial_2R_f$ are nonzero.
\\
It remains to compute
$\mathcal K=\partial_{12}\log|R_f|$.
Using $\partial_1r=1-f_1$, we obtain
\[
\partial_1\log|R_f|
=-\frac{2(1-f_1)}{r}
=-\frac{2g(f)r}{W}.
\]
To differentiate this expression with respect to $y$,
observe that
\[
\partial_2(g(f)r)
=\bigl(a''(f)r-g(f)\bigr)f_2
\]
and
\[
\partial_2W
=\bigl(a''(f)r^2-2g(f)r\bigr)f_2.
\]
The quotient rule then gives
\[
\begin{aligned}
\mathcal K
&=-2\partial_2\left(\frac{g(f)r}{W}\right)=-\frac{2f_2}{W^2}
\left[
\bigl(a''(f)r-g(f)\bigr)W
-g(f)r\bigl(a''(f)r^2-2g(f)r\bigr)
\right].
\end{aligned}
\]
Substituting $f_2=r^2/W$ and $W=g(f)r^2+T$ yields
\begin{equation}\label{spectral:curvature}
\mathcal K
=\frac{2r^2}{W^3}
\left(Tg(f)-Ta''(f)r-g(f)^2r^2\right).
\end{equation}
On the component parametrized by $r$ with the same sign as $r_0$,
the chosen ordered eigenvalue satisfies $f(x,y)=\tau$. Hence
\[
Tg(f)-Ta''(f)r-g(f)^2r^2
=
Tg(\tau)-Ta''(\tau)r-g(\tau)^2r^2
\longrightarrow Tg(\tau)>0
\]
as $r\to0$. Consequently, $\mathcal K\ne0$ for
all sufficiently small nonzero $r$. Together with
\eqref{spectral:first-derivatives}, this gives
\[
f_1f_2\,\partial_2R_f\,\mathcal K\ne0
\]
at the corresponding points
\[
(x,y)=\left(\tau+r,a(\tau)+\frac{T}{r}\right).
\]
For either sign of $r_0$, when $0<|r_0|$ is sufficiently small, this gives nondegeneracy at $(x_0,y_0)$ for the branch $\lambda_{j(r_0)}^{\mathbf b}$ through that point. The branch can differ between the two signs.

\medskip
\noindent\emph{Identifying the ordered eigenvalue.}
We must still show that the points constructed above
give nondegeneracy for every ordered eigenvalue map.

Let $m$ be the number of eigenvalues of $D$ strictly
below $\tau$. Write $n_-(L)$ for the number of strictly
negative eigenvalues of a real symmetric matrix $L$,
counted with multiplicity.
\\
For $n>2$, Haynsworth's inertia additivity formula~\cite{Haynsworth68},
applied to the invertible block $D-\tau\mathrm{Id}_{n-2}$ and its
Schur complement, gives
\[
n_-(H_{\mathbf b}(x,y)-\tau\mathrm{Id}_n)
=m+n_-\begin{pmatrix}r&t_1\\ t_1&T/r\end{pmatrix}.
\]
The same formula holds for $n=2$ with $m=0$.
The displayed $2\times2$ matrix has determinant zero and trace $r+T/r$.
Its only nonzero eigenvalue has the sign of $r$, since $T>0$.
It follows that
\[
n_-(H_{\mathbf b}(x,y)-\tau\mathrm{Id}_n)
=m+\mathbf1_{\{r<0\}}.
\]
Since the eigenvalues of $H_{\mathbf b}(x,y)$ are simple,
the index of the eigenvalue $\tau$ is one more than the
number of eigenvalues strictly below it. Consequently,
\[
\tau=
\begin{cases}
\lambda_{m+1}^{\mathbf b}(x,y),&r>0,\\
\lambda_{m+2}^{\mathbf b}(x,y),&r<0.
\end{cases}
\]
The gaps of the tail spectrum realize every count
$m=0,\ldots,n-2$. When $n=2$, the only gap is $\R$
and $m=0$. Choosing an energy in each gap and taking
both signs of sufficiently small nonzero $r$ therefore
gives a nondegenerate point for every ordered eigenvalue
map $\lambda_j^{\mathbf b}$, with $1\le j\le n$.
\\
Since $s_1+s_2>3/2$, Lemma~\ref{spectral:bivariate-frostman-ac}
now gives
\[
\rho_j=(\lambda_j^{\mathbf b})_*
(\nu_1\otimes\nu_2)\ll\Leb
\]
for every $1\le j\le n$.

\medskip
\noindent\textbf{Case 2: Three or more random sites.}
We apply Theorem~\ref{thm:balanced-collision-products}(i)
to each ordered eigenvalue map. The main step is to find
a balanced partition of the collision variables for which
the Monge--Amp\`ere determinant is nonzero at some point.

Fix an ordered eigenvalue map $f:\R^k\to\R$ and let
$v(\mathbf x)$ be its real unit eigenvector with positive first coordinate.
For $1\le i\le n$, let $e_i$ be the $i$th standard basis vector
of $\R^n$ and set $v_i(\mathbf x)=\langle v(\mathbf x),e_i\rangle$. Let $D$ be the fixed principal
submatrix on sites $k+1,\ldots,n$, with empty spectrum
when $k=n$.

To express the eigenvector in terms of the diagonal entries, we use the left Dirichlet solution. Given a spectral parameter $\tau\in\R$ and
$\mathbf x\in\R^k$, define
\begin{equation}\label{spectral:recurrence}
\begin{gathered}
p_0(\tau,\mathbf x)=0,\qquad p_1(\tau,\mathbf x)=1,\\
p_{i+1}(\tau,\mathbf x)
=
\frac{
(\tau-(H_{\mathbf b}(\mathbf x))_{ii})p_i(\tau,\mathbf x)
-t_{i-1}p_{i-1}(\tau,\mathbf x)
}{t_i},
\qquad 1\le i<n,
\end{gathered}
\end{equation}
where the term involving $t_0$ is omitted.
At $\tau=f(\mathbf x)$, the recurrence gives
\[
p_i(f(\mathbf x),\mathbf x)=\frac{v_i(\mathbf x)}
{v_1(\mathbf x)}\qquad\text{for every }1\le i\le n.
\]
Since $\|v(\mathbf x)\|=1$ and $v_1(\mathbf x)>0$,
we obtain
\[
v_i(\mathbf x)=
\frac{p_i(f(\mathbf x),\mathbf x)}
{\left(
\sum_{j=1}^n p_j(f(\mathbf x),\mathbf x)^2
\right)^{1/2}}.
\]
The functions $p_i$ are polynomials in $(\tau,\mathbf x)$,
and the denominator is positive because $p_1=1$.
Thus $v$ depends real analytically on $\mathbf x$.
The Hellmann--Feynman formula~\cite[Chapter~II,Theorem~5.4]{Kato95} gives
\[
\partial_i f(\mathbf x)
=\langle e_ie_i^*v,v\rangle
=|\langle v,e_i\rangle|^2
=v_i^2
\]
for every $1\le i\le k$.
We use this identity to express the derivatives of
the collision phase in terms of the eigenvector.

\medskip
\noindent\emph{Choosing a nondegenerate point.}
We first find a point $\mathbf x^0$ at which
\[
v_i(\mathbf x^0)\ne0
\quad\text{for every }1\le i\le k,
\qquad
f(\mathbf x^0)\notin\operatorname{spec}(D).
\]
The recurrence, applied from either endpoint, shows
that $v_1$ and $v_n$ never vanish. For each
$2\le i\le\min\{k,n-1\}$, let $D_i$ be the principal
submatrix of $H_{\mathbf b}(\mathbf x)$ on sites
$i+1,\ldots,n$.
Write $v_{>i}=(v_{i+1},\ldots,v_n)^T$.
The eigenvalue equation on these sites is
\[
(D_i-f(\mathbf x)\mathrm{Id}_{n-i})v_{>i}
=-t_i v_i e_1^{(n-i)},
\]
where $e_1^{(n-i)}$ is the first standard basis
vector of $\R^{n-i}$. Since $v_n\ne0$, we have
$v_{>i}\ne0$. In particular,
\[
v_i(\mathbf x)=0
\quad\Longrightarrow\quad
f(\mathbf x)\in\operatorname{spec}(D_i).
\]
Fix $x_2,\ldots,x_k$. The matrices $D$ and $D_i$
are then independent of $x_1$. On the other hand,
\[
\partial_1f(\mathbf x)=v_1(\mathbf x)^2>0,
\]
so the map $x_1\mapsto f(x_1,\ldots,x_k)$ is
strictly increasing.
Each energy in the finite set
\[
\operatorname{spec}(D)
\cup
\bigcup_{i=2}^{\min\{k,n-1\}}
\operatorname{spec}(D_i)
\]
can therefore be attained for at most one value of
$x_1$. Choosing $x_1$ outside the resulting finite
set of preimages gives a point $\mathbf x^0$ satisfying the
required conditions.

\medskip
\noindent\emph{The balanced partition.}
Define the collision phase $\Phi:\R^k\times\R^k\to\R$ by
\[
\Phi(\mathbf x,\mathbf y)=f(\mathbf x)-f(\mathbf y).
\]
We will verify \eqref{shell:local-bordered-condition}
at $(\mathbf x,\mathbf y)=(\mathbf x^0,\mathbf x^0)$.

The three-term recurrence dictates the partition: two consecutive fixed sites in one copy force the eigenvalue variation to vanish, while alternating fixed sites determine the remaining eigenvector coordinates.

Set
\[
A=\{2j-1:1\le j\le\lfloor k/2\rfloor\},
\qquad
C=\{1\}\cup
\{2j:1\le j\le\lfloor(k-1)/2\rfloor\}.
\]
For any index set $J\subset\{1,\ldots,k\}$, write
$\mathbf x_J=(x_j)_{j\in J}$ and
$\mathbf y_J=(y_j)_{j\in J}$, with coordinates
listed in increasing order.

Use the notation of Section~\ref{sec:fourier-shell}: set
$\sigma=(\sigma_L,\sigma_R)$, where $\sigma_L=A\cup\{k+j:j\in C\}$
and $\sigma_R=\{1,\ldots,2k\}\setminus\sigma_L$. The indices refer to the concatenated vector $(\mathbf x,\mathbf y)$. The coordinate vectors are
\[
x_{\sigma_L}=(\mathbf x_A,\mathbf y_C),
\qquad
x_{\sigma_R}=(\mathbf x_{A^c},\mathbf y_{C^c}),
\]
where complements are taken in $\{1,\ldots,k\}$.
Since
\[
|A|+|C|
=\lfloor k/2\rfloor+1+\lfloor(k-1)/2\rfloor=k,
\]
both groups contain exactly $k$ variables.
For example, when $k=4$, the partition is
\[
x_{\sigma_L}=(x_1,x_3,y_1,y_2),
\qquad
x_{\sigma_R}=(x_2,x_4,y_3,y_4).
\]
Let $(x_{\sigma_L}^0,x_{\sigma_R}^0)=\pi_\sigma(\mathbf x^0,\mathbf x^0)$
and write $\Phi_\sigma=\Phi\circ\pi_\sigma^{-1}$. With $x_{\sigma_L}=x_{\sigma_L}^0$ fixed, consider
the map
\[
\mathcal T(x_{\sigma_R},\theta)
=
\bigl(
\theta\nabla_{x_{\sigma_L}}\Phi_\sigma(x_{\sigma_L}^0,x_{\sigma_R}),
\Phi_\sigma(x_{\sigma_L}^0,x_{\sigma_R})
\bigr),
\qquad
\mathcal T:\R^k\times(\R\setminus\{0\})
\to\R^{k+1}.
\]
The derivative with respect to $(x_{\sigma_R},\theta)$ is the matrix
below. Its entries are formed by differentiating $\Phi_\sigma$ before
evaluating at $(x_{\sigma_L}^0,x_{\sigma_R})$.
\[
D\mathcal T=
\begin{pmatrix}
\theta D_{x_{\sigma_L}}D_{x_{\sigma_R}}\Phi_\sigma
&\nabla_{x_{\sigma_L}}\Phi_\sigma\\
(\nabla_{x_{\sigma_R}}\Phi_\sigma)^T&0
\end{pmatrix}.
\]
At $\theta=1$, this is precisely the matrix defining
the Monge--Amp\`ere determinant in
\eqref{shell:local-bordered-condition}.
It therefore suffices to show that
\[
\ker D\mathcal T(x_{\sigma_R}^0,\theta)=\{0\}.
\]

\medskip
\noindent\emph{Linearization of the eigenvalue equation.}
Let $(\dot x_{\sigma_R},\dot\theta)$ belong to the
kernel of $D\mathcal T(x_{\sigma_R}^0,\theta)$ for a fixed $\theta\in\R\setminus\{0\}$.
The corresponding variations of the two diagonal
vectors satisfy
\[
\dot x_i=0\quad\text{for every }i\in A,
\qquad
\dot y_i=0\quad\text{for every }i\in C.
\]
We write $\dot f$, $\dot H$, and $\dot v$ for the
induced variations in each copy.
Since the last component of $\mathcal T$ is $\Phi_\sigma$,
the kernel condition gives
\[
\dot f^{(\mathbf x)}-\dot f^{(\mathbf y)}=0.
\]
Denote their common value by $\dot \tau$.
At the diagonal point, both copies have the same
matrix, eigenvalue, and eigenvector:
\[
H_0=H_{\mathbf b}(\mathbf x^0),
\qquad
\tau=f(\mathbf x^0),
\qquad
v=v(\mathbf x^0),
\qquad
H_0v=\tau v.
\]
Their variations $\dot H$ and $\dot v$ may differ.
The same derivative formula gives
\[
\nabla_{x_{\sigma_L}}\Phi_\sigma(\pi_\sigma(\mathbf x,\mathbf y))
=
\left(
(v_i(\mathbf x)^2)_{i\in A},
-(v_i(\mathbf y)^2)_{i\in C}
\right).
\]
Since the variation of
$\theta\nabla_{x_{\sigma_L}}\Phi_\sigma$ vanishes,
\[
\dot\theta\,v_i^2+2\theta v_i\dot v_i=0
\]
for $i\in A$ in the first copy and for $i\in C$ in the second.
Using $v_i\ne0$, we obtain
\[
\dot v_i=hv_i,
\qquad
h=-\frac{\dot\theta}{2\theta}.
\]
In particular, $h$ is common to both copies.
In each copy, define $w=\dot v-hv$.
Then
\[
w_i=\dot H_{ii}=0
\]
at the sites indexed by $A$ in the first copy and by $C$ in the second. The remaining fixed diagonal entries have $\dot H_{ii}=0$, but those coordinates of $w$ are determined below.
Differentiating the eigenvalue equation
$H_{\mathbf b}(\mathbf x)v(\mathbf x)=f(\mathbf x)v(\mathbf x)$ at the base point in each copy gives
\[
\dot H v+H_0\dot v
=\dot \tau v+\tau\dot v.
\]
Since $(H_0-\tau\mathrm{Id}_n)v=0$, subtracting
$hv$ from $\dot v$ does not change the left-hand side
of the resulting equation. Consequently,
\begin{equation}\label{spectral:linearized-recurrence}
(H_0-\tau\mathrm{Id}_n)w
=(\dot \tau\mathrm{Id}_n-\dot H)v.
\end{equation}
Moreover, differentiating $\|v\|^2=1$ gives
$\langle v,\dot v\rangle=0$, so
\begin{equation}\label{spectral:linearized-normalization}
\langle v,w\rangle
=\langle v,\dot v-hv\rangle=-h.
\end{equation}
We now show that these identities force
$\dot \tau=h=0$ and $w=0$ in both copies.

\medskip
\noindent\emph{Vanishing of the eigenvalue variation.}
The first two sites are fixed in the second copy,
since $\{1,2\}\subset C$. Thus
\[
w_1=w_2=0,
\qquad
\dot H_{11}=\dot H_{22}=0.
\]
The first row of
\eqref{spectral:linearized-recurrence} is
\[
((H_0)_{11}-\tau)w_1+t_1w_2
=(\dot \tau-\dot H_{11})v_1.
\]
It follows that $\dot \tau v_1=0$.
Since $v_1>0$, we obtain $\dot \tau=0$ in both copies.
At a site $i\in A$ in the first copy or $i\in C$ in the second copy, the $i$th row of
\eqref{spectral:linearized-recurrence} now reads
\[
t_{i-1}w_{i-1}
+((H_0)_{ii}-\tau)w_i
+t_iw_{i+1}
=-\dot H_{ii}v_i.
\]
Since $w_i=\dot H_{ii}=0$ at that site,
\begin{equation}\label{spectral:fixed-row-recurrence}
t_{i-1}w_{i-1}+t_iw_{i+1}=0,
\end{equation}
with the first term omitted when $i=1$.

\medskip
\noindent\emph{Propagation through the first $k$ sites.}
In the first copy, the fixed coordinates satisfy
\[
w_{2j-1}=0
\qquad\text{for every }
1\le j\le\lfloor k/2\rfloor.
\]
Applying \eqref{spectral:fixed-row-recurrence}
successively at these sites gives
\[
t_1w_2=0,\qquad
t_{2j-2}w_{2j-2}+t_{2j-1}w_{2j}=0
\quad\text{for }2\le j\le\lfloor k/2\rfloor .
\]
Since every $t_i$ is nonzero, it follows that
\[
w_1=\cdots=w_{2\lfloor k/2\rfloor}=0.
\]
In the second copy, the fixed coordinates satisfy
\[
w_1=0,
\qquad
w_{2j}=0
\quad\text{for every }
1\le j\le\lfloor(k-1)/2\rfloor.
\]
Applying \eqref{spectral:fixed-row-recurrence}
successively at the even sites gives
\[
t_{2j-1}w_{2j-1}+t_{2j}w_{2j+1}=0
\quad\text{for }1\le j\le\lfloor(k-1)/2\rfloor .
\]
Hence $w_1=\cdots=w_{2\lfloor(k-1)/2\rfloor+1}=0$.
We conclude that
\[
\begin{array}{c|cc}
&\text{First copy}&\text{Second copy}\\ \hline
k\text{ even}
&w_1=\cdots=w_k=0
&w_1=\cdots=w_{k-1}=0\\
k\text{ odd}
&w_1=\cdots=w_{k-1}=0
&w_1=\cdots=w_k=0.
\end{array}
\]
Thus all the first $k$ coordinates vanish in one copy,
while only $w_k$ may remain nonzero in the other.

\medskip
\noindent\emph{Eliminating the remaining coordinate.}
Suppose first that $k<n$. Decompose
\[
v=(v_1,\ldots,v_k,v_{\mathrm{tail}}),
\qquad
w=(w_1,\ldots,w_k,w_{\mathrm{tail}}),
\]
where $v_{\mathrm{tail}},w_{\mathrm{tail}}\in\R^{n-k}$.
Let $e_1^{\mathrm{tail}}$ be the first standard
basis vector of $\R^{n-k}$.
Since the tail is fixed, its diagonal variations
vanish. The eigenvalue equation and
\eqref{spectral:linearized-recurrence}, with
$\dot \tau=0$, therefore give
\[
\begin{aligned}
(D-\tau\mathrm{Id}_{n-k})v_{\mathrm{tail}}
&=-t_kv_ke_1^{\mathrm{tail}},\\
(D-\tau\mathrm{Id}_{n-k})w_{\mathrm{tail}}
&=-t_kw_ke_1^{\mathrm{tail}}.
\end{aligned}
\]
By the choice of $\mathbf x^0$, the matrix
$D-\tau\mathrm{Id}_{n-k}$ is invertible. Hence
\[
\begin{aligned}
v_{\mathrm{tail}}
&=-t_kv_k(D-\tau\mathrm{Id}_{n-k})^{-1}
e_1^{\mathrm{tail}},\\
w_{\mathrm{tail}}
&=-t_kw_k(D-\tau\mathrm{Id}_{n-k})^{-1}
e_1^{\mathrm{tail}}.
\end{aligned}
\]
Since $v_k\ne0$, we obtain
\[
w_{\mathrm{tail}}
=\frac{w_k}{v_k}v_{\mathrm{tail}}.
\]
In the copy whose first $k$ coordinates vanish,
$w_k=0$, so $w_{\mathrm{tail}}=0$ as well.
Thus $w=0$ in that copy, and
\eqref{spectral:linearized-normalization} gives $h=0$.
In the other copy, using $w_1=\cdots=w_{k-1}=0$
and $h=0$, the same normalization identity gives
\[
\begin{aligned}
0
&=\langle v,w\rangle=v_kw_k+
\langle v_{\mathrm{tail}},w_{\mathrm{tail}}\rangle=\frac{w_k}{v_k}
\left(v_k^2+\|v_{\mathrm{tail}}\|^2\right).
\end{aligned}
\]
Since $v_k\ne0$, it follows that $w_k=0$.
Consequently, $w_{\mathrm{tail}}=0$ and hence
$w=0$ in both copies.
\\
When $k=n$, there is no tail. The copy with all
coordinates vanishing gives $h=0$, and the
normalization identity in the other copy reduces to
\[
0=v_kw_k.
\]
Again, $v_k\ne0$ implies $w_k=0$.
We have therefore established
\[
w=0,\qquad \dot \tau=0,\qquad h=0
\]
in both copies. Returning to
\eqref{spectral:linearized-recurrence}, we obtain $\dot H v=0$.
In the first and second copies, respectively, this
identity reads
\[
\operatorname{diag}(\dot x_1,\ldots,\dot x_k,0,\ldots,0)v=0
\]
and
\[
\operatorname{diag}(\dot y_1,\ldots,\dot y_k,0,\ldots,0)v=0.
\]
Since $v_i\ne0$ for every $1\le i\le k$, we conclude
that
\[
\dot x_i=\dot y_i=0
\qquad\text{for every }1\le i\le k.
\]
Moreover,
\[
h=-\frac{\dot\theta}{2\theta}=0
\quad\Longrightarrow\quad
\dot\theta=0.
\]
Thus the kernel vector $(\dot x_{\sigma_R},\dot\theta)$
vanishes, proving that $D\mathcal T$ is invertible.

\medskip
\noindent\emph{Absolute continuity of the eigenvalue laws.}
Taking $\theta=1$ in the preceding calculation gives
\[
\left.
\det
\begin{pmatrix}
D_{x_{\sigma_L}}D_{x_{\sigma_R}}\Phi_\sigma
&\nabla_{x_{\sigma_L}}\Phi_\sigma\\
(\nabla_{x_{\sigma_R}}\Phi_\sigma)^T&0
\end{pmatrix}
\right|_{(x_{\sigma_L},x_{\sigma_R})=\pi_\sigma(\mathbf x^0,\mathbf x^0)}
\ne0.
\]
This verifies the Monge--Amp\`ere condition
\eqref{analytic:balanced-product-condition} for
the balanced partition constructed above.
\\
Choose compact intervals $I_i$ containing
$\supp\nu_i$ for every $1\le i\le k$.
Since $f$ is real analytic on $\R^k$, the nondegenerate
point $\mathbf x^0$ need not belong to the product
$I_1\times\cdots\times I_k$.
Theorem~\ref{thm:balanced-collision-products}(i),
together with
\[
\sum_{i=1}^k s_i>\frac{k+1}{2},
\]
therefore gives $f_*\boldsymbol\nu\ll\Leb$.
Since $f$ was an arbitrary ordered eigenvalue map,
\[
\rho_j=(\lambda_j^{\mathbf b})_*\boldsymbol\nu
\ll\Leb
\qquad\text{for every }1\le j\le n.
\]
The spectral-weight bound at the start of the proof now yields
$\overline\sigma_u\ll\Leb$ for every $u\in\mathbb C^n$.
\end{proof}

\subsection{Spectral averaging for half-line Jacobi operators}
\label{spectral:proof-halfline}

Unlike in the finite-dimensional setting, where we study
individual eigenvalue maps, the proof for half-line Jacobi
operators proceeds by fixing the spectral parameter and
studying the dependence of the boundary ratio on the random
diagonal entries. We first establish a uniform density
estimate for the boundary ratio, from which we derive
spectral averaging estimates for finite restrictions.
Passing to the half-line limit and showing that the
exceptional energies carry no averaged spectral mass
will complete the proof.

\begin{proof}[Proof of Theorem~\ref{spectral:halfline}]
Recall that $m=2\ell\ge4$ and that the first $m$ diagonal
entries are independent random variables $X_1,\ldots,X_m$
with compactly supported $s_i$-Frostman probability laws
$\mu_i$, satisfying
\[
0<s_i\le1\quad\text{for every }1\le i\le m,
\qquad
S_m:=\sum_{i=1}^m s_i>\frac{m}{2}+1=\ell+1.
\]
The hopping coefficients $(t_i)_{i\ge1}$ are fixed, real,
nonzero, and uniformly bounded. The full diagonal sequence is bounded almost surely, and the entries beyond the
first $m$ sites are either deterministic or independent
of $X_1,\ldots,X_m$.

By independence, it suffices to prove the theorem for
each fixed bounded tail $\mathbf b=(b_i)_{i>m}$.
For $\mathbf x=(x_1,\ldots,x_m)\in\R^m$, let
$J(\mathbf x)$ be the Jacobi operator on $\ell^2(\mathbb N)$
given by
\[
J(\mathbf x)=
\begin{pmatrix}
x_1 & t_1 & 0 & \cdots & & & \\
t_1 & x_2 & t_2 & \ddots & & & \\
0 & t_2 & \ddots & \ddots & & & \\
\vdots & \ddots & \ddots & x_m & t_m & & \\
& & & t_m & b_{m+1} & t_{m+1} & \\
& & & & t_{m+1} & b_{m+2} & \ddots\\
& & & & & \ddots & \ddots
\end{pmatrix}.
\]
The dependence on the fixed tail $\mathbf b$ is
suppressed in the notation.

Let $\boldsymbol\mu=\bigotimes_{i=1}^m\mu_i$ be the joint
law of $X_1,\ldots,X_m$. Our goal is to establish absolute
continuity of the averaged scalar spectral measure
\[
\overline{\sigma}_u(B)
:=
\int_{\R^m}
\left\langle\mathbf1_B(J(\mathbf x))u,u\right\rangle
\,\dd\boldsymbol\mu(\mathbf x)
\]
for every $u\in\ell^2(\mathbb N)$ and Borel set $B\subset\R$.
Since $e_1$ is a cyclic vector for every $J(\mathbf x)$,
it suffices to establish
\[
\overline{\sigma}_{e_1}
=\int_{\R^m}\sigma_{J(\mathbf x),e_1}
\,\dd\boldsymbol\mu(\mathbf x)
\ll\Leb.
\]
On sets of uniformly bounded diagonal entries,
coordinatewise convergence implies strong operator
convergence of the corresponding Jacobi operators.
By polynomial approximation, the associated spectral
measures depend weakly continuously on the diagonal
entries. Since the diagonal is almost surely bounded,
the spectral measures form Borel kernels on a set of
full probability, justifying the applications of
Tonelli's theorem below.

\medskip
\noindent\textbf{Step 1: The boundary ratio and its derivatives.}
Write $d_i(\mathbf x)=x_i$ for $1\le i\le m$ and
$d_i(\mathbf x)=b_i$ for $i>m$ for the diagonal entries of $J(\mathbf x)$.
Let $\tau\in\R$ denote the spectral parameter.
The left Dirichlet solution
$(p_i(\tau,\mathbf x))_{i\ge0}$ is determined by
\[
p_0(\tau,\mathbf x)=0,
\qquad
p_1(\tau,\mathbf x)=1,
\]
and the recurrence
\[
p_2=\frac{\tau-x_1}{t_1},\qquad
p_{i+1}=\frac{(\tau-d_i(\mathbf x))p_i-t_{i-1}p_{i-1}}{t_i}
\quad\text{for every }i\ge2.
\]
In particular, $p_{i+1}$ depends only on $\tau$ and
the first $i$ diagonal entries.
Define the regular set
\[
U=\left\{(\tau,\mathbf x)\in\R\times\R^m:
p_2(\tau,\mathbf x)\cdots p_m(\tau,\mathbf x)\ne0
\right\}.
\]
For $1\le i\le m$, define the ratios $q_i:U\to\R$ by
\[
q_i(\tau,\mathbf x)
=-t_i\frac{p_{i+1}(\tau,\mathbf x)}{p_i(\tau,\mathbf x)},
\qquad 1\le i\le m.
\]
The identity for $p_2$ gives $q_1=x_1-\tau$. For $2\le i\le m$,
dividing the recurrence by $p_i$ and using
$q_{i-1}=-t_{i-1}p_i/p_{i-1}$ gives
\[
q_i=x_i-\tau-\frac{t_{i-1}^2}{q_{i-1}}.
\]
Thus the boundary ratio can also be written as
\[
q_m(\tau,\mathbf x)
=
x_m-\tau-
\cfrac{t_{m-1}^2}{
x_{m-1}-\tau-
\cfrac{t_{m-2}^2}{
\ddots-\cfrac{t_1^2}{x_1-\tau}
}}.
\]
In particular, $q_m$ is independent of the tail.
We compute the derivatives of $q_m$.
For $i=1$ use $q_1=x_1-\tau$; for $2\le i\le m$, use the fact that
$q_{i-1}$ is independent of $x_i$. Thus, for $1\le i\le m$, $\partial_{x_i}q_i=1$.
Moreover, whenever $1\le j<i\le m$,
\[
\begin{aligned}
\partial_{x_j}q_i
&=\frac{t_{i-1}^2}{q_{i-1}^2}
\partial_{x_j}q_{i-1}=\frac{p_{i-1}^2}{p_i^2}
\partial_{x_j}q_{i-1}.
\end{aligned}
\]
Iterating this identity gives, for every $1\le j\le m$,
\[
\partial_{x_j}q_m
=
\prod_{i=j+1}^m
\frac{p_{i-1}^2}{p_i^2}
=\frac{p_j^2}{p_m^2}.
\]
The continued-fraction formula also gives
\[
q_m(\tau+t,\mathbf x+t(1,\ldots,1))
=q_m(\tau,\mathbf x).
\]
Differentiating at $t=0$, we obtain
\[
\partial_\tau q_m+\sum_{j=1}^m\partial_{x_j}q_m=0.
\]
Consequently,
\begin{equation}\label{spectral:boundary-ratio-derivatives}
\partial_{x_j}q_m=\frac{p_j^2}{p_m^2},
\qquad
-\partial_\tau q_m
=\frac{\sum_{j=1}^m p_j^2}{p_m^2}.
\end{equation}
The first identity will be used to verify the
Monge--Amp\`ere condition, while the second will
relate the boundary ratio to the spectral weights.

\medskip
\noindent\textbf{Step 2: Nondegeneracy and the uniform density bound.}
The dependence of $p_i$ on the preceding diagonal
variables suggests partitioning the first $m=2\ell$
sites into odd and even indices. Set
\[
A=\{1,3,\ldots,m-1\},
\qquad
B=\{2,4,\ldots,m\},
\]
and write
\[
\mathbf x_A=(x_1,x_3,\ldots,x_{m-1}),
\qquad
\mathbf x_B=(x_2,x_4,\ldots,x_m).
\]
Both vectors belong to $\R^\ell$.
At fixed energy $\tau$, consider the Monge--Amp\`ere matrix
\[
\mathcal M(\tau,\mathbf x)=
\begin{pmatrix}
D^2_{AB}q_m&\nabla_Aq_m\\
(\nabla_Bq_m)^T&0
\end{pmatrix},
\]
where $\nabla_A$ and $\nabla_B$ denote the gradients in the ordered
coordinate vectors $\mathbf x_A$ and $\mathbf x_B$, respectively, and
\[
D^2_{AB}q_m
=
\left(
\partial_{x_i}\partial_{x_j}q_m
\right)_{i\in A,\;j\in B}.
\]
We show that $\det\mathcal M(\tau,\mathbf x)\ne0$
throughout $U$.
By \eqref{spectral:boundary-ratio-derivatives},
the row indexed by $i\in A$ is
\[
\left(
\left(
\partial_{x_j}\frac{p_i^2}{p_m^2}
\right)_{j\in B},
\frac{p_i^2}{p_m^2}
\right).
\]
Since $p_1=1$, subtracting $p_i^2$ times the first
row from each odd row $i>1$ gives
\[
\partial_{x_j}\left(\frac{p_i^2}{p_m^2}\right)
-p_i^2\partial_{x_j}\left(\frac1{p_m^2}\right)
=\frac{\partial_{x_j}(p_i^2)}{p_m^2}.
\]
The resulting matrix is
\[
\begin{pmatrix}
(\nabla_B(p_m^{-2}))^T&p_m^{-2}\\
p_m^{-2}\bigl(\partial_{x_{2a}}(p_{2j+1}^2)\bigr)_{1\le j<\ell,\,1\le a\le\ell}&\mathbf 0_{(\ell-1)\times 1}\\
(\nabla_Bq_m)^T&0
\end{pmatrix}.
\]
Expanding along the last column and factoring
$p_m^{-2}$ from the remaining $\ell-1$ odd rows yields
\[
|\det\mathcal M|
=
|p_m|^{-2\ell}
\left|
\det D_B\bigl((p_{2j+1}^2)_{j=1}^{\ell-1},q_m\bigr)
\right|,
\]
where $D_B$ denotes the Jacobian with respect to
$\mathbf x_B$.
We compute this Jacobian using the recurrence.
For $1\le j\le\ell-1$, $p_{2j+1}$ depends only on
$x_1,\ldots,x_{2j}$,
\[
\partial_{x_{2a}}p_{2j+1}=0
\qquad\text{for }1\le j\le\ell-1\text{ and }j<a\le\ell.
\]
When $m=4$, the Jacobian has only the $p_3^2$ and $q_4$ rows.
Thus the Jacobian is lower triangular:
\[
D_B\bigl((p_{2j+1}^2)_{j=1}^{\ell-1},q_m\bigr)
=
(K_{ja})_{1\le j,a\le\ell},\qquad
K_{ja}=\begin{cases}
\partial_{x_{2a}}(p_{2j+1}^2),&1\le j<\ell,\\
\partial_{x_{2a}}q_m,&j=\ell,
\end{cases}
\quad 1\le a\le\ell.
\]
For $1\le j\le\ell-1$, the recurrence at site $2j$ gives
\[
t_{2j}p_{2j+1}
=(\tau-x_{2j})p_{2j}-t_{2j-1}p_{2j-1}.
\]
Neither $p_{2j}$ nor $p_{2j-1}$ depends on $x_{2j}$.
Consequently,
\[
\partial_{x_{2j}}p_{2j+1}
=-\frac{p_{2j}}{t_{2j}},
\qquad
\partial_{x_{2j}}(p_{2j+1}^2)
=-\frac{2p_{2j+1}p_{2j}}{t_{2j}}.
\]
Using also $\partial_{x_m}q_m=1$, we obtain
\begin{equation}\label{spectral:boundary-ratio-determinant}
\begin{aligned}
|\det\mathcal M|
&=
|p_m|^{-2\ell}
\prod_{j=1}^{\ell-1}
\left|
\partial_{x_{2j}}(p_{2j+1}^2)
\right|=
|p_m|^{-2\ell}
\prod_{j=1}^{\ell-1}
\left|
\frac{2p_{2j+1}p_{2j}}{t_{2j}}
\right|
>0.
\end{aligned}
\end{equation}
All the factors are nonzero on $U$. Hence the
odd--even partition satisfies the Monge--Amp\`ere
condition at every regular point.

\medskip
\noindent\emph{A uniform density bound for the boundary ratio.}
Choose a sufficiently small product neighborhood
$I\times Q\Subset U$, where $I,I_1,\ldots,I_m\subset\R$ are
bounded open intervals and $Q=I_1\times\cdots\times I_m$.
Define
\[
\boldsymbol\mu_Q
=\boldsymbol\mu|_Q
=\bigotimes_{i=1}^m(\mu_i|_{I_i}),
\qquad
\gamma_\tau=(q_m(\tau,\cdot))_*\boldsymbol\mu_Q.
\]
On this neighborhood, our goal is a density bound uniform in the energy $\tau$:
\[
\gamma_\tau\ll\Leb,
\qquad
\sup_{\tau\in I}
\left\|\frac{\dd\gamma_\tau}{\dd z}\right\|_{L^\infty(\R)}
<\infty.
\]
Writing $\xi\in\R$ for the Fourier frequency, the definition of the
pushforward gives
\[
\widehat\gamma_\tau(\xi)
=
\int_Q e^{-2\pi i\xi q_m(\tau,\mathbf x)}
\,\dd\boldsymbol\mu(\mathbf x).
\]
Choose $\beta\in C_c^\infty(\R^m)$ equal to one on $\overline Q$
and supported in a slightly larger product box $Q'$
such that
$\overline I\times\overline{Q'}\subset U$.
For the pairing lemma, multiply $q_m$ by a smooth cutoff
supported in $U$ and equal to one near
$\overline I\times\supp\beta$, then extend by zero. This does not
change the kernel or the derivatives on the amplitude support.

Let $\psi_0,\psi\in C_c^\infty(\R)$ be nonnegative cutoffs
satisfying
\[
1=\psi_0(\xi)+\sum_{j\ge1}\psi(2^{-j}\xi),
\qquad \xi\in\R,
\]
where $\psi_0$ is compactly supported and
$\psi$ is supported in $\{1/2\le|\theta|\le2\}$.
\\
For $z\in\R$ and dyadic $R=2^j$ with $j\ge1$, the corresponding
Fourier inversion term is
\[
\begin{aligned}
&\int_\R
e^{2\pi iz\xi}\widehat\gamma_\tau(\xi)
\psi(\xi/R)\,\dd\xi=
\int_Q\int_\R
e^{2\pi i\xi(z-q_m(\tau,\mathbf x))}
\psi(\xi/R)\,\dd\xi\,\dd\boldsymbol\mu(\mathbf x)=
\int K_{R,\tau,z}(\mathbf x)\,
\dd\boldsymbol\mu_Q(\mathbf x),
\end{aligned}
\]
where the change of variables $\xi=R\theta$ gives
\[
K_{R,\tau,z}(\mathbf x)
=
R\int_\R
e^{2\pi iR\theta(z-q_m(\tau,\mathbf x))}
\beta(\mathbf x)\psi(\theta)\,\dd\theta.
\]
Using the odd--even partition, regard $K_{R,\tau,z}$
as the kernel of the operator
\[
T_{R,\tau,z}:L^2(\R^\ell)\longrightarrow L^2(\R^\ell),
\]
defined by
\[
T_{R,\tau,z}h(\mathbf x_A)
=
\int_{\R^\ell}
K_{R,\tau,z}(\mathbf x_A,\mathbf x_B)
h(\mathbf x_B)\,\dd\mathbf x_B.
\]
The phase in Lemma~\ref{bi:fractal-fio-pairing} is
\[
\Phi_{\tau,z}(\mathbf x_A,\mathbf x_B)
=q_m(\tau,\mathbf x)-z.
\]
By \eqref{spectral:boundary-ratio-determinant},
\[
\det
\begin{pmatrix}
D^2_{AB}\Phi_{\tau,z}&\nabla_A\Phi_{\tau,z}\\
(\nabla_B\Phi_{\tau,z})^T&0
\end{pmatrix}
\ne0.
\]
Moreover, since $1\in A$ and $m\in B$,
\[
\partial_{x_1}\Phi_{\tau,z}
=p_m^{-2}>0,
\qquad
\partial_{x_m}\Phi_{\tau,z}=1.
\]
Thus both the Monge--Amp\`ere condition
\eqref{eq:pairing-bordered-condition} and the gradient
condition \eqref{eq:pairing-gradient-condition} hold.
Their lower bounds are uniform on
$\overline I\times\operatorname{supp}\beta\Subset U$.
\\
Define the restricted product measures
\[
\mu_A=\bigotimes_{i\in A}(\mu_i|_{I_i}),
\qquad
\mu_B=\bigotimes_{i\in B}(\mu_i|_{I_i}),
\]
and set
\[
S_A=\sum_{i\in A}s_i,
\qquad
S_B=\sum_{i\in B}s_i.
\]
By the product Frostman bound of Section~\ref{sec:prelim},
$\mu_A$ and $\mu_B$ are $S_A$- and $S_B$-Frostman measures on $\R^\ell$,
respectively; restricting the factors does not increase their Frostman constants.

Applying the pairing estimate
\eqref{bi:fractal-fio-bound} with
$\lambda=R$, $d=\ell$, $\alpha_1=S_A$, and $\alpha_2=S_B$,
we obtain
\[
\begin{aligned}
\left|
\int K_{R,\tau,z}(\mathbf x)\,
\dd\boldsymbol\mu_Q(\mathbf x)
\right|
&=
\left|
\iint K_{R,\tau,z}(\mathbf x_A,\mathbf x_B)
\,\dd\mu_A(\mathbf x_A)\dd\mu_B(\mathbf x_B)
\right|\\
&\lesssim R^{(\ell+1-S_A-S_B)/2}\\
&=R^{(\ell+1-S_m)/2}.
\end{aligned}
\]
If either restricted measure has zero mass, the estimate
is immediate; otherwise, one may apply the pairing
estimate to the normalized restrictions and restore
their masses.

The bound is uniform in $\tau\in I$ and $z$ in a compact
neighborhood of the range of $q_m$ on
$I\times\operatorname{supp}\beta$.
Outside this neighborhood, we have
\[
|z-q_m(\tau,\mathbf x)|\ge c>0
\]
on the support of $\beta$. Using
\[
e^{2\pi iR\theta(z-q_m)}
=
\frac{1}{2\pi iR(z-q_m)}
\partial_\theta
e^{2\pi iR\theta(z-q_m)},
\]
repeated integration by parts in $\theta$ gives $|K_{R,\tau,z}(\mathbf x)|\lesssim_L R^{-L}$
for every $L>0$. Hence
\begin{equation}\label{spectral:boundary-ratio-shell-bound}
\sup_{\tau\in I}\sup_{z\in\R}
\left|
\int_\R e^{2\pi iz\xi}
\widehat\gamma_\tau(\xi)\psi(\xi/R)\,\dd\xi
\right|
\lesssim R^{(\ell+1-S_m)/2}.
\end{equation}
By the dimension hypothesis
\eqref{spectral:intro-infinite-sum}, $S_m>\ell+1$,
and therefore
\[
\sum_{j\ge1}2^{j(\ell+1-S_m)/2}<\infty.
\]
The low-frequency term is also bounded, since
\[
\left|
\int_\R e^{2\pi iz\xi}
\widehat\gamma_\tau(\xi)\psi_0(\xi)\,\dd\xi
\right|
\le\gamma_\tau(\R)\|\psi_0\|_{L^1}.
\]
Consequently, the series
\[
\begin{aligned}
g_\tau(z)
&=
\int_\R e^{2\pi iz\xi}
\widehat\gamma_\tau(\xi)\psi_0(\xi)\,\dd\xi+
\sum_{j\ge1}
\int_\R e^{2\pi iz\xi}
\widehat\gamma_\tau(\xi)\psi(2^{-j}\xi)\,\dd\xi
\end{aligned}
\]
converges uniformly in $\tau\in I$ and $z\in\R$.
Testing against Schwartz functions and using the
dyadic partition of unity gives
$\gamma_\tau=g_\tau(z)\,\dd z$. In particular,
\begin{equation}\label{spectral:boundary-ratio-density}
\sup_{\tau\in I}
\left\|\frac{\dd\gamma_\tau}{\dd z}\right\|_{L^\infty(\R)}
<\infty.
\end{equation}
The bound depends on the regular neighborhood, the
fixed hopping coefficients, and the Frostman constants,
but is independent of the tail.

\medskip
\noindent\textbf{Step 3: A localized estimate for finite restrictions.}
The density estimate so far concerns the boundary ratio at a fixed energy. To obtain a spectral averaging estimate, we relate that ratio to the first-site spectral measure of a finite restriction.

Let $N>m$, and let $J_N(\mathbf x)$ be the Dirichlet
restriction of $J(\mathbf x)$ to sites $1,\ldots,N$.
Write $D_N$ for its principal submatrix on sites
$m+1,\ldots,N$, so that
\[
J_N(\mathbf x)=
\begin{pmatrix}
H_m(\mathbf x)&t_m e_m e^*\\
t_m e e_m^*&D_N
\end{pmatrix},
\]
where $H_m(\mathbf x)$ is the restriction to the
first $m$ sites, $e_m$ is the last standard basis
vector of $\R^m$, and $e$ is the first standard
basis vector of $\R^{N-m}$.

For $\tau\notin\operatorname{spec}(D_N)$, define
\[
a_N(\tau)=t_m^2
\left\langle
(D_N-\tau\mathrm{Id}_{N-m})^{-1}e,e
\right\rangle.
\]
This function encodes the contribution of the tail
to the eigenvalue equation at site $m$. Unlike $a(\tau)$ in
\eqref{spectral:tail-function}, it contains no additive $\tau$ term.

Suppose that $(\tau,\mathbf x)\in I\times Q$ and that
$\tau$ is an eigenvalue of $J_N(\mathbf x)$. Scale
the corresponding real eigenvector $v$ so that its first coordinate is $1$. The recurrence gives
\[
v=(p_1(\tau,\mathbf x),\ldots,p_N(\tau,\mathbf x))^T.
\]
Writing
\[
v_{\mathrm{head}}=(p_1,\ldots,p_m)^T,
\qquad
v_{\mathrm{tail}}=(p_{m+1},\ldots,p_N)^T,
\]
the eigenvalue equation becomes
\[
\begin{aligned}
(H_m(\mathbf x)-\tau\mathrm{Id}_m)v_{\mathrm{head}}
&=-t_m p_{m+1}e_m,\\
(D_N-\tau\mathrm{Id}_{N-m})v_{\mathrm{tail}}
&=-t_mp_me.
\end{aligned}
\]
We first show that $\tau\notin\operatorname{spec}(D_N)$.
Otherwise, choosing a nonzero eigenvector $w$ of
$D_N$ at energy $\tau$ gives
\[
\begin{aligned}
0
&=\langle v_{\mathrm{tail}},
(D_N-\tau\mathrm{Id}_{N-m})w\rangle\\
&=\langle(D_N-\tau\mathrm{Id}_{N-m})
v_{\mathrm{tail}},w\rangle\\
&=-t_mp_m\langle e,w\rangle.
\end{aligned}
\]
Since $t_m\ne0$ and $\langle e,w\rangle\ne0$,
we obtain $p_m=0$, contradicting $(\tau,\mathbf x)\in U$.
Hence
\[
v_{\mathrm{tail}}
=-t_mp_m(D_N-\tau\mathrm{Id}_{N-m})^{-1}e.
\]
Taking the first coordinate gives
\[
p_{m+1}
=-t_mp_m
\left\langle
(D_N-\tau\mathrm{Id}_{N-m})^{-1}e,e
\right\rangle.
\]
Dividing by $p_m$ and multiplying by $-t_m$, we obtain
\[
q_m(\tau,\mathbf x)=a_N(\tau).
\]
Conversely, if this identity holds away from the
tail spectrum, the vector
\[
\begin{pmatrix}
p_1\\
\vdots\\
p_m\\
-t_mp_m(D_N-\tau\mathrm{Id}_{N-m})^{-1}e
\end{pmatrix}
\]
satisfies the eigenvalue equation for $J_N(\mathbf x)$.
Thus, for $(\tau,\mathbf x)\in I\times Q$ with $\tau\notin\operatorname{spec}(D_N)$,
\[
\tau\in\operatorname{spec}(J_N(\mathbf x))
\quad\Longleftrightarrow\quad
q_m(\tau,\mathbf x)=a_N(\tau).
\]
We next compute the spectral weight at such an eigenvalue.
Differentiating $a_N$ gives
\[
\begin{aligned}
a_N'(\tau)
&=t_m^2
\left\langle
(D_N-\tau\mathrm{Id}_{N-m})^{-2}e,e
\right\rangle=t_m^2
\left\|(D_N-\tau\mathrm{Id}_{N-m})^{-1}e\right\|^2.
\end{aligned}
\]
Consequently,
\[
\|v_{\mathrm{tail}}\|^2=p_m^2a_N'(\tau).
\]
Combining this identity with
\eqref{spectral:boundary-ratio-derivatives}, we obtain
\[
\begin{aligned}
\|v\|^2
&=\sum_{j=1}^m p_j^2+\|v_{\mathrm{tail}}\|^2=p_m^2\bigl(-\partial_\tau q_m+a_N'(\tau)\bigr)=-p_m^2\partial_\tau (q_m-a_N).
\end{aligned}
\]
Moreover,
\[
-\partial_\tau (q_m-a_N)
=
\frac{\sum_{j=1}^m p_j^2}{p_m^2}
+a_N'(\tau)>0.
\]
Thus every root of $q_m-a_N$ is simple.

Since the corresponding unit eigenvector is $v/\|v\|$,
its first-site spectral weight is
\[
\begin{aligned}
\sigma_{J_N(\mathbf x),e_1}(\{\tau\})
&=
\left|
\left\langle\frac{v}{\|v\|},e_1\right\rangle
\right|^2=\frac1{\|v\|^2}=\frac1{p_m^2|\partial_\tau (q_m-a_N)|}.
\end{aligned}
\]
In particular, for each fixed $\mathbf x\in Q$ and for any nonnegative continuous function $\varphi$ supported in $I$,
\[
\int_I\varphi(\tau)\,
\dd\sigma_{J_N(\mathbf x),e_1}(\tau)
=
\sum_{\substack{\tau\in I\\q_m(\tau,\mathbf x)=a_N(\tau)}}
\frac{\varphi(\tau)}
{p_m(\tau,\mathbf x)^2
|\partial_\tau (q_m-a_N)(\tau,\mathbf x)|}.
\]

To apply the boundary-ratio density bound to this identity, we approximate a point mass at each root. Let
$0\le\chi\in C_c(I\times Q)$ and
$0\le\varphi\in C_c(I)$. Choose a nonnegative function
$\delta\in C_c^\infty(\R)$ with $\int\delta=1$, and set $\delta_\varepsilon(z)
=\varepsilon^{-1}\delta(z/\varepsilon)$.
At each root $\tau_0$ of $q_m-a_N$, the change of variables $z=q_m(\tau,\mathbf x)-a_N(\tau)$
gives
\[
\begin{aligned}
&\lim_{\varepsilon\downarrow0}
\int_{\tau_0-\eta}^{\tau_0+\eta}
\frac{\varphi(\tau)\chi(\tau,\mathbf x)}
{p_m(\tau,\mathbf x)^2}
\delta_\varepsilon(q_m(\tau,\mathbf x)-a_N(\tau))
\,\dd \tau\\
&\qquad=
\frac{\varphi(\tau_0)\chi(\tau_0,\mathbf x)}
{p_m(\tau_0,\mathbf x)^2
|\partial_\tau (q_m-a_N)(\tau_0,\mathbf x)|},
\end{aligned}
\]
where $\eta>0$ is sufficiently small that the change
of variables is valid and no tail pole is present.
\\
Summing over the roots and applying Fatou's lemma gives
\[
\begin{aligned}
&\int_Q\int_I
\varphi(\tau)\chi(\tau,\mathbf x)\,
\dd\sigma_{J_N(\mathbf x),e_1}(\tau)
\,\dd\boldsymbol\mu(\mathbf x)\\
&\quad\le
\liminf_{\varepsilon\downarrow0}
\int_I\varphi(\tau)
\int_Q
\frac{\chi(\tau,\mathbf x)}{p_m(\tau,\mathbf x)^2}
\delta_\varepsilon(q_m(\tau,\mathbf x)-a_N(\tau))
\,\dd\boldsymbol\mu(\mathbf x)\,\dd \tau.
\end{aligned}
\]
The finitely many tail eigenvalues are omitted from
the energy integral, and, as shown above, none of them
can be an eigenvalue of $J_N(\mathbf x)$ on $I\times Q$.
\\
Since $p_m$ is bounded away from zero on
$\operatorname{supp}\chi$,
\[
\left\|\frac{\chi}{p_m^2}\right\|_\infty<\infty.
\]
Moreover, the definition of $\gamma_\tau$ gives
\[
\begin{aligned}
\int_Q
\delta_\varepsilon(q_m(\tau,\mathbf x)-a_N(\tau))
\,\dd\boldsymbol\mu(\mathbf x)
&=
\int_\R\delta_\varepsilon(z-a_N(\tau))\,
\dd\gamma_\tau(z)\\
&=
\int_\R\delta_\varepsilon(z-a_N(\tau))g_\tau(z)\,\dd z\\
&\le\|g_\tau\|_{L^\infty(\R)}.
\end{aligned}
\]
Therefore, by \eqref{spectral:boundary-ratio-density},
\begin{equation}\label{spectral:localized-finite-averaging}
\begin{aligned}
&\int_Q\int_I
\varphi(\tau)\chi(\tau,\mathbf x)\,
\dd\sigma_{J_N(\mathbf x),e_1}(\tau)
\,\dd\boldsymbol\mu(\mathbf x)\le
\left\|\frac{\chi}{p_m^2}\right\|_\infty
\sup_{\tau\in I}\|g_\tau\|_{L^\infty(\R)}
\int_I\varphi(\tau)\,\dd \tau.
\end{aligned}
\end{equation}
The right-hand side is independent of $N$ and of
the tail entries.

\medskip
\noindent\textbf{Step 4: The half-line limit.}
Let $\Pi_N$ be the orthogonal projection onto
$\operatorname{span}\{e_1,\ldots,e_N\}$.
Extending $J_N(\mathbf x)$ by zero on the orthogonal
complement, we have
\[
J_N(\mathbf x)\oplus0
=\Pi_NJ(\mathbf x)\Pi_N.
\]
For each fixed $\mathbf x\in Q$, set $A_N=J_N(\mathbf x)\oplus0$
and $J=J(\mathbf x)$. Since $\Pi_N\to\mathrm{Id}$ strongly and $J$ is
bounded, $A_N=\Pi_NJ\Pi_N\to J$ strongly. For $z\in\mathbb C\setminus\R$,
the resolvent identity
\[
(A_N-z\mathrm{Id})^{-1}-(J-z\mathrm{Id})^{-1}
=(A_N-z\mathrm{Id})^{-1}(J-A_N)(J-z\mathrm{Id})^{-1}
\]
and $\|(A_N-z\mathrm{Id})^{-1}\|\le|\operatorname{Im}z|^{-1}$ give strong
resolvent convergence. By Reed and Simon~\cite[Theorem~VIII.20(b)]{ReedSimon80},
the first-site spectral measures therefore converge weakly. Apply this to
$h(\tau)=\varphi(\tau)\chi(\tau,\mathbf x)$. Since the first-site
spectral measures have mass one, dominated convergence
gives
\[
\begin{aligned}
&\int_Q\int_I
\varphi(\tau)\chi(\tau,\mathbf x)\,
\dd\sigma_{J(\mathbf x),e_1}(\tau)
\,\dd\boldsymbol\mu(\mathbf x)=
\lim_{N\to\infty}
\int_Q\int_I
\varphi(\tau)\chi(\tau,\mathbf x)\,
\dd\sigma_{J_N(\mathbf x),e_1}(\tau)
\,\dd\boldsymbol\mu(\mathbf x).
\end{aligned}
\]
By \eqref{spectral:localized-finite-averaging}, we obtain
\begin{equation}\label{spectral:localized-halfline-averaging}
\begin{aligned}
&\int_Q\int_I
\varphi(\tau)\chi(\tau,\mathbf x)\,
\dd\sigma_{J(\mathbf x),e_1}(\tau)
\,\dd\boldsymbol\mu(\mathbf x)\le C_{I,Q,\chi}
\int_I\varphi(\tau)\,\dd \tau,
\end{aligned}
\end{equation}
where
\[
C_{I,Q,\chi}
=
\left\|\frac{\chi}{p_m^2}\right\|_\infty
\sup_{\tau\in I}\|g_\tau\|_{L^\infty(\R)}
\]
is independent of the tail.
In particular, the localized energy measure
\[
\dd\omega_\chi(\tau)
=
\int_Q\chi(\tau,\mathbf x)
\,\dd\sigma_{J(\mathbf x),e_1}(\tau)
\,\dd\boldsymbol\mu(\mathbf x)
\]
satisfies $\omega_\chi\le C_{I,Q,\chi}\Leb$.

\medskip
\noindent\textbf{Step 5: Exceptional energies and global averaging.}
The estimate \eqref{spectral:localized-halfline-averaging} applies on the regular set $U$. It remains to show that its complement carries no first-site spectral mass after averaging.
Fix the tail and write
\[
\mathbf x=(x_1,\widehat{\mathbf x}),
\qquad
\widehat{\mathbf x}=(x_2,\ldots,x_m).
\]
For each fixed $\mathbf x$, define
\[
\mathcal B(\mathbf x)
=
\{\tau\in\R:
p_2(\tau,\mathbf x)\cdots p_m(\tau,\mathbf x)=0\}.
\]
The recurrence gives
\[
p_i(\tau,\mathbf x)
=
\frac{\tau^{i-1}}{t_1\cdots t_{i-1}}
+\text{terms of lower degree in }\tau.
\]
Hence each $p_i$, $2\le i\le m$, is a nonzero
polynomial of degree $i-1$. In particular,
\[
|\mathcal B(\mathbf x)|
\le\sum_{i=2}^m(i-1)
=\frac{m(m-1)}2.
\]
Thus only atoms of $\sigma_{J(\mathbf x),e_1}$
can contribute to $\mathcal B(\mathbf x)$.
Fix $\widehat{\mathbf x}$. Suppose that
$\tau\in\mathcal B(\mathbf x)$ is an eigenvalue of
$J(\mathbf x)$ with eigenvector $u$. Since every
eigenvector is a multiple of the left Dirichlet solution,
\[
u_i=\langle u,e_1\rangle p_i(\tau,\mathbf x).
\]
Therefore, $u_i=0$ for some $2\le i\le m$.
\\
Let $D_i$ denote the right Dirichlet restriction
of $J(\mathbf x)$ to sites $i+1,i+2,\ldots$.
Writing $u_{>i}$ for the restriction of $u$ to these
sites, the eigenvalue equation gives
\[
(D_i-\tau\mathrm{Id})u_{>i}
=-t_i u_i e_{i+1}=0,
\]
where $e_{i+1}$ is the first coordinate vector in
the tail space. Moreover, $u_{>i}\ne0$, since two
consecutive zero coordinates would force the whole
eigenvector to vanish by the recurrence. Hence
\[
\tau\in\sigma_{\mathrm p}(D_i),
\]
where $\sigma_{\mathrm p}$ denotes the point spectrum.
Since the operators $D_i$ are independent of $x_1$,
all possible exceptional eigenvalues, as $x_1$ varies,
belong to the fixed countable set
\[
\mathcal E(\widehat{\mathbf x})
=\bigcup_{i=2}^m\sigma_{\mathrm p}(D_i).
\]
Each point spectrum is countable by self-adjointness
and separability of the underlying Hilbert space.

We claim that a fixed energy $\tau$ can be an eigenvalue
of $J(x_1,\widehat{\mathbf x})$ for at most one
value of $x_1$. Indeed, suppose
\[
J(x_1,\widehat{\mathbf x})u=\tau u,
\qquad
J(x_1',\widehat{\mathbf x})v=\tau v.
\]
By self-adjointness,
\[
\begin{aligned}
0
&=\langle J(x_1,\widehat{\mathbf x})u,v\rangle
-\langle u,J(x_1',\widehat{\mathbf x})v\rangle=
\left\langle
\bigl(J(x_1,\widehat{\mathbf x})
-J(x_1',\widehat{\mathbf x})\bigr)u,v
\right\rangle.
\end{aligned}
\]
Since
\[
J(x_1,\widehat{\mathbf x})
-J(x_1',\widehat{\mathbf x})
=(x_1-x_1')e_1e_1^*,
\]
we obtain
\[
(x_1-x_1')
\langle u,e_1\rangle
\langle e_1,v\rangle=0.
\]
The first coordinate of every nonzero Jacobi
eigenvector is nonzero. Hence $x_1=x_1'$.
By nonatomicity of $\mu_1$, it follows that
\[
\int_\R
\sigma_{J(x_1,\widehat{\mathbf x}),e_1}(\{\tau\})
\,\dd\mu_1(x_1)=0
\]
for every fixed $\tau\in\mathcal E(\widehat{\mathbf x})$. Consequently,
\[
\begin{aligned}
&\int_\R
\sigma_{J(\mathbf x),e_1}(\mathcal B(\mathbf x))
\,\dd\mu_1(x_1)\le
\sum_{\tau\in\mathcal E(\widehat{\mathbf x})}
\int_\R
\sigma_{J(\mathbf x),e_1}(\{\tau\})
\,\dd\mu_1(x_1)=0.
\end{aligned}
\]
Since $\{(\tau,\mathbf x):\tau\in\mathcal B(\mathbf x)\}=U^c$
is Borel, Tonelli's theorem gives
\[
\int_{\R^m}
\sigma_{J(\mathbf x),e_1}(\mathcal B(\mathbf x))
\,\dd\boldsymbol\mu(\mathbf x)=0.
\]

\medskip
\noindent\emph{Combining the local estimates.}
Define the joint measure $\mathfrak m$ on
$\R\times\R^m$ by
\[
\dd\mathfrak m(\tau,\mathbf x)
=
\dd\sigma_{J(\mathbf x),e_1}(\tau)
\,\dd\boldsymbol\mu(\mathbf x).
\]
The exceptional set contributes no mass, as the preceding argument gives
\[
\mathfrak m(U^c)
=
\int_{\R^m}
\sigma_{J(\mathbf x),e_1}(\mathcal B(\mathbf x))
\,\dd\boldsymbol\mu(\mathbf x)=0.
\]
Choose a countable cover of $U$ by regular product
neighborhoods $I_j\times Q_j\Subset U$ and a
nonnegative smooth partition of unity $(\chi_j)_j$
with each $\chi_j$ compactly supported in one of
these neighborhoods. Thus
\[
\sum_j\chi_j(\tau,\mathbf x)=1
\qquad\text{on }U.
\]
Let
\[
\pi_\tau:\R\times\R^m\longrightarrow\R,
\qquad
\pi_\tau(\tau,\mathbf x)=\tau.
\]
By \eqref{spectral:localized-halfline-averaging}, $(\pi_\tau)_*(\chi_j\mathfrak m)\ll\Leb$
for every $j$.
Since $\mathfrak m(U^c)=0$, we have
\[
\mathfrak m=\sum_j\chi_j\mathfrak m.
\]
Lemma~\ref{lem:local-to-global-ac} therefore gives
\[
(\pi_\tau)_*\mathfrak m
=
\sum_j(\pi_\tau)_*(\chi_j\mathfrak m)
\ll\Leb.
\]
Equivalently,
\[
\int_{\R^m}
\sigma_{J(\mathbf x),e_1}
\,\dd\boldsymbol\mu(\mathbf x)
\ll\Leb.
\]
This holds for every fixed bounded tail. If the tail
is random and independent of the first $m$ diagonal
entries, write $\vartheta$ for its law and
$J_{\mathbf b}(\mathbf x)$ for the corresponding
operator. Then, for every Lebesgue-null Borel set
$Z\subset\R$,
\[
\begin{aligned}
\mathbb E\,\sigma_{J(\mathbf X),e_1}(Z)
&=
\int
\left(
\int_{\R^m}
\sigma_{J_{\mathbf b}(\mathbf x),e_1}(Z)
\,\dd\boldsymbol\mu(\mathbf x)
\right)\dd\vartheta(\mathbf b)=0.
\end{aligned}
\]
Thus $\mathbb E\,\sigma_{J(\mathbf X),e_1}\ll\Leb$.

\medskip
\noindent\textbf{Step 6: Cyclicity and the conclusion.}
Cyclicity lets us pass from the first-site spectral measure to that of an arbitrary vector $u\in\ell^2(\mathbb N)$.

Fix a Lebesgue-null Borel set $Z\subset\R$.
By the preceding absolute continuity,
\[
\begin{aligned}
0
&=\mathbb E\,\sigma_{J(\mathbf X),e_1}(Z)=\mathbb E\,
\langle\mathbf1_Z(J(\mathbf X))e_1,e_1\rangle=\mathbb E\,
\|\mathbf1_Z(J(\mathbf X))e_1\|^2.
\end{aligned}
\]
Hence
\[
\mathbf1_Z(J(\mathbf X))e_1=0
\qquad\text{almost surely}.
\]
For every realization with bounded diagonal, $e_1$ is cyclic by the half-line spectral
representation~\cite[Theorem~2.12]{Teschl00}: it is sent to the constant
function $1$, and the polynomials are dense for the compactly supported
first-site spectral measure.
Since $\mathbf1_Z(J)$ commutes with $J$,
\[
\mathbf1_Z(J)J^je_1
=J^j\mathbf1_Z(J)e_1=0
\]
for every $j\ge0$. By cyclicity,
\[
\mathbf1_Z(J(\mathbf X))=0
\qquad\text{almost surely}
\]
for the fixed null set $Z$. Thus, for every
$u\in\ell^2(\mathbb N)$,
\[
\begin{aligned}
\mathbb E\,\sigma_{J(\mathbf X),u}(Z)
&=\mathbb E\,
\langle\mathbf1_Z(J(\mathbf X))u,u\rangle=0.
\end{aligned}
\]
As $Z$ was an arbitrary Lebesgue-null Borel set,
we conclude that $\mathbb E\,\sigma_{J(\mathbf X),u}\ll\Leb$
for every $u\in\ell^2(\mathbb N)$.
\\
Finally, if the first $m$ site laws are all
$s$-Frostman with $s>1/2$, the dimension condition is
\[
ms>\frac m2+1
\quad\Longleftrightarrow\quad
m>\frac{2}{2s-1}.
\]
Choosing an even integer $m$ satisfying this inequality
proves the final assertion.
\end{proof}

\section*{Acknowledgement}
C.-Y. Shen was partially supported by NSTC under Grant No.~111-2115-M-002-010-MY5.

\end{document}